\documentclass[reqno,10pt,centertags]{amsart}
\usepackage{amsmath,amsthm,amscd,amssymb,latexsym,upref,stmaryrd}

\usepackage{hyperref}
\newcommand*{\mailto}[1]{\href{mailto:#1}{\nolinkurl{#1}}}

\newcommand{\bbN}{{\mathbb{N}}}
\newcommand{\bbR}{{\mathbb{R}}}

\newcommand{\bbC}{{\mathbb{C}}}

\newcommand{\cA}{{\mathcal A}}
\newcommand{\cB}{{\mathcal B}}

\newcommand{\cD}{{\mathcal D}}

\newcommand{\cH}{{\mathcal H}}

\newcommand{\cK}{{\mathcal K}}

\newcommand{\cM}{{\mathcal M}}

\newcommand{\no}{\notag}
\newcommand{\dott}{\,\cdot\,}
\newcommand{\lb}{\label}
\newcommand{\bi}{\bibitem}

\newcommand{\ol}{\overline}

\newcommand{\wti}{\widetilde}
\newcommand{\hatt}{\widehat}

\newcommand{\f}{\frac}
\newcommand{\supp}{\operatorname{supp}}

\renewcommand{\Re}{\mathop\mathrm{Re}}

\renewcommand{\le}{\leqslant}

\DeclareMathOperator{\dom}{dom}
\DeclareMathOperator{\ran}{ran} 
\DeclareMathOperator{\esssup}{ess\,sup}
\DeclareMathOperator{\essinf}{ess\,inf}

\DeclareMathOperator*{\slim}{s-lim}
\DeclareMathOperator*{\wlim}{w-lim}
\mathchardef\mhyphen="2D

\newcommand{\si}{\sigma}

\newcommand{\la}{\lambda}

\allowdisplaybreaks 
\numberwithin{equation}{section}

\newtheorem{theorem}{Theorem}[section]

\newtheorem{lemma}[theorem]{Lemma}
\newtheorem{corollary}[theorem]{Corollary}
\newtheorem{hypothesis}[theorem]{Hypothesis}

\newtheorem{example}[theorem]{Example}
\theoremstyle{remark}
\newtheorem{remark}[theorem]{Remark}

\begin{document}

\title[Abstract Wave Equations and Dirac-Type Operators]{Abstract Wave 
Equations and Associated Dirac-Type Operators}

\author[F.\ Gesztesy]{Fritz Gesztesy}
\address{Department of Mathematics,
University of Missouri,
Columbia, MO 65211, USA}
\email{\mailto{gesztesyf@missouri.edu}}
\urladdr{\url{http://www.math.missouri.edu/personnel/faculty/gesztesyf.html}}

\author[J.\ A.\ Goldstein]{Jerome A.\ Goldstein}
\address{Department of Mathematical Sciences, University of Memphis, Memphis, Tennessee 38152, USA}
\email{\mailto{jgoldste@memphis.edu}}
\urladdr{\url{http://www.msci.memphis.edu/faculty/goldsteinj.html}}

\author[H.\ Holden]{Helge Holden}
\address{Department of Mathematical Sciences,
Norwegian University of
Science and Technology, NO--7491 Trondheim, Norway}
\email{\mailto{holden@math.ntnu.no}}
\urladdr{\url{http://www.math.ntnu.no/~holden/}}

\author[G.\ Teschl]{Gerald Teschl}
\address{Faculty of Mathematics\\ University of Vienna\\ 
Nordbergstrasse 15\\ 1090 Wien\\ Austria\\ and International
Erwin Schr\"odinger
Institute for Mathematical Physics\\ Boltzmanngasse 9\\ 1090 Wien\\ Austria}
\email{\mailto{Gerald.Teschl@univie.ac.at}}
\urladdr{\url{http://www.mat.univie.ac.at/~gerald/}}

\thanks{{\it Ann. Mat. Pura Appl.} {\bf 191}, 631--676 (2012)}
\thanks{{\it Research supported by in part by the Research Council of Norway 
and the Austrian Science Fund $($FWF$\,)$ under Grant No.\ Y330}}

\date{\today}

%%%%%%%%%%%%%%%%%%%%%%%%%%%%%%%%%%%%%%%%
%%%%%%%%%%%%%%%%%%%%%%%%%%%%%%%%%%%%%%%%
\begin{abstract}
We discuss the unitary equivalence of generators $G_{A,R}$ associated with abstract damped wave equations of the type $\ddot{u} + R \dot{u} + A^*A u = 0$ 
in some Hilbert space $\mathcal{H}_1$ and certain non-self-adjoint Dirac-type 
operators $Q_{A,R}$ (away from the nullspace of the latter) in $\mathcal{H}_1 \oplus \mathcal{H}_2$. The operator $Q_{A,R}$ represents a non-self-adjoint perturbation of 
a supersymmetric self-adjoint Dirac-type operator. Special emphasis is devoted to the case where $0$ belongs to the continuous spectrum of $A^*A$. 

In addition to the unitary equivalence results concerning $G_{A,R}$ and 
$Q_{A,R}$, we provide a detailed study of the domain of the generator 
$G_{A,R}$, consider spectral properties of the underlying quadratic operator 
pencil $M(z) = |A|^2 - iz R - z^2 I_{\mathcal{H}_1}$, $z\in\mathbb{C}$, derive a family of 
conserved quantities for abstract wave equations in the absence of damping, 
and prove equipartition of energy for supersymmetric self-adjoint 
Dirac-type operators. 

The special example where $R$ represents an appropriate function of $|A|$ is treated in depth and the semigroup growth bound for this example is explicitly 
computed and shown to coincide with the corresponding spectral bound for the underlying generator and also with that of the corresponding Dirac-type 
operator. 

The cases of undamped ($R=0$) and damped ($R \neq 0$) abstract wave equations as well as the cases $A^* A \geq \varepsilon I_{\mathcal{H}_1}$ for some 
$\varepsilon > 0$ and $0 \in \sigma (A^* A)$ (but $0$ not an eigenvalue of 
$A^*A$) are separately studied in detail. 
\end{abstract}
%%%%%%%%%%%%%%%%%%%%%%%%%%%%%%%%%%%%%%%%
%%%%%%%%%%%%%%%%%%%%%%%%%%%%%%%%%%%%%%%%

\maketitle

%%%%%%%%%%%%%%%%%%%%%%%%%%%%%%%%%%%%%%%%
%%%%%%%%%%%%%%%%%%%%%%%%%%%%%%%%%%%%%%%%
\section{Introduction}
\label{s1}
%%%%%%%%%%%%%%%%%%%%%%%%%%%%%%%%%%%%%%%%
%%%%%%%%%%%%%%%%%%%%%%%%%%%%%%%%%%%%%%%%

We are interested in an abstract version of the damped wave equation of 
the form
\begin{equation}
\ddot u (t) + R \dot u (t) + A^*A u (t) = 0, \quad u(0) = f_0, \;  
\dot u (0) = f_1, \quad  t \geq 0,    \lb{1.1}
\end{equation}
where $A$ is a densely defined closed operator in a separable Hilbert space $\cH$, $f_j \in \cH$, 
$j=0,1$, are chosen appropriately, $R$ is a certain perturbation of $A^*A$ to be specified in more detail in Section \ref{s3}, and we used the abbreviations $\dot u = (d/dt) u$, $\ddot u = (d^2/dt^2) u$. (In the main body of this paper we will employ a two-Hilbert space approach where $A$ maps its domain, a dense subspace of the Hilbert space $\cH_1$ into a Hilbert space $\cH_2$.)
 
Traditionally, one rewrites \eqref{1.1} in the familiar first-order form
\begin{equation}
\f{d}{dt} \begin{pmatrix} u \\ \dot u \end{pmatrix} =\begin{pmatrix} 0 & I_{\cH}  
\\ - A^* A & -R \end{pmatrix}
\begin{pmatrix} u \\ \dot u \end{pmatrix}, \quad 
\begin{pmatrix} u(0) \\ \dot u (0) \end{pmatrix} = \begin{pmatrix} f_0 \\ f_1 \end{pmatrix}, 
\quad t \geq 0.   \lb{1.2}
\end{equation}

Our principal result centers around a unitary equivalence between an appropriate operator realization of the formal generator $G_{A,R}$ of \eqref{1.2}, 
\begin{equation}
G_{A,R} = \begin{pmatrix} 0 & I_{\cH} \\ - A^* A & -R \end{pmatrix},    \lb{1.3}
\end{equation}
in an associated energy space $\cH_A \oplus \cH$ to be determined in 
Section \ref{s2}, and the operator 
\begin{equation}
Q_{A,R} (I_{\cH} \oplus [I_{\cH} - P_{\ker(A^*)}])
 = \begin{pmatrix} -i \, R & A^*[I_{\cH} - P_{\ker(A^*)}] \\ A & 0 \end{pmatrix},   
 \lb{1.4}
\end{equation}
with $Q_{A,R}$ a perturbed supersymmetric Dirac-type operator in $\cH \oplus \cH$,
\begin{equation}
Q_{A,R} = \begin{pmatrix} - i \, R & A^* \\  A & 0 \end{pmatrix}, \quad 
\dom (Q_{A,R}) = \dom(A) \oplus \dom(A^*) \subseteq \cH \oplus \cH.   
\lb{1.5} 
\end{equation} 

More precisely, we will first establish the unitary equivalence between the self-adjoint operators 
$i \, G_{A,0}$ in $\cH_A \oplus \cH$ and $Q_{A,0}$ in $\cH_1 \oplus \cH_2$ 
and then treat the damping terms 
$\left(\begin{smallmatrix} 0 & 0 \\ 0 & -R \end{smallmatrix}\right)$ and 
$\left(\begin{smallmatrix} -i \, R & 0 \\ 0 & 0 \end{smallmatrix}\right)$ 
perturbatively, keeping the same unitary equivalence between $i \, G_{A,R}$ 
and $Q_{A,R}$ for $R \neq 0$. 

Particular attention is devoted to domain properties of the generator $G_{A,R}$. Moreover, we carefully distinguish the cases of undamped ($R=0$) and damped 
($R\neq 0$) abstract wave equations, and the cases where 
$A^* A \geq \varepsilon I_{\cH}$ for some $\varepsilon > 0$ and the far more subtle situation where $ 0 \in \sigma (A^* A)$ (but $0$ is not an eigenvalue of 
$A^* A$). 

More precisely, the case where $A^* A \geq \varepsilon I_{\cH}$ for some 
$\varepsilon > 0$ and no damping, that is, the situation $R=0$, is treated in 
Section \ref{s2}. The unitary equivalence 
of the generator $G_{A,0}$ and the supersymmetric self-adjoint Dirac-type 
operator $Q_{A,0}$ (away from its nullspace) is the centerpiece of this section. Section \ref{s2} concludes with a discussion of the special case where $A$ is 
replaced by the self-adjoint operator $|A|$. Section \ref{s3} then considers the 
more general case where $0 \in \sigma (A^* A)$ (but $0$ is not an eigenvalue 
of $A$). After establishing the appropriate extension of the unitary equivalence 
of the generator $G_{A,0}$ and the supersymmetric self-adjoint Dirac-type 
operator $Q_{A,0}$ (away from its nullspace) in this case, we provide a detailed 
study of the domain of the generator $G_{A,0}$. Abstract damped linear wave 
equations, assuming $A^* A \geq \varepsilon I_{\cH}$ for some 
$\varepsilon > 0$, are studied in Section \ref{s4}. In this section we also 
compute the resolvent of $Q_{|A|,R}$ in terms of the quadratic operator pencil 
$M(z) = |A|^2 - iz R - z^2 I_{\cH_1}$, $\dom(M(z)) = \dom\big(|A|^2\big)$, 
$z \in\bbC$, and relate the spectrum of $Q_{|A|,R}$ with that of the pencil 
$M(\cdot)$. This section once more derives the unitary equivalence results 
between $Q_{|A|,R}$ and $G_{A,R}$ and similarly, between $Q_{A,R}$ 
(away from its nullspace) and $G_{A,R}$. We also briefly revisit classical 
solutions for the abstract first-order and second-order Cauchy problems. 
Section \ref{s4} concludes with a detailed discussion of the example where the 
damping term $R= 2 F(|A|) \geq 0$ is an appropriate function of $|A|$. 
Employing the spectral theorem for the self-adjoint operator $|A|$, the 
semigroup growth bound for $e^{G_{A,2F(|A|)}t}$, $t\geq 0$, is explicitly 
computed and shown to coincide with the corresponding spectral bound for 
the underlying generator $G_{A,2F(|A|)}$ and hence also with that of 
$-i Q_{|A|,2F(|A|)}$. The most general case of abstract damped wave 
equations where $0 \in \sigma (A^* A)$ (but $0$ is not an eigenvalue of 
$A$) is considered in Section \ref{s5}. Again we compute the resolvent of 
$Q_{|A|,R}$ in terms of the quadratic operator pencil 
$M(z) = |A|^2 - iz R - z^2 I_{\cH_1}$, $\dom(M(z)) = \dom\big(|A|^2\big)$, 
$z \in\bbC$, and relate the spectrum of $Q_{|A|,R}$ with that of the pencil 
$M(\cdot)$. In addition, we once more derive the unitary equivalence results 
between $Q_{|A|,R}$ and $G_{A,R}$ and similarly, between $Q_{A,R}$ 
(away from its nullspace) and $G_{A,R}$. Section \ref{s5} concludes with a derivation of a family conserved quantities for the abstract wave equation in 
the absence of damping. In Section \ref{s6} we prove equipartition of energy 
for the supersymmetric self-adjoint Dirac-type operator $Q = Q_{A,0}$. 
Appendix \ref{sA} summarizes well-known results on supersymmetric 
Dirac-type operators used throughout the bulk of this manuscript and 
Appendix \ref{sB} studies adjoints and closures of products of linear operators.
 
Concluding this introduction, we briefly summarize some of the notation used in this paper. Let $\cH$ be a separable complex Hilbert space, $(\cdot,\cdot)_{\cH}$ the scalar product in $\cH$ (linear in the second factor), and $I_{\cH}$ the identity operator in $\cH$. Next, let $T$ be a linear operator mapping (a subspace of) a
Hilbert space into another, with $\dom(T)$, $\ran(T)$, and $\ker(T)$ denoting the
domain, range, and kernel (i.e., null space) of $T$, respectively. The closure of a closable 
operator $S$ in $\cH$ is denoted by $\ol S$. The spectrum, essential spectrum, point spectrum, discrete spectrum, and resolvent set of a closed linear operator 
in $\cH$ will be denoted by 
$\sigma(\cdot)$, $\sigma_{\rm ess}(\cdot)$, $\sigma_{\rm p}(\cdot)$, 
$\sigma_{\rm d}(\cdot)$, and 
$\rho(\cdot)$, respectively. The Banach space of bounded linear operators 
in $\cH$ is 
denoted by $\cB(\cH)$; the analogous notation $\cB(\cH_1,\cH_2)$ will be 
used for bounded operators between two Hilbert spaces $\cH_1$ and $\cH_2$. 
The norm in $\cH_1 \oplus \cH_2$ is defined as usual by 
$\|f\|_{\cH_1 \oplus \cH_2} = \big[\|f_1\|_{\cH_1}^2 + \|f_2\|_{\cH_2}^2\big]^{1/2}$ 
for $f =(f_1 \; f_2)^\top \in \cH_1 \oplus \cH_2$. The symbols $\slim$ 
(resp., $\wlim$)  denote the strong (resp., weak) limits either in the context of Hilbert space vectors or in the context of bounded operators between two 
Hilbert spaces. Finally, $P_{\cM}$ denotes the orthogonal projection onto 
a closed, linear subspace $\cM$ of $\cH$.

%%%%%%%%%%%%%%%%%%%%%%%%%%%%%%%%%%%%%%
%%%%%%%%%%%%%%%%%%%%%%%%%%%%%%%%%%%%%%
\section{Abstract Linear Wave Equations in the Absence of Damping. \\ 
The Case $A^* A \geq \varepsilon I_{\cH}$ for some $\varepsilon >0$}
\label{s2}
%%%%%%%%%%%%%%%%%%%%%%%%%%%%%%%%%%%%%%
%%%%%%%%%%%%%%%%%%%%%%%%%%%%%%%%%%%%%%

In this section we consider self-adjoint realizations of $i\, G_{A,0}$ modeling abstract linear wave equations in the absence of damping and study their unitary equivalence to self-adjoint supersymmetric Dirac-type operators.

To set the stage we first introduce the following assumptions used throughout this section.

%%%%%%%%%%%
\begin{hypothesis} \lb{h2.1}
Let $\cH_j$, $j=1,2$, be complex separable Hilbert spaces. 
Assume that $A: \dom(A) \subseteq \cH_1 \to \cH_2$ is a densely 
defined, closed, linear operator such that 
\begin{equation}
A^* A \geq \varepsilon I_{\cH_1}     \lb{2.2}
\end{equation} 
for some $\varepsilon >0$. 
\end{hypothesis}
%%%%%%%%%%% 

To illustrate the implications of Hypothesis \ref{h2.1}, we briefly digress a bit. 
Let $T : \dom(T)\subseteq \cH_1\to \cH_2$ be a densely 
defined, closed, linear operator. We recall the definition of the self-adjoint operator 
$|T| = (T^* T)^{1/2}$ in $\cH_1$ and note that 
\begin{align}
& \dom(|T|) = \dom(T), \quad \ker(|T|) = \ker(T^*T) = \ker(T), \quad 
\ran(|T|) = \ran(T^*),   \lb{2.3} \\
& \||T| f\|_{\cH_1} = \|T f\|_{\cH_2}, \; f \in \dom(T).     \lb{2.4}
\end{align}
The latter fact immediately follows from the polar decomposition of $T$ 
(cf.\ \eqref{A.4a}--\eqref{A.4e}).

Thus, Hypothesis \ref{h2.1} is equivalent to 
\begin{equation}
|A| \geq \varepsilon^{1/2} I_{\cH_1},     \lb{2.7}
\end{equation}
and hence equivalent to 
\begin{equation}
|A|^{-1} \in \cB(\cH_1), \, \text{ or equivalently, to } \, 0 \in \rho(|A|).    \lb{2.8}
\end{equation} 
In particular, it implies that
\begin{equation}
\ker(A) = \{0\}.  \lb{2.8a}
\end{equation}

Since $A$ is closed and $|A| \geq \varepsilon^{1/2} I_{\cH_1}$, the norm $\| \dott \|_{A}$ on the 
subspace $\dom(A)$ of $\cH_1$ defined by  
\begin{equation}
\| f \|_{A} = \|A f \|_{\cH_2}, \quad f\in\dom(A),    \lb{2.9}
\end{equation}
and the graph norm $\interleave \cdot \interleave_{A}$ on $\dom(A)$ defined by 
\begin{align}
\begin{split}
\interleave f \interleave_{A} = \|Af\|_{\cH_2} + \|f\|_{\cH_1} \Big(  \text{or alternatively, by }
\big[\|Af\|_{\cH_2}^2 + \|f\|_{\cH_1}^2\big]^{1/2} \Big),& \\ 
f \in \dom(A),&     \lb{2.10}
\end{split}
\end{align}
are equivalent norms on $\dom(A)$. In particular, one verifies that 
\begin{equation}
\f{\varepsilon}{1+\varepsilon} \big[\|Af\|_{\cH_2} + \|f\|_{\cH_1}\big] \leq \| f \|_{A} 
\leq \big[\|Af\|_{\cH_2} + \|f\|_{\cH_1}\big], 
\quad f \in \dom(A).     \lb{2.11}
\end{equation}
Associated with the norm $\| \dott \|_{A}$ we also introduce the corresponding scalar product 
$(\dott,\dott)_A$ on $\dom(A)$ by
\begin{equation}
(f,g)_A = (Af,Ag)_{\cH_2}, \quad f, g \in \dom(A).    \lb{2.12}
\end{equation}
Consequently, equipping the linear space $\dom(A)$ with the scalar product $(\dott,\dott)_A$, one arrives at a Hilbert space denoted by $\cH_A$,
\begin{equation}
\cH_A = \big(\dom(A); (\dott,\dott)_A\big) \subseteq \cH_1.     \lb{2.13}
\end{equation} 

We emphasize that while Hypothesis \ref{h2.1} implies $|A|^{-1}\in\cB(\cH_1)$, it does {\it not} 
imply that $A$ is boundedly invertible on all of $\cH_1$ (mapping into $\cH_2$), as one can see from the following typical example.

%%%%%%%%%%%%%
\begin{example} \lb{e2.2}
Consider the operator $B$ in the Hilbert space $L^2([0,1]; dx)$ defined by
\begin{align}
\begin{split}
& B f = -i f', \quad 
 f\in\dom(B) = \big\{g \in L^2([0,1]; dx) \,\big|\, g \in AC([0,1]);   \\
 & \hspace*{4.7cm} g(0)=g(1)=0; \, g' \in L^2([0,1]; dx)\big\}.   \lb{2.15}
 \end{split} 
\end{align}
Then $B$ is symmetric, its adjoint is given by 
\begin{align}
\begin{split}
& B^* f = -i f',    \lb{2.16} \\
& f\in\dom(B^*) = \big\{g \in L^2([0,1]; dx) \,\big|\, g \in AC([0,1]); \, g' \in L^2([0,1]; dx)\big\},  
 \end{split} 
\end{align}
and the deficiency indices $n_\pm (B)$ of $B$ are given by
\begin{equation}
n_\pm (B) = 1.     \lb{2.17}
\end{equation}
Consequently, 
\begin{equation}
\sigma(B) = \bbC,     \lb{2.18}
\end{equation}
in particular, $B$ is not boundedly invertible on $L^2([0,1]; dx)$. On the other hand,  
\begin{align}
\begin{split} 
& B^* B f = - f'', \quad \dom(B^* B) = \big\{g \in L^2([0,1]; dx) \,\big|\, 
g, g' \in AC([0,1]);   \\
 & \hspace*{5.1cm} g(0)=g(1)=0; \, g'' \in L^2([0,1]; dx)\big\},    \lb{2.19}
\end{split}
\end{align}
$($implying also $g' \in L^2([0,1]; dx)$$)$ and hence 
\begin{equation}
|B| \geq \pi I_{L^2([0,1]; dx)}, \quad  |B|^{-1} \in \cB\big(L^2([0,1]; dx)\big).    \lb{2.20}
\end{equation} 
In fact, one has $|B|^{-1} \in \cB_p \big(L^2([0,1]; dx)\big)$ for all $p>1$. $($Here 
$\cB_p\big(L^2([0,1]; dx)\big)$, $p>0$, denote the $\ell^p(\bbN)$-based trace 
ideals of $L^2([0,1]; dx)$.$)$
\end{example}
%%%%%%%%%%%%%

In this context we note that by \eqref{2.3} and \eqref{2.4}, one has of course
\begin{equation}
\cH_A = \big(\dom(A); (\dott,\dott)_A\big) = 
\cH_{|A|} = \big(\dom(|A|); (\dott,\dott)_{|A|}\big) \subseteq \cH_1,     \lb{2.20a}
\end{equation} 
which is of some significance since under Hypothesis \ref{h2.1} we always have 
\begin{equation}
0 < \varepsilon^{-1/2} I_{\cH_1} \leq |A|^{-1} \in \cB(\cH_1),     \lb{2.20b}
\end{equation}
while in general (cf.\ Example \ref{e2.2}), $A$ is not boundedly invertible on all 
of $\cH_1$ (mapping into $\cH_2$).

The following result is well-known, but for convenience, we provide its short proof.

%%%%%%%%%%%%%
\begin{lemma} \lb{l2.3}
Assume Hypothesis \ref{h2.1}. Then $\ran(A)$ is a closed linear subspace of $\cH_2$.
\end{lemma}
%%%%%%%%%%%%%
\begin{proof}
Let $\{g_n\}_{n\in\bbN} \subset \ran(A)$ be a Cauchy sequence, that is, $g_n = A f_n$, $n\in\bbN$, 
for some sequence $\{f_n\}_{n\in\bbN} \subset \dom(A)$, and hence suppose that 
$\lim_{n\to\infty}\|g_n - g\|_{\cH_2}=0$ for some $g\in\cH_2$. Since by \eqref{2.11},
\begin{align}
\|g_n - g_m\|_{\cH_2} & = \|Af_n - Af_m\|_{\cH_2} = \|f_n - f_m\|_{A}   \no \\
& \geq \f{\varepsilon}{1+\varepsilon} \big[\|A(f_n - f_m)\|_{\cH_2} + \|f_n - f_m\|_{\cH_1}\big], 
\quad m, n \in \bbN,     \lb{2.20c}
\end{align}
$\{f_n\}_{n\in\bbN}$ and $\{A f_n\}_{n\in\bbN}$ are Cauchy sequences in $\cH_1$ and $\cH_2$, respectively. In particular, there exists $f\in\cH_1$ such that $\lim_{n\to\infty}\|f_n - f\|_{\cH_1}=0$. Since $A$ is closed, one infers $f\in\dom(A)$ and $g = \slim_{n\to\infty} A f_n = Af$, and hence $\ran(A)$ is closed in $\cH_2$. 
\end{proof}
%%%%%%%%%%%%%

Given Lemma \ref{l2.3}, we can now introduce the Hilbert space 
\begin{equation}
\cK_{A} = \ran(A) = \ol{\ran(A)} = \ker(A^*)^{\bot} \subseteq \cH_2,     \lb{2.21}
\end{equation}
and the associated projection operator $P_{\cK_{A}}$ in $\cH_2$,
\begin{equation}
P_{\cK_{A}} = [I_{\cH_2} - P_{\ker(A^*)}]. 
\end{equation} 

Next we state the following elementary result.

%%%%%%%%%%%%%
\begin{lemma} \lb{l2.4}
Assume Hypothesis \ref{h2.1} and introduce the operator
\begin{equation}
\wti A : \begin{cases} \cH_A \to \cK_{A}, \\  f \mapsto Af. \end{cases}    \lb{2.22}
\end{equation}
Then 
\begin{equation} 
\wti A \in \cB(\cH_A, \cK_{A}) \, \text{ is unitary,} 
\end{equation} 
and hence,   
\begin{equation}
\big({\wti A}\big)^{-1} : \begin{cases} \cK_{A} \to \cH_A, \\  g \mapsto A^{-1}g, \end{cases}    
\quad \big(\wti A\big)^{-1} \in \cB(\cK_{A}, \cH_A) \, \text{  is unitary.}  \lb{2.22a}
\end{equation}
\end{lemma}
%%%%%%%%%%%%%
\begin{proof}
First we note that $\ker\big(\wti A\big) = \ker(A) =\{0\}$. Next, one infers that 
\begin{equation}
\big\|\wti A f\big\|_{\cK_{A}} = \|Af\|_{\cH_2} = \|f\|_A 
= \|f\|_{\cH_A}, \quad f \in \dom(A),     \lb{2.23}
\end{equation} 
and hence $\wti A$ is isometric. Since $\ran\big(\wti A\big) = \ran(A) = \cK_{A}$, $\wti A$ is 
unitary. 
\end{proof}
%%%%%%%%%%%%%

%%%%%%%%%%%%%
\begin{lemma} \lb{l2.5}
Assume Hypothesis \ref{h2.1} and introduce the $2\times 2$ block matrix operator 
\begin{equation}
U_{\wti A} = \begin{pmatrix} 0 & I_{\cH_1} \\ -i \, \wti A & 0 \end{pmatrix}: \cH_A \oplus \cH_1 
\to \cH_1 \oplus \cK_{A}.     \lb{2.24}
\end{equation} 
Then 
\begin{equation}
U_{\wti A} \in \cB(\cH_A \oplus \cH_1, \cH_1 \oplus \cK_{A}) \, \text{ is unitary,}
\end{equation}
and hence,  
\begin{equation}
U_{\wti A}^{-1} = \begin{pmatrix} 0 & i \big({\wti A}\big)^{-1} \\ I_{\cH_1}  & 0 
\end{pmatrix} 
\in \cB(\cH_1 \oplus \cK_{A}, \cH_A \oplus \cH_1)  \, \text{ is unitary.}     \lb{2.25}
\end{equation} 
\end{lemma}
%%%%%%%%%%%%%
\begin{proof}
Assuming $f \in \cH_A$ and $g\in\cH_1$, one infers that 
\begin{align}
& \big\|U_{\wti A} (f \;\; g )^\top\big\|^2_{\cH_1 \oplus \cK_{A}} = 
\big\|U_{\wti A} (f \;\; g )^\top\big\|^2_{\cH_1 \oplus \cH_2} = 
\big\|(g \;\; -i  \, Af)^\top \big\|_{\cH_1 \oplus \cH_2}^2   \no \\
& \quad = \|g\|_{\cH_1}^2 + \|A f\|_{\cH_2}^2 = \big\|(f \;\; g)^\top\|_{\cH_A \oplus \cH_1}^2.     \lb{2.26}
\end{align} 
Thus, $\ker\big(U_{\wti A}\big) = \{0\} \oplus \ker\big(\wti A\big) = \{0\}$,  
and hence $\wti A$ is isometric. Since $\ran\big(U_{\wti A}\big) = \cH_1 \oplus \cK_{A}$, $U_{\wti A}$ is 
unitary. In addition, $U_{\wti A} U_{\wti A}^{-1} = I_{\cH_1 \oplus \cK_{A}}$ and 
$U_{\wti A}^{-1} U_{\wti A} = I_{\cH_A \oplus \cH_1}$ follow from \eqref{2.24}, \eqref{2.25}, and 
Lemma \ref{l2.4}. 
\end{proof}
%%%%%%%%%%%%%

Next, we explicitly introduce the continuous embedding operator $\iota_A$ effecting 
$\cH_A \hookrightarrow \cH_1$ by
\begin{equation}
\iota_A: \begin{cases} \cH_A \to \cH_1, \\  f \mapsto f, \end{cases} 
\end{equation} 
such that
\begin{equation}
\dom(\iota_A) = \cH_A, \quad \ran(\iota_A) = \dom(A) \subseteq \cH_1. 
\end{equation}
Then \eqref{2.7} implies
\begin{equation}
\iota_A \in \cB(\cH_A, \cH_1), \quad \|\iota_A\|_{\cB(\cH_A,\cH_1)} \leq \varepsilon^{-1/2}.   
\end{equation}
In addition, we consider
\begin{equation}
J_A = \iota_A^{-1} : \begin{cases} \cH_1 \supseteq \dom(A) \to \cH_A, \\  
 f \mapsto f. \end{cases} 
\end{equation}
We briefly summarize some properties of $J_A$.

%%%%%%%%%%%%%
\begin{lemma} \lb{l2.6}
Assume Hypothesis \ref{h2.1}. Then $J_A$ is densely defined, closed, and bijective. Moreover, 
$J_A$ is bounded if and only if $A$ is bounded, in particular, $J_A \in \cB(\cH_1,\cH_A)$ if and only 
if $A \in \cB(\cH_1,\cH_2)$. 
\end{lemma}
%%%%%%%%%%%%%
\begin{proof}
Since $J_A$ is injective and $J_A^{-1} = \iota_A \in \cB(\cH_A, \cH_1)$ is closed, so is 
$J_A$ (cf.\ \cite[p.\ 89]{We80}). 

Boundedness of $J_A$ is then equivalent to the existence of $C \in (0,\infty)$ such that 
\begin{equation}
\|J_A f\|_{\cH_A} = \|f\|_{\cH_A} = \|A f\|_{\cH_2} \leq C \|f\|_{\cH_1}, \quad f \in \dom(A),
\end{equation}
which is equivalent to $A$ being bounded. 
\end{proof}
%%%%%%%%%%%%%

With the introduction of $\iota_A$ and $J_A = \iota_A^{-1}$ one obtains
\begin{equation}
\wti A = A \, \iota_A, \quad A = \wti A \, \iota_A^{-1} = \wti A \, J_A   \lb{2.26A}
\end{equation} 
and
\begin{equation}
\dom (A^* A \iota_A) = \dom \big(A^* \wti A\big) = J_A \dom(A^* A).    \lb{2.26B} 
\end{equation}
Moreover, the following result holds.

%%%%%%%%%%%%%
\begin{lemma} \lb{l2.7}
Assume Hypothesis \ref{h2.1}. Then 
\begin{equation}
\big(A^* \wti A\big)^* = \big(A^* A \, \iota_A\big)^* = J_A,  \quad  J_A^* = A^* \wti A = A^* A \, \iota_A. 
\end{equation}
\end{lemma}
%%%%%%%%%%%%%
\begin{proof}
For brevity we denote $T = A^* \wti A = A^* A \, \iota_A$. Then $T^*$ is given by
\begin{align}
& \dom(T^*) = \{f \in \cH_1 \,|\, \text{there exists } g \in \cH_A: \, (f, T h)_{\cH_1} = (g,h)_{\cH_A}    \no \\
& \hspace*{7cm} \text{for all } h\in\dom(T)\},     \lb{2.26a} \\
& T^* f = g,      \no 
\end{align}
where 
\begin{equation}
\dom(T) = \{h \in \cH_A \,|\, A \, \iota_A h \in \dom(A^*)\}.
\end{equation}
Given $h\in \dom(T)$ and $g \in \cH_A$ as in $\dom(T^*)$ in \eqref{2.26a}, one concludes
\begin{equation}
(g,h)_{\cH_A} = (A \iota_A g, A \iota_A h)_{\cH_2} = (\iota_A g, A^* A \iota_A h)_{\cH_1} 
= (f, T h)_{\cH_1} = (f, A^* A \iota_A h)_{\cH_1},
\end{equation}
that is,
\begin{equation}
f = \iota_A g = \iota_A T^* f, \quad f \in \dom(A), 
\end{equation}
since $\ran(A^* A) = \cH$. Thus,
\begin{equation}
\dom(T^*) = \dom(A) \, \text{ and } \, T^* = \iota_A^{-1} = J_A. 
\end{equation}
Consequently, 
\begin{equation}
J_A^* = \ol T = T
\end{equation}
since $T = A^* A \, \iota_A$ is closed as $A^* A$ is closed in $\cH$, $(A^* A)^{-1} \in \cB(\cH_1)$, and 
$\iota_A \in \cB(\cH_A, \cH_1)$ (cf.\ \cite[p.\ 164]{Ka80}).
\end{proof}
%%%%%%%%%%%%%

Assuming Hypothesis \ref{h2.1}, we next introduce the operator $G_{A,0}$ in $\cH_A \oplus \cH_1$ by
\begin{align}
\begin{split} 
& G_{A,0} = \begin{pmatrix} 0 & J_A \\ - J_A^* & 0 \end{pmatrix} 
= \begin{pmatrix} 0 & J_A \\ - A^* \wti A & 0 \end{pmatrix},    \\
& \dom(G_{A,0}) = \dom \big(A^* \wti A \big) \oplus \dom(A) \subseteq \cH_A \oplus \cH_1,    \lb{2.27}
\end{split}
\end{align}
where, 
\begin{equation}
\dom \big(A^* \wti A \big) = \big\{ f \in \cH_A \, \big| \, \wti A f \in \dom(A^*) \big\}. 
\end{equation}
In particular, one infers
\begin{equation}
G_{A,0} = \begin{pmatrix} 0 & J_A \\ - A^* A \, \iota_A & 0 \end{pmatrix}   
=  \begin{pmatrix} 0 & I_{\cH_A} \\ - A^* A & 0 \end{pmatrix}  
\begin{pmatrix} \iota_A & 0 \\  0 & J_A \end{pmatrix}. 
\end{equation}
We recall that the Hilbert space $\cH_A \oplus \cH_1$ in connection with $G_{A,0}$ is sometimes called the {\it energy space}. 

In addition, still assuming Hypothesis \ref{h2.1}, we introduce the supersymmetric Dirac-type operator $Q_{A,0}$ in $\cH_1 \oplus \cH_2$ by
\begin{equation}
Q_{A,0} = \begin{pmatrix} 0 & A^* \\ A & 0 \end{pmatrix}, \quad 
\dom(Q_{A,0}) = \dom(A) \oplus \dom(A^*) \subseteq \cH_1 \oplus \cH_2.    \lb{2.28}
\end{equation}
As discussed in Appendix \ref{sA}, $Q_{A,0}$ is self-adjoint in $\cH_1 \oplus \cH_2$. Moreover, \eqref{A.21} and $\ker(A) = \{0\}$ yield
\begin{align}
Q_{A,0} [I_{\cH_1\oplus\cH_2} - P_{\ker(Q_{A,0})}] &= \begin{pmatrix} 0 & A^* \\ A & 0 \end{pmatrix} 
\begin{pmatrix} I_{\cH} - P_{\ker(A)} & 0 \\ 0 &  I_{\cH} - P_{\ker(A^*)} \end{pmatrix}   \no \\
& = \begin{pmatrix} 0 & A^*[I_{\cH} - P_{\ker(A^*)}] \\ A & 0 \end{pmatrix}
= \begin{pmatrix} 0 & A^*P_{\cK_{A}} \\ A & 0 \end{pmatrix}.     \lb{2.29}
\end{align} 
Clearly, 
\begin{align}
\begin{split}
Q_{A,0} [I_{\cH_1\oplus\cH_2} - P_{\ker(Q_{A,0})}] &= [I_{\cH_1\oplus\cH_2} - P_{\ker(Q_{A,0})}] Q_{A,0}   \\
&= [I_{\cH_1\oplus\cH_2} - P_{\ker(Q_{A,0})}] Q_{A,0} [I_{\cH_1\oplus\cH_2} - P_{\ker(Q_{A,0})}]     \lb{2.30}
\end{split}
\end{align}
is self-adjoint in $\cH_1 \oplus \cH_2$, with $\cH_1 \oplus \cK_{A}$ a reducing (i.e., invariant) subspace 
for $Q_{A,0} [I_{\cH_1\oplus\cH_2} - P_{\ker(Q_{A,0})}]$. In this context we also note that 
\begin{equation}
[I_{\cH} - P_{\ker(A^*)}] A = P_{\cK_{A}} A = A.   \lb{2.31}
\end{equation}

At this point we are in position to formulate our first principal result and establish the following 
remarkable connection between the generator $G_{A,0}$ and the abstract supersymmetric 
Dirac-type operator $Q_{A,0}$.

%%%%%%%%%%%%%%
\begin{theorem} \lb{t2.8}
Assume Hypothesis \ref{h2.1}. Then
\begin{equation}
Q_{A,0} [I_{\cH_1\oplus\cH_2} - P_{\ker(Q_{A,0})}] 
=  U_{\wti A} \, i \, G_{A,0} U_{\wti A}^{-1}.    
\lb{2.32}
\end{equation}
In particular, the operator $i \, G_{A,0}$ is self-adjoint in the energy space 
$\cH_A \oplus \cH_1$ and hence $G_{A,0}$ generates a unitary group 
$e^{G_{A,0} t}$, $t\in\bbR$, in $\cH_A \oplus \cH_1$. Moreover, 
$G_{A,0}$ is unitarily equivalent to $- G_{A,0}$. 
\end{theorem}
%%%%%%%%%%%%%%
\begin{proof}
Self-adjointness of $i \, G_{A,0}$ is an immediate consequence of $J_A^* = A^* \wti A$ in 
Lemma \ref{l2.7} and the first equality in \eqref{2.27}, that is,
\begin{equation}
G_{A,0} = \begin{pmatrix} 0 & J_A \\ - A^* \wti A & 0 \end{pmatrix} 
= \begin{pmatrix} 0 & J_A \\ - J_A^* & 0 \end{pmatrix}     \lb{2.32a}
\end{equation}
in $\cH_A \oplus \cH_1$ and the fact that $J_A$ is closed by Lemma \ref{l2.6}. 

Employing the fact that $\dom (G_{A,0}) = \dom \big(A^* \wti A) \oplus \dom(A) \subseteq \cH_A \oplus \cH_1$, one first obtains
\begin{align}
& U_{\wti A} \dom (G_{A,0}) = \begin{pmatrix} 0 & I_{\cH_1} \\ -i \, \wti A & 0 \end{pmatrix} 
\bigg\{\begin{pmatrix} f \\ g \end{pmatrix} \in \cH_A \oplus \cH_1 \, \bigg| \, f \in \dom \big(A^* \wti A \big), 
\, g \in \dom(A) \bigg\}   \no \\
& \quad = \big\{(g \; -i \, \wti A f \big)^\top \in \cH_1 \oplus \cK_{A} \, \big| \, f \in \dom \big(\wti A\big), 
\, \wti A f \in \dom(A^*) , \, g \in \dom (A)\big\}    \no \\
& \quad = \dom (A) \oplus \{h \in \cK_{A} \, | \, h \in \dom(A^*) \}    \no \\
& \quad = \dom (A) \oplus \dom (A^* [I_{\cH_2} - P_{\ker(A^*)}])   \no \\
& \quad = \dom (A) \oplus \dom (A^* P_{\cK_{A}})    \no \\
& \quad = \dom (Q_{A,0} [I_{\cH_1\oplus\cH_2} - P_{\ker(Q_{A,0})}]).     \lb{2.33}
\end{align}
Next, one computes for $(f \; g)^{\top} \in \cH_1 \oplus \cK_{A}$ such that 
$U_{\wti A}^{-1} (f \; g)^{\top} \in \dom (G_{A,0}) = \dom \big(A^* \wti A) \oplus \dom(A)$, 
\begin{align}
& U_{\wti A}^{-1} \begin{pmatrix} f \\ g \end{pmatrix} 
= \begin{pmatrix} -i \big({\wti A}\big)^{-1} g \\ f \end{pmatrix} 
\in \dom (G_{A,0}) = \dom \big(A^* \wti A) \oplus \dom(A)  \no \\
& \quad \text{if and only if } \, f \in \dom(A), \; 
\big({\wti A}\big)^{-1} g \in \dom \big(A^* \wti A\big), \; g \in \cK_{A} = \ran \big(\wti A\big)   \no \\
& \quad \text{if and only if } \, f \in \dom(A), \; g \in \cK_{A}, \; g \in \dom (A^*)    \no \\
& \quad \text{if and only if } \, f \in \dom(A), \; g \in \dom (A^* P_{\cK_{A}}) 
= \dom (A^* [I_{\cH_2} - P_{\ker(A^*)}]),    \lb{2.34}
\end{align}
where we used the fact that $\cK_{A} = P_{\cK_{A}} \cH_2 = [I_{\cH_2} - P_{\ker(A^*)}] \cH_2$ 
reduces $A^*$. Thus, 
\begin{align}
U_{\wti A} i \, G_{A,0} U_{\wti A}^{-1} 
 &= i \begin{pmatrix} 0 & I_{\cH_1} \\ -i \, \wti A & 0 \end{pmatrix} 
\begin{pmatrix} 0 & J_A \\ - A^* \wti A & 0 \end{pmatrix} 
\begin{pmatrix} 0 & i \big({\wti A}\big)^{-1} \\ I_{\cH_1} & 0 \end{pmatrix}  \no \\
& = i \begin{pmatrix} 0 & I_{\cH_1} \\ -i \, \wti A & 0 \end{pmatrix} 
\begin{pmatrix} J_A & 0 \\ 0 & -i \, A^* P_{\cK_{A}} \end{pmatrix}   \no \\
& = \begin{pmatrix} 0 & A^* P_{\cK_{A}} \\ \wti A \, J_A & 0 \end{pmatrix}  
= \begin{pmatrix} 0 & A^* P_{\cK_{A}} \\ \wti A \, \iota_A^{-1} & 0 \end{pmatrix}   \no \\
& = \begin{pmatrix} 0 & A^* P_{\cK_{A}} \\  A & 0 \end{pmatrix},   \lb{2.35}
\end{align}
using $\wti A \, J_A = \wti A \, \iota_A^{-1} = A$ by \eqref{2.26A}.

An alternative proof of the self-adjointness of $i \, G_{A,0}$ then follows from \eqref{2.32} and the self-adjointness of $Q_{A,0}$ (cf.\ \eqref{A.2}) and hence that of 
$Q_{A,0} [I_{\cH_1 \oplus \cH_2} - P_{\ker(Q_{A,0})}]$. 

Finally, the unitary equivalence of $Q_{A,0}$ to $- Q_{A,0}$ in \eqref{A.22} together with 
\eqref{A.21}, which implies the unitary equivalence of the operators 
$Q_{A,0} [I_{\cH_1\oplus\cH_2} - P_{\ker(Q_{A,0})}]$ 
and $-Q_{A,0} [I_{\cH_1\oplus\cH_2} - P_{\ker(Q_{A,0})}]$, and \eqref{2.32} then prove the unitary equivalence of $G_{A,0}$ and $- G_{A,0}$.
\end{proof}
%%%%%%%%%%%%%%

%%%%%%%%%%%%%%
\begin{remark} \lb{r2.9}
$(i)$ Given Hypothesis \ref{h2.1}, the self-adjointness of $G_{A,0}$ on 
$\dom (G_{A,0}) = \dom \big(A^* \wti A) \oplus \dom(A)$ in the energy space $\cH_A \oplus \cH_1$ 
is of course well-known. We refer, for instance, to the monographs \cite[Sect.\ VI.3]{EN00}, 
\cite[Sect.\ 2.7]{Go85}, \cite[p.\ 2, 3]{Re76}, \cite[Sect.\ X.13]{RS75}. These sources typically employ a combination of semigroup methods and the spectral theorem for self-adjoint operators. Our proof of 
\eqref{2.32a} closely follows the pattern displayed in the Klein--Gordon context in 
\cite[Subsect.\ 5.5.3]{Th92}.   
Our proof based on the unitary equivalence to the self-adjoint Dirac-type operator 
$Q_{A,0} [I_{\cH_1\oplus\cH_2} - P_{\ker(Q_{A,0})}]$ in $\cH_1 \oplus \cH_2$ appears to be a new twist in this context. \\
$(ii)$ The observation that $G_{A,0}$ (and more generally, $G_{A,R}$) in the energy space 
$\cH_A \oplus \cH_1$ is related to a Dirac-type operator in $\cH_1 \oplus \cH_2$ has recently been 
made in the context of trace formulas for the damped string equation \cite{GH10}. However, this observation is not new and has already been made in \cite{KL78} (under more restricted assumptions of compactness of $A^*A$ and self-adjointness and boundedness of $R$) and \cite{GS79}, 
and was subsequently also discussed in \cite{GLS93}, \cite{GS82a}, 
\cite{GS82}, \cite{Hu88}, and \cite[Subsect.\ 5.5.3]{Th92}. We have not been able to find the precise unitary equivalence result 
\eqref{2.32} in Theorem \ref{t2.8} in the literature. The fact that $G_{A,0}$ 
and $- G_{A,0}$ are similar 
operators has been noted in \cite[p.\ 382]{EN00}. 
\end{remark}
%%%%%%%%%%%%%%

Still assuming the basic Hypothesis \ref{h2.1}, we now briefly summarize the basic results derived thus far if $A$ and $A^*$ in the factorization $A^* A$ are both systematically replaced by $|A|$ using the fact that $A^* A = |A|^2$. This case is of considerable interest and used in practice as 
$0 < \varepsilon^{-1/2} I_{\cH_1} \leq |A|^{-1} \in \cB(\cH_1)$, whereas $A$ is in general not boundedly invertible as discussed in 
Example \ref{e2.2}. Since this is a special case of the discussion thus far, we now focus on some of the simplifications that arise in this context and present the results without proofs as the latter parallel those that have already been presented in great detail.  

We start by noting that in this special case 
\begin{align}
\cH_{|A|} &= \cH_A,   \lb{2.59a} \\
\iota_{|A|} & = \iota_A,     \lb{2.60a} \\
J_{|A|} & = J_A,    \lb{2.61} \\
\cK_{|A|} &= \ol{\ran(|A|)} = \ker(|A|)^\bot = \cH_1.   \lb{2.62}
\end{align}
In addition\footnote{We emphasize that $\wti{|A|}$ as defined in \eqref{2.63} 
differs of course from $\big|\wti A\big| = \big(\big(\wti A\big)^* \wti A\big)^{1/2}$
(noting the different order of operations). In fact, since $\wti A$ as defined 
in \eqref{2.22} is unitary, one has $\big|\wti A\big| = I_{\cH_A}$.}, 
\begin{align}
& \wti {|A|} : \begin{cases} \cH_{|A|} \to \cH_1, \\  f \mapsto |A| f, \end{cases} \quad  
\big({\wti{|A|}}\big)^{-1} : \begin{cases} \cH_1 \to \cH_{|A|}, \\  g \mapsto |A|^{-1}g, \end{cases}    
\lb{2.63}  \\
& \wti {|A|} \in \cB(\cH_{|A|}, \cH_1), \; \big({\wti{|A|}}\big)^{-1} \in \cB(\cH_1, \cH_{|A|}) \, 
\text{ are both unitary,}     \lb{2.64a} \\
& U_{\wti{|A|}} = \begin{pmatrix} 0 & I_{\cH_1} \\ -i \, \wti {|A|} & 0 \end{pmatrix} 
\in \cB(\cH_{|A|} \oplus \cH_1, \cH_1 \oplus \cH_1) \, \text{ is unitary,}     \lb{2.65} \\
& U_{\wti{|A|}}^{-1} = \begin{pmatrix} 0 & i \big({\wti{|A|}}\big)^{-1} \\ I_{\cH_1}  & 0 \end{pmatrix} 
\in \cB(\cH_1 \oplus \cH_1, \cH_{|A|} \oplus \cH_1)  \, \text{ is unitary,}     \lb{2.66} \\
& J_{|A|}^* = |A|^2 \iota_{|A|} = A^*A \iota_A = J_A^*,    \lb{2.67}  \\
& \wti{|A|} = |A| \iota_A, \quad \big({\wti{|A|}}\big)^{-1} = J_A |A|^{-1},     \lb{2.67a} \\
& G_{|A|,0} = \begin{pmatrix} 0 & J_{|A|} \\ - |A|^2 \iota_{|A|} & 0 \end{pmatrix} 
= \begin{pmatrix} 0 & J_A \\ - |A|^2 \iota_A & 0 \end{pmatrix} 
= \begin{pmatrix} 0 & J_A \\ - A^* A \iota_A & 0 \end{pmatrix} = G_{A,0},   
\lb{2.68} \\
& \dom(G_{|A|,0}) = \dom \big(|A|^2 \iota_A \big) \oplus \dom(A) = \dom(G_{A,0}) 
\subseteq \cH_{|A|} \oplus \cH_1,   \lb{2.69} \\
& Q_{|A|,0} = \begin{pmatrix} 0 & |A| \\ |A| & 0 \end{pmatrix}, \quad 
\dom (Q_{|A|,0}) = \dom(|A|) \oplus \dom(|A|) \subseteq \cH_1 \oplus \cH_1.   \lb{2.70} 
\end{align}
Consequently, one obtains as in Theorem \ref{t2.8} that 
\begin{equation}
Q_{|A|,0} =  U_{\wti{|A|}} \, i \, G_{A,0} U_{\wti{|A|}}^{-1},  \quad 
\dom (Q_{|A|,0}) = U_{\wti{|A|}} \dom (G_{|A|,0}).    \lb{2.71} 
\end{equation}

\smallskip

We emphasize that $Q_{|A|,0}$ in \eqref{2.71} does not involve any additional projection as opposed 
to $Q_{A,0} [I_{\cH_1\oplus\cH_2} - P_{\ker(Q_{A,0})}]$ in \eqref{2.32}. Still, the two operators are of course unitarily equivalent. Indeed, equation \eqref{2.71} implies 
\begin{equation}
Q_{|A|,0} =  \big[U_{\wti{|A|}} U_{\wti A}^{-1}\big] 
Q_{A,0} [I_{\cH_1\oplus\cH_2} - P_{\ker(Q_{A,0})}] 
\big[U_{\wti{|A|}} U_{\wti A}^{-1}\big]^{-1},   \lb{2.73}
\end{equation}
where
\begin{align}
U_{\wti{|A|}} U_{\wti A}^{-1} & 
= \begin{pmatrix} I_{\cH_1} & 0 \\ 0 & \wti{|A|} \big(\wti{A}\big)^{-1} \end{pmatrix} 
= \begin{pmatrix} I_{\cH_1} & 0 \\ 0 & |A| A^{-1} \end{pmatrix} 
 = \begin{pmatrix} I_{\cH_1} & 0 \\ 0 & (V_A)^* \end{pmatrix}  \no \\
& = \begin{pmatrix} I_{\cH_1} & 0 \\ 0 & V_{A^*} \end{pmatrix} 
\in \cB(\cH_1 \oplus \cK_{A},\cH_1 \oplus \cH_1) \, \text{ is unitary}, 
\end{align}
using, \eqref{2.3}, \eqref{A.4a}--\eqref{A.4e}, and employing the fact that the 
initial set of $V_{A^*}$ coincides with $\ol{\ran (A)} = \cK_{A}$. 

We note that the \eqref{2.71} is a special case of a result observed by 
Huang \cite{Hu97} in connection with his Proposition\ 3.1 (the latter also 
includes a damping term $R$, see also Theorem \ref{t4.4}).

%%%%%%%%%%%%%%%%%%%%%%%%%%%%%%%%%%%%%%
%%%%%%%%%%%%%%%%%%%%%%%%%%%%%%%%%%%%%%
\section{Abstract Linear Wave Equations in the Absence of Damping. \\ 
The Case $\inf(\sigma(A^* A)) = 0$}
\label{s3}
%%%%%%%%%%%%%%%%%%%%%%%%%%%%%%%%%%%%%%
%%%%%%%%%%%%%%%%%%%%%%%%%%%%%%%%%%%%%%

In this section we indicate how to extend the results of the previous section to the case 
$\inf(\sigma(A^* A)) = 0$. This case will to a large extend parallel the case 
$A^* A \geq \varepsilon I_{\cH_1}$ for some $\varepsilon > 0$, and hence we will 
mainly focus on the differences between these two situations. 

Our basic hypothesis throughout this section now reads as follows. 

%%%%%%%%%%%
\begin{hypothesis} \lb{h3.1}
Let $\cH_j$, $j=1,2$, be complex separable Hilbert spaces.
Assume that $A: \cH_1 \supseteq \dom(A) \to \cH_2$ is a densely 
defined, closed, linear operator satisfying 
\begin{equation}
\ker(A) = \{0\}   \lb{3.1}
\end{equation}
and 
\begin{equation}
\inf(\sigma(A^* A)) = 0.     \lb{3.2}
\end{equation}
\end{hypothesis}
%%%%%%%%%%% 

As in the previous case we can equip $\dom(A)$ with the norm $\| \dott \|_{A}$, but since the stronger Hypothesis \ref{h2.1} is no longer assumed, the resulting space will in general not be complete. Hence we
denote by $\cH_A$ its completion, 
\begin{equation}
\cH_A = \ol{(\dom(A); (\cdot,\cdot)_A)}, \quad (f,g)_A = (Af, Ag)_{\cH_2}, \; 
f,g \in \dom(A) \subseteq \cH_1.     \lb{3.3}
\end{equation}
In general (cf.\ Example \ref{e5.7}), 
\begin{equation}
\cH_A \nsubseteq \cH_1 \, \text{ and } \, \cH_1 \nsubseteq \cH_A.   \lb{3.3a}
\end{equation}   
Moreover, Lemma \ref{l2.3} will also fail in general and
consequently we now define
\begin{equation}
\cK_{A} = \ol{\ran(A)} = \ker(A^*)^{\bot} \subseteq \cH_2.     \lb{3.4}
\end{equation}
Next, Lemma \ref{l2.4} also requires some modifications.

%%%%%%%%%%%%%
\begin{lemma} \lb{l3.2}
Assume Hypothesis \ref{h3.1} and introduce the operator
\begin{equation}
A_0 : \begin{cases} \cH_A \supseteq \dom(A) \to \cK_{A}, \\  f \mapsto Af. \end{cases}    
\lb{3.5}
\end{equation}
Then there exists a $($unique$)$ unitary extension 
$\wti A = \ol{A_0} \in \cB(\cH_A, \cK_{A})$ of $A_0$. \\
\end{lemma}
%%%%%%%%%%%%%
\begin{proof}
As in the proof of Lemma \ref{l2.4} one infers that $A_0$ is isometric. Since
$\ran (A_0) = \ran(A) \subseteq \cK_{A}$ is dense, there is a unique unitary extension $\wti A$ of $A_0$ given by the closure $\ol{A_0}$ of $A_0$. 
\end{proof}
%%%%%%%%%%%%%

Consequently, Lemma \ref{l2.5} extends without further modifications to the 
present setting.

%%%%%%%%%%%%%
\begin{lemma} \lb{l3.3}
Assume Hypothesis \ref{h3.1} and introduce the $2\times 2$ block matrix operator 
\begin{equation}
U_{\wti A} = \begin{pmatrix} 0 & I_{\cH_1} \\ -i \, \wti A & 0 \end{pmatrix}: \cH_A \oplus \cH_1 
\to \cH_1 \oplus \cK_{A}. 
\end{equation}
Then 
\begin{equation}
U_{\wti A} \in \cB(\cH_A \oplus \cH_1, \cH_1 \oplus \cK_{A}) \, \text{ is unitary,}
\end{equation}
and hence,  
\begin{equation}
U_{\wti A}^{-1} = \begin{pmatrix} 0 & i \big({\wti{A}}\big)^{-1} \\ I_{\cH_1}  & 0 \end{pmatrix} 
\in \cB(\cH_1 \oplus \cK_{A}, \cH_A \oplus \cH_1)  \, \text{ is unitary.} 
\end{equation} 
\end{lemma}
%%%%%%%%%%%%%

We can also introduce the embedding operator $\iota_A$ effecting the embedding 
$ \cH_A  \supseteq \dom(A) \hookrightarrow \cH_1$ by
\begin{equation}
\iota_A: \begin{cases} \cH_A \supseteq \dom(A) \to \cH_1, \\  f \mapsto f, \end{cases} 
\end{equation} 
such that
\begin{equation}
\dom(\iota_A) = \dom(A) \subseteq \cH_A , \quad \ran(\iota_A) = \dom(A) \subseteq \cH_1. 
\end{equation}
In particular, we note that $\iota_A$ is no longer a bounded operator unless 
Hypothesis \ref{h2.1} holds. In addition, we consider
\begin{equation}
J_A = \iota_A^{-1} : \begin{cases}\cH_1 \supseteq  \dom(A) \to \cH_A, \\  
f \mapsto f. \end{cases} 
\end{equation}
Both $\iota_A$ and $J_A$ are densely defined, closed, and bijective.

With the introduction of $\iota_A$ and $J_A = \iota_A^{-1}$ one obtains
\begin{equation}
A_0 = A \iota_A, \quad \wti A = \ol{A \, \iota_A}, 
\quad A = \wti A \, \iota_A^{-1} = \wti A \, J_A,    \lb{3.9}
\end{equation} 
and the analog of Lemma \ref{l2.7} holds. 

%%%%%%%%%%%%%
\begin{lemma} \lb{l3.4}
Assume Hypothesis \ref{h3.1}. Then 
\begin{equation}
\big(A^* \wti A\big)^* = J_A,  \quad  J_A^* = A^* \wti A. 
\end{equation}
\end{lemma}
%%%%%%%%%%%%%
\begin{proof}
Since $\wti A$ is unitary and $A$ is closed, one computes 
(cf.\ \cite[Exercise 4.18]{We80})
\begin{equation}
\big(A^* \wti A\big)^* = \big(\wti A\big)^* A = \big(\wti{A}\big)^{-1} A
=  \big(\ol{A \, \iota_A}\big)^{-1} A = \ol{\big(A \, \iota_A\big)^{-1}} A =
\ol{J_A A^{-1}} A.
\end{equation}
In addition, since $\dom(J_A A^{-1})=\ran(A)$, one can drop the closure in
the last equation which finally yields $\big(A^* \wti A\big)^* = J_A A^{-1} A
= J_A$. 

Hence one also obtains $J_A^* = \big(A^* \wti A\big)^{**} = \ol{A^* \wti A} = A^* \wti A$
since $A^* \wti A$ is closed as $\wti A$ is unitary and $A^*$ is closed.
\end{proof}
%%%%%%%%%%%%%%%

Assuming Hypothesis \ref{h3.1}, we again introduce the operator $G_{A,0}$ in 
$\cH_A \oplus \cH_1$ by
\begin{align}
\begin{split} 
& G_{A,0} = \begin{pmatrix} 0 & J_A \\ - J_A^* & 0 \end{pmatrix}
= \begin{pmatrix} 0 & J_A \\ - A^* \wti A & 0 \end{pmatrix},    \\
& \dom(G_{A,0}) = \dom \big(A^* \wti A \big) \oplus \dom(A) \subseteq \cH_A \oplus \cH_1,    \lb{3.12}
\end{split}
\end{align}
and also introduce the supersymmetric Dirac-type operator $Q_{A,0}$ in 
$\cH_1 \oplus \cH_2$ by
\begin{equation}
Q_{A,0} = \begin{pmatrix} 0 & A^* \\ A & 0 \end{pmatrix}, \quad 
\dom(Q_{A,0}) = \dom(A) \oplus \dom(A^*) \subseteq \cH_1 \oplus \cH_2.    \lb{3.13}
\end{equation}
As discussed in Appendix \ref{sA}, $Q_{A,0}$ is self-adjoint in 
$\cH_1 \oplus \cH_2$.

The analog of Theorem \ref{t2.8} then reads as follows.

%%%%%%%%%%%%%%
\begin{theorem} \lb{t3.5}
Assume Hypothesis \ref{h3.1}. Then
\begin{equation}
Q_{A,0} [I_{\cH_1\oplus\cH_2} - P_{\ker(Q_{A,0})}] =  U_{\wti A} \, i \, 
G_{A,0} U_{\wti A}^{-1}.    
\lb{3.14}
\end{equation}
In particular, the operator $i \, G_{A,0}$ is self-adjoint in the energy space $\cH_A \oplus \cH_1$ and hence generates a unitary group $e^{G_{A,0} t}$, $t\in\bbR$, in $\cH_A \oplus \cH_1$. Moreover, 
$G_{A,0}$ is unitarily equivalent to $- G_{A,0}$. 
\end{theorem}
%%%%%%%%%%%%%%

Next, we further analyze the domain of $G_{A,0}$, more precisely, the domain 
of $A^* {\wti A}$ (cf.\ \eqref{3.12}), applying some results 
discussed in Appendix \ref{sB}. Since $\wti A = \ol{A_0}$, and $A^* \wti A$ 
is known to be a closed operator (cf.\ \eqref{3.12}), the natural question arises whether or not $A^* \wti A = A^* \ol{A_0}$ coincides with the closure 
$\ol{A^* A_0}$ of $A^* A_0$. This is a somewhat intricate question, an answer to which is given in Theorem \ref{t3.8} below.

We start with the following elementary result.

%%%%%%%%%%%%%%
\begin{lemma} \lb{l3.6}
Suppose $S$ is self-adjoint in the complex separable Hilbert space $\cH$ with 
$\ker(S) = \{0\}$. Then
\begin{equation}
\dom(S) \cap \ran(S) = \dom(S) \cap \dom\big(S^{-1}\big) \, 
\text{ is dense in $\cH$ and a core for $S$ and $S^{-1}$.}   \lb{3.15}
\end{equation}
\end{lemma}
%%%%%%%%%%%%%%
\begin{proof}
Since $\ker(S) = \{0\}$, the operator $S^{-1}$ exists and is self-adjoint (and also 
$\ker\big(S^{-1}\big) = \{0\}$). For any $g \in \cH$, $g_n = 
E_S ([-n, - n^{-1}] \cup [n^{-1},n]) g \in \dom(S) \cap \dom\big(S^{-1}\big)$, 
$n\in\bbN$, and hence 
\begin{equation}
\lim_{n\to\infty} \|g_n - g\|_{\cH} 
= \lim_{n\to\infty} \|[E_S ([-n, - n^{-1}] \cup [n^{-1},n]) - I_{\cH}] g\|_{\cH}= 0
\end{equation}
proves that $\ol{\dom(S) \cap \dom\big(S^{-1}\big)} = \cH$. Here $E_S (\cdot)$ denotes the strongly right continuous family of spectral projections associated 
with $S$. 

Next, let $f \in \dom(S)$ 
and introduce $f_n = E_S ((-\infty, - n^{-1}] \cup [n^{-1},\infty)) f \in \dom(S) 
\cap \dom\big(S^{-1}\big)$, $n\in\bbN$. Then 
\begin{align}
\begin{split}
& \lim_{n\to\infty} \|f_n - f\|_{\cH} 
= \lim_{n\to\infty} \|[E_S ((-\infty, - n^{-1}] \cup [n^{-1},\infty)) - I_{\cH}] f\|_{\cH}= 0,   
\lb{3.16} \\ 
& \lim_{n\to\infty} \|S f_n - S f\|_{\cH}
= \lim_{n\to\infty} \|[E_S ((-\infty, - n^{-1}] \cup [n^{-1},\infty)) - I_{\cH}]S f\|_{\cH} = 0
\end{split}
\end{align}  
prove that $\dom(S) \cap \dom\big(S^{-1}\big)$ is a core for $S$ since 
$f \in \dom(S)$ was arbitrary. By symmetry between $S$ and $S^{-1}$, 
$\dom(S) \cap \dom\big(S^{-1}\big)$ is also a core for $S^{-1}$. 
\end{proof}
%%%%%%%%%%%%%%

The next lemma is of an auxiliary nature and together with Lemma \ref{l3.6} the basic ingredient for the proof of Theorem \ref{t3.8} below.

%%%%%%%%%%%%%%
\begin{lemma} \lb{l3.7}
Assume Hypothesis \ref{h3.1} and denote by $P_{\cK_A}$ the orthogonal projection onto $\cK_A$ in $\cH_2$. \\
$(i)$ Suppose that $\ol{\dom(A^* P_{\cK_A}) \cap \ran(A)} = \cK_A$, then 
\begin{equation}
\ol{A^* A_0} = \ol{A^* P_{\cK_A}\big|_{\dom(A^* P_{\cK_A}) \cap \ran(A)}} 
\, \ol{A_0}.     \lb{3.19}
\end{equation}
$(ii)$ Assume that $\ker(A^*) = \{0\}$. Then $\cK_A = \cH_2$, 
$P_{\cK_A} = I_{\cH_2}$, 
\begin{align} 
& \ol{\dom(A^*) \cap \ran(A)} = \cH_2,    \lb{3.19a} \\
& \dom(A^*) \cap \ran(A) \, \text{ is a core for $A^*$,}   \lb{3.19c}
\end{align}
and 
\begin{equation}
\ol{A^* A_0} = A^* \ol{A_0}.   \lb{3.20} 
\end{equation}
In particular, if $A$ is self-adjoint in $\cH_1$ satisfying \eqref{3.1} and \eqref{3.2}, 
then \eqref{3.19a}--\eqref{3.20} hold with $\cH_2 = \cK_A = \cH_1$.
\end{lemma} 
%%%%%%%%%%%%%%
\begin{proof}
$(i)$ By general principles, $A^* A_0 \subseteq A^* \ol{A_0}$ implies 
\begin{equation}
\ol{A^* A_0} \subseteq \ol{A^* \ol{A_0}} = A^* \ol{A_0},     \lb{3.21} 
\end{equation}
as the latter is a closed operator (cf.\ \eqref{3.12}). The reverse inclusion is 
more subtle, though.

Since $\wti A = \ol{A_0}$ is unitary, 
$\ran(A_0) = \ran(A)$, and $A^*$ 
is closed, one can apply Lemma \ref{lB.1}\,$(iv)$ to obtain
\begin{equation}
(A^* A_0)^* = (A^* P_{\cK_A} A_0)^* = A_0^* \Big(\hatt{A^* P_{\cK_A}}\Big)^*
= A_0^* \big(A^* P_{\cK_A}\big|_{\dom(A^* P_{\cK_A}) \cap \ran(A)}\big)^*. 
\lb{3.22}
\end{equation}
Using unitarity of $A_0^* = \big(\ol{A_0}\big)^*$, and applying 
Lemma \ref{lB.1}\,$(ii)$ one finally obtains
\begin{equation}
\ol{A^* A_0} = ((A^* P_{\cK_A} A_0)^*)^* 
= \ol{A^* P_{\cK_A}\big|_{\dom(A^* P_{\cK_A}) \cap \ran(A)}} \, 
\ol{A_0},   \lb{3.23}
\end{equation}
employing $\ol T = (T^*)^*$, whenever $T$ is densely defined and closable. \\
$(ii)$ Next, one recalls the 
fact that for any densely defined closed operator $T$ in $\cH_1$ mapping into 
$\cH_2$ one has (cf.\ \cite[p.\ 335]{Ka80}, \cite[Theorem IV.3.2]{EE89})
\begin{equation}
\ran(T) = \ran (|T^*|).    \lb{3.29}
\end{equation}
Equation \eqref{3.29} is a consequence of the polar decompositions for $T$ and 
$T^*$, more precisely, of 
\begin{equation}
T = |T^*| U_T, \quad |T^*| = T U_T^*, \quad |T^*| = U_T |T| U_T^*,   \lb{3.30} 
\end{equation}
where $U_T$ is a partial isometry with initial set $\ol{\ran(|T|)}$ and final set 
$\ol{\ran(T)}$ (and hence $U_T^*$ is a partial isometry with initial set 
$\ol{\ran(T)}$ and final set  $\ol{\ran(|T|)}$). Using the fact that 
$\dom(T) = \dom(|T|)$ and applying \eqref{3.29} to $T=A$, 
one concludes from Lemma \ref{l3.6} and the fact that by hypothesis 
$\ker(A^*) = \ker(|A^*|) = \{0\}$ and hence $P_{\cK_A} = I_{\cH_2}$,   
\begin{equation}
\dom(A^*) \cap \ran(A) = \dom(|A^*|) \cap \ran (|A^*|) \, 
\text{ is dense in $\cH_2$ and a core for $|A^*|$.}     \lb{3.31}
\end{equation}
The polar decomposition for $T^*$, $T^* = U_T^* |T^*|$ then immediately yields 
that 
\begin{equation}
\cD_0 = \dom(A^*) \cap \ran(A) = \dom(|A^*|) \cap \ran (|A^*|) \, 
\text{ is a core for $A^*$.}     \lb{3.32} 
\end{equation} 
Indeed, if $f \in \dom(A^*)=\dom(|A^*|)$, there exists $f_n \in \cD_0$, $n \in\bbN$, 
such that 
\begin{equation}
\lim_{n\to\infty} \|f_n - f\|_{\cH} = 0 \, \text{ and } \, 
\lim_{n\to\infty} \||A^*| f_n - |A^*| f\|_{\cH} = 0,     \lb{3.33}
\end{equation} 
and hence also,
\begin{equation}
\lim_{n\to\infty} \|A^* f_n - A^* f\|_{\cH} = 
\| U_A^*|A^*| f_n - U_A^* |A^*| f\|_{\cH} = 0,     \lb{3.34} 
\end{equation} 
proving that $\cD_0$ is a core for $A^*$. Thus, \eqref{3.19} then yields  
\begin{equation}
\ol{A^* A_0} = \ol{A^*\big|_{\dom(A^*) \cap \ran(A)}} \, \ol{A_0}
= A^* \, \ol{A_0},   \lb{3.35} 
\end{equation}
and hence proves \eqref{3.20}. 
\end{proof}
%%%%%%%%%%%%%%

Since $\cH_2 = \ol{\ran(A)} \oplus \ker(A^*) = \cK_A \oplus \ker(A^*)$, one can 
 introduce the operator 
 \begin{equation}
 B : \begin{cases} \cH_1 \supseteq \dom(A) \to \cK_A, \\
 f \mapsto A f,   \end{cases}    \lb{3.36}
\end{equation}
and then concludes that 
\begin{equation}
\ker(B^*) = \ker(|B^*|) = \ran(B)^\perp = \{0\},     \lb{3.37} 
\end{equation}
and that
\begin{equation} 
A = P_{\cK_A} B.   \lb{3.38} 
\end{equation} 
Thus, an application of Lemma \ref{lB.1}\,$(ii)$ yields 
\begin{equation} 
A^* = B^* P_{\cK_A}.    \lb{3.39}
\end{equation}
Moreover, one verifies that 
\begin{equation} 
\cK_A = \cK_B, \quad B_0 = A_0,  \text{ and hence, } \, 
\wti B = \wti A.     \lb{3.40}
\end{equation} 

Given the preparatory Lemmas \ref{l3.6} and \ref{l3.7} we finally are in a 
position to formulate the following result, a resolution of the question posed in 
the paragraph preceding Lemma \ref{l3.6}. 

%%%%%%%%%%%%%%
\begin{theorem} \lb{t3.8}
Assume Hypothesis \ref{h3.1}. Then 
\begin{equation}
 \ol{A^* A_0} = A^* \ol{A_0}.     \lb{3.41}
 \end{equation} 
\end{theorem} 
%%%%%%%%%%%%%%
\begin{proof}
Since $\ker(|B^*|) =  \{0\}$ and $|B^*|$ is self-adjoint in $\cK_B = \cK_A$, 
\eqref{3.31} and \eqref{3.32} apply and yield  
\begin{align}
\begin{split}
& \dom(B^*) \cap \ran(B) = \dom(|B^*|) \cap \ran (|B^*|) \, 
\text{ is dense in $\cK_B=\cK_A$} \\
& \quad \text{and a core for $|B^*|$ and $B^*$.}     \lb{3.42}
\end{split} 
\end{align}
Equation \eqref{3.42} together with \eqref{3.40} then yields 
 \begin{align}
 \begin{split}
 \ol{A^* A_0} &= \ol{B^* P_{\cK_A} A_0} = \ol{B^* B_0} = B^* \ol{B_0} 
 = B^* \wti B = B^* P_{\cK_A} \wti B = A^* \wti B = A^* \wti A   \lb{3.46} \\
 & = A^* \ol{A_0}. 
 \end{split}
 \end{align} 
Here we used Lemma \ref{l3.7}\,$(ii)$ (applied with $A$ replaced 
by $B$) in the third equality. This proves \eqref{3.41}. 
\end{proof}
%%%%%%%%%%%%%%

%%%%%%%%%%%%%%
\begin{remark} \lb{r3.8}
$(i)$ We note that  $H_\mathrm{en}$ used in Goldstein and Wacker \cite{GW03} 
coincides with $\cH_{|A|} = \cH_A$ used in the present paper. Moreover, it is 
noted in \cite[Proposition\ 2.1]{GW03} that $G_{|A|,0}$ (denoted by $\cA$ 
in \cite{GW03}) generates a strongly continuous unitary group. In addition, some 
properties of the domain $G_{|A|,0}$, amounting to the validity of \eqref{3.20} 
(with $A$ replaced by $|A|$ and hence also $A^*$ replaced by $|A|$), are mentioned without proof. The last part of Lemma \ref{l3.7} and of course 
Theorem \ref{t3.8} now explicitly provide such a proof. \\
$(ii)$ In connection with the operator $G_{A,0}$ in \eqref{3.12} and the 
second-order Cauchy problem (${\rm ACP}_2$) considered in the next 
Section \ref{s4}, we recall that 
$J_A^* = A^* \wti A = A^* \ol{A \iota_A} = A^* \ol{A_0} = \ol{A^* A_0} 
= \ol{A^* A \iota_A}$ as just shown in Theorem \ref{t3.8}. In the simpler 
situation where 
$A^* A \geq \varepsilon I_{\cH_1}$ for some $\varepsilon > 0$, one notes 
(cf.\ Lemma \ref{l2.7}) that $J^* = A^* \wti A = A^* A \iota_A$. 
The actual choice of $A$ in the factorization of the 
self-adjoint operator $S\geq 0$ into $S = A^* A$ is of course highly non-unique.   
 In particular, the self-adjoint factorization $S = S^{1/2} S^{1/2}$ (i.e., 
 $A = A^* = |A| = S^{1/2}$) is always possible, but may not be the most natural one as the following standard example shows. Let $S = -\Delta$ on 
 $\dom(S) = H^2(\bbR^n)$ be the usual self-adjoint Laplacian in 
 $\cH_1 = L^2\big(\bbR^n; d^n x\big)$, $n\in\bbN$ (with $H^m(\bbR^n)$, 
 $m\in\bbN$, the standard Sobolev spaces on $\bbR^n$). Then 
 \begin{equation}
 S = (-\Delta)^{1/2} (-\Delta)^{1/2} = \nabla^* \nabla,
 \end{equation} 
 with the last factorization being more natural for some purposes. Here 
 $\dom(\nabla) = H^1(\bbR^n)$, 
 $\cH_2 = \big[L^2\big(\bbR^n; d^n x\big)\big]^n$, and 
 $\nabla^* = - {\rm div} (\cdot)$ with 
 $\dom(\nabla^*) = \big[H^1(\bbR^n)\big]^n$. \\
 $(iii)$ We finally note that given the results \eqref{3.36}--\eqref{3.40} 
 and \eqref{3.42}, the result \eqref{3.19} in Lemma \ref{l3.7}\,$(i)$ can be 
 improved as follows and underscores the preliminary nature of the latter. 
 Since $\ran(A) = \ran(B)$ and 
 $A^* = B^* P_{\cK_A}$, and hence also 
\begin{equation} 
A^* = A^* P_{\cK_A} = B^* P_{\cK_A},     \lb{3.43} 
\end{equation} 
one concludes from \eqref{3.42} that 
\begin{align}
\begin{split} 
& \ol{\dom(A^* P_{\cK_A}) \cap \ran(A)} 
= \ol{\dom(B^* P_{\cK_A}) \cap \ran(B)}   \\
& \quad = \ol{\dom(B^*) \cap \ran(B)} = \cK_B = \cK_A.    \lb{3.44}
\end{split}
\end{align}
Thus, Lemma \ref{l3.7}\,$(i)$ applies and \eqref{3.19} can be amended to read 
\begin{equation}
\ol{A^* A_0} = \ol{A^* P_{\cK_A}\big|_{\dom(A^* P_{\cK_A}) \cap \ran(A)}} 
\, \ol{A_0}
= \ol{A^*\big|_{\dom(A^*) \cap \ran(A)}} 
\, \ol{A_0}.     \lb{3.45}
\end{equation}
Of course, Theorem \ref{t3.8} further improves on \eqref{3.45} and yields 
the final and optimal result \eqref{3.41}. 
 \end{remark}
%%%%%%%%%%%%%%

%%%%%%%%%%%%%%%%%%%%%%%%%%%%%%%%%%%%%%
%%%%%%%%%%%%%%%%%%%%%%%%%%%%%%%%%%%%%%
\section{Abstract Linear Damped Wave Equations. 
The Case $A^* A \geq \varepsilon I_{\cH}$ \\ for some $\varepsilon >0$}
\label{s4}
%%%%%%%%%%%%%%%%%%%%%%%%%%%%%%%%%%%%%%
%%%%%%%%%%%%%%%%%%%%%%%%%%%%%%%%%%%%%%

In this section we now introduce abstract damped wave equations employing 
appropriate perturbation techniques for Dirac-type operators. 

We first treat the case $A^* A \geq \varepsilon I_{\cH_1}$ for some $\varepsilon >0$   
 and hence introduce the following assumptions.

%%%%%%%%%%%
\begin{hypothesis} \lb{h4.1}
Let $\cH_j$, $j=1,2$, be complex separable Hilbert spaces. \\
$(i)$ Assume that $A: \dom(A) \subseteq \cH_1 \to \cH_2$ is a densely 
defined, closed, linear operator such that 
\begin{equation}
A^* A \geq \varepsilon I_{\cH_1}     \lb{4.2}
\end{equation} 
for some $\varepsilon >0$. \\
$(ii)$ Let $R$ be a densely defined, closable operator in $\cH_1$ satisfying
\begin{equation}
\dom (R) \supseteq \dom (A).    \lb{4.3}
\end{equation}
\end{hypothesis}
%%%%%%%%%%% 

We emphasize that closability of $R$ and the assumption \eqref{4.3} imply 
\begin{equation}
R (|A| - z I_{\cH_1})^{-1} \in \cB(\cH_1), \quad z \in \rho(|A|)   \lb{4.4}
\end{equation}
(see, e,g., \cite[p.\ 191]{Ka80}). 

In the following we intend to introduce the operator $i \, G_{A,R}$ in 
$\cH_A \oplus \cH_1$ and study its properties by utilizing its unitary 
equivalence to the Dirac-type operator $Q_{|A|,R}$ in $\cH_1 \oplus \cH_1$. 

We start by introducing $Q_{|A|,R}$ in $\cH_1 \oplus \cH_1$ assuming Hypothesis \ref{h4.1}.
\begin{equation}
Q_{|A|,R} = \begin{pmatrix} - i \, R & |A| \\ |A| & 0 \end{pmatrix}, \quad 
\dom (Q_{|A|,R}) = \dom(|A|) \oplus \dom(|A|) \subseteq \cH_1 \oplus \cH_1.   \lb{4.5} 
\end{equation} 

Next, we recall that an operator $T$ in the complex separable Hilbert 
space $\cH$ is called {\it accretive} if
\begin{equation}
\Re ((f, T f)_{\cH}) \geq 0, \quad f \in \dom(T).    \lb{4.5a}
\end{equation}
If in addition, $T$ is closed and $\ran(T + \alpha I_{\cH})$ is dense in $\cH$ 
for some $\alpha > 0$ , then $T$ is called $m$-accretive. Moreover 
(cf.\ \cite[Proposition\ C.7.2]{Ha06}, \cite[p.\ 279]{Ka80}) 
\begin{align} 
& \text{$T$ is $m$-accretive if and only if} \\
& \quad \text{$(-\infty,0) \subset \rho(T)$ and 
$\|(T - z I_{\cH})^{-1}\| \le -[\Re(z)]^{-1}$ for $\Re(z)<0$.}
\end{align}
In particular, an $m$-accretive operator is equivalent to a closed, densely 
defined, maximal accretive operator. Finally, $T$ is $m$-accretive if and 
only if $T^*$ is.(We note that one also calls $T$ ($m$-){\it dissipative} 
whenever $-T$ is ($m$-)accretive. However, since this definition is not 
universally accepted in the literature, we shall not adopt it here.) 

For the following it is convenient to introduce the quadratic operator pencil 
$M(z)$ in $\cH_1$, 
\begin{equation}
M(z) = |A|^2 - iz R - z^2 I_{\cH_1}, \quad \dom(M(z)) = \dom\big(|A|^2\big), 
\; z \in\bbC.     \lb{4.5b}
\end{equation}

%%%%%%%%%%%%
\begin{lemma} \lb{l4.2}
Assume Hypothesis \ref{h4.1}. \\
$(i)$ Then $M(z)$, $z\in\bbC$, is a densely defined, closed operator in $\cH_1$. \\
$(ii)$ If in addition $R^*$ satisfies 
\begin{equation}
\dom (R^*) \supseteq \dom (A),    \lb{4.5c}
\end{equation}
then  
\begin{equation}
M(z)^* = |A|^2 + i \ol{z} R^* - \ol{z}^2 I_{\cH_1}, \quad 
\dom(M(z)^*) = \dom\big(|A|^2\big), \; z \in\bbC.     \lb{4.5d} 
\end{equation}
\end{lemma}
%%%%%%%%%%%
\begin{proof}
Since by hypothesis $R$ is bounded with respect to $|A|$, it is relatively bounded with relative bound equal to zero with respect to $|A|^p$ for 
any $p>1$ (cf.\ \cite[Theorem\ 9.11\,(a)]{We80}). Thus, 
for each $z \in \bbC$, $M(z)$ is a closed operator in $\cH_1$ by a 
Kato--Rellich-type result (cf.\ \cite[Theorem\ IV.1.1]{Ka80}, \cite[Theorem\ 5.5]{We80}). 
Since by \eqref{4.5c} also $R^*$ is relatively bounded with relative bound 
equal to zero with respect to $|A|^p$ for any $p>1$, one also obtains \eqref{4.5d} 
by a Kato--Rellich-type argument discussed in \cite[p.\ 111]{We80}.  
\end{proof}
%%%%%%%%%%%

The spectrum and resolvent set of $M(\cdot)$, denoted by $\sigma(M(\cdot))$ and 
$\rho(M(\cdot))$, respectively, are then defined by
\begin{align}
\sigma(M(\cdot)) &= \{\lambda \in \bbC\,|\, 0 \in \sigma (M(\lambda))\},   \\
\rho(M(\cdot)) &= \{z \in\bbC \,|\, 0 \in \rho(M(z))\}  
= \{z \in\bbC \,|\, M(z)^{-1} \in \cB(\cH_1)\}   \no \\ 
& = \bbC \backslash \sigma(M(\cdot)). 
\end{align}

%%%%%%%%%%%%
\begin{theorem} \lb{t4.3}
Assume Hypothesis \ref{h4.1}. Then $Q_{|A|,R}$ is injective and closed, and  
\begin{align}
&(Q_{|A|,R} - z I_{\cH_1\oplus \cH_1})^{-1} =  
 \begin{pmatrix} z M(z)^{-1} 
& |A|^{-1} + M(z)^{-1} (i z R + z^2 I_{\cH_1}) |A|^{-1} \\
|A| M(z)^{-1} & |A| M(z)^{-1} (iR + z I_{\cH_1}) |A|^{-1} \end{pmatrix},   \no \\
& \hspace*{2.7cm} z \in \rho(Q_{|A|,R}) 
= \{\zeta \in\bbC \,|\, M(\zeta)^{-1} \in \cB(\cH_1)\} = \rho(M(\cdot)),    \lb{4.5e} \\
& \sigma(Q_{|A|,R}) = \sigma (M(\cdot)).    \lb{4.5f}
\end{align}
In particular, 
\begin{equation}
Q_{|A|,R}^{-1} = \begin{pmatrix} 0 & |A|^{-1} \\ |A|^{-1} & i \, |A|^{-1} R |A|^{-1} \end{pmatrix} 
\in \cB(\cH_1 \oplus \cH_1),     \lb{4.6}
\end{equation}
and hence 
\begin{equation}
0 \in \rho (Q_{|A|,R}).    \lb{4.7}
\end{equation}
Suppose in addition that $R$ is accretive. Then also $i \, Q_{|A|,R}$ is accretive 
and $- i  \, Q_{|A|,R}$ generates a contraction semigroup in $\cH_1 \oplus \cH_1$, denoted by $\exp (-i \, Q_{|A|,R} \, t)$, $t\geq 0$. 
\end{theorem}
%%%%%%%%%%%%
\begin{proof}
To prove injectivity of $Q_{|A|,R}$, assume $f, g \in \dom(|A|)$ and   
$Q_{|A|,R} (f \;\; g)^\top = 0$. Then 
$- i \, R f + |A| g = 0$ and $|A| f = 0$ imply $f = 0$ since $\ker(A) = \ker(|A|) = \{0\}$ by \eqref{4.2} and hence also $|A| g = 0$, implying $g = 0$ as well. That $Q_{|A|,R}$ is 
closed in $\cH_1 \oplus \cH_1$ follows from an application of 
\cite[Corollary\ 2.2.11\,$(ii)$]{Tr08}. (Alternatively, one can first establish \eqref{4.6} 
directly and then use again the fact that since $Q_{|A|,R}$ is injective, 
$Q_{|A|,R}^{-1} \in \cB(\cH_1 \oplus \cH_1)$ implies closedness of $Q_{|A|,R}^{-1}$ 
and hence that of $Q_{|A|,R}$ by \cite[p.\ 81]{We80}.) 

Denoting temporarily the right-hand side of \eqref{4.5e} by $S_{|A|,R}(z)$, one notes that 
$S_{|A|,R}(z) \in  \cB(\cH_1 \oplus \cH_1)$, $z \in \rho(M(\cdot))$, 
by \eqref{4.2} and \eqref{4.4}. A simple computation then yields that
\begin{align}
& (Q_{|A|,R} - z I_{\cH_1\oplus \cH_1}) \, S_{|A|,R}(z) = I_{\cH_1 \oplus \cH_1},  \no \\
& S_{|A|,R}(z) \, (Q_{|A|,R} - z I_{\cH_1\oplus \cH_1}) 
= I_{\cH_1 \oplus \cH_1}\big|_{\dom(|A|) \oplus \dom(|A|)},  \\ 
& \hspace*{2.35cm} z \in \{\zeta \in\bbC \,|\, M(\zeta)^{-1} \in \cB(\cH_1)\} = \rho(M(\cdot)),  \no 
\end{align} 
implying 
\begin{equation}
\rho(M(\cdot)) \subseteq \rho(Q_{|A|,R}).
\end{equation}
Conversely, let $z\in \rho(Q_{|A|,R})$. Then the resolvent of $Q_{|A|,R}$ is necessarily of 
the $2 \times 2$ block operator form with respect to $\cH_1 \oplus \cH_2$, 
\begin{equation}
(Q_{|A|,R} - z I_{\cH_1 \oplus \cH_2})^{-1} = 
\begin{pmatrix} S_{1,1} (z) & S_{1,2} (z) \\ S_{2,1}  (z) & S_{2,2} (z) \end{pmatrix}, 
\quad z \in \rho(Q_{|A|,R}),     \lb{4.8} 
\end{equation}
where $S_{j,k} (z) \in \cB(\cH_k,\cH_j)$, $j,k \in \{1,2\}$, $ z \in \rho(Q_{|A|,R})$. Thus,
\begin{align}
& \begin{pmatrix} I_{\cH_1} & 0 \\ 0 & I_{\cH_2} \end{pmatrix} = 
\begin{pmatrix} -i R - z I_{\cH_1} & |A| \\ |A| & - z I_{\cH_2} \end{pmatrix} 
\begin{pmatrix} S_{1,1} (z) & S_{1,2} (z) \\ S_{2,1}  (z) & S_{2,2} (z) \end{pmatrix}   \no \\
& \quad = \begin{pmatrix} (-i R - z I_{\cH_1}) S_{1,1} (z) + |A| S_{2,1} (z)
& (-i R - z I_{\cH_1}) S_{1,2} (z) + |A| S_{2,2} (z) \\ 
|A| S_{1,1} (z) - z S_{2,1}  (z) & |A| S_{1,2} (z) - z S_{2,2} (z) \end{pmatrix},   \no \\
& \hspace*{9cm}  z \in \rho(Q_{|A|,R}), 
\end{align}
in particular, 
\begin{equation}
S_{2,1} (z) = z^{-1} |A| S_{1,1} (z), \quad z \in \rho(Q_{|A|,R})\backslash \{0\}, 
\end{equation}
and hence 
\begin{equation}
z^{-1} \big(|A|^2 - i z R - z^2 I_{\cH_1}\big) S_{1,1} (z) = I_{\cH_1}, \quad 
z \in \rho(Q_{|A|,R})\backslash \{0\}.
\end{equation}
Thus, $z^{-1} S_{1,1} (z)$ is a bounded right-inverse of $M(z)$, 
$z \in \rho(Q_{|A|,R})\backslash \{0\}$. An analogous computation yields 
\begin{equation}
z^{-1} S_{1,1} (z) \big(|A|^2 - i z R - z^2 I_{\cH_1}\big)\big|_{\dom(|A|^2)} 
= I_{\cH_1}\big|_{\dom(|A|^2)}, \quad z \in \rho(Q_{|A|,R})\backslash \{0\}, 
\end{equation}
and hence $z^{-1} S_{1,1} (z)$ is also a bounded left-inverse of $M(z)$. Thus,
\begin{equation}
\rho(Q_{|A|,R})\backslash \{0\} \subseteq \rho(M(\cdot)).    \lb{4.8a}
\end{equation}
Since by hypothesis \eqref{4.2}, $0 \in \rho(Q_{|A|,R}) \cap \rho(M(\cdot))$, one concludes 
$\rho(Q_{|A|,R}) = \rho(M(\cdot))$ and hence \eqref{4.5e}--\eqref{4.7}. 

Finally, assuming $f, g \in \dom(|A|)$ one computes
\begin{equation}
\Re \big(\big(( f \; g)^\top, i \, Q_{|A|,R} \, (f \; g)^\top 
\big)_{\cH_1 \oplus \cH_1}\big) = \Re \big((f, R f)_{\cH_1} \big) \geq 0, 
\end{equation}
since $|A|$ is self-adjoint. This proves that $i \, Q_{|A|,R}$ is accretive. 
Equation \eqref{4.7} yields a sufficiently small open disk with center at $0$ in the 
resolvent set of $- i \, Q_{|A|,R}$ and this fact combined with the Lumer--Phillips 
Theorem \cite{LP61} (cf.\ \cite[Theorem\ II.3.15]{EN00}) then 
proves that $- i  \, Q_{|A|,R}$ generates a contraction semigroup.
\end{proof}
%%%%%%%%%%%

We note that block operator matrices and their inverses, and more specifically, spectral properties of $2 \times 2$ block operator matrices  have been studied extensively in the literature. We refer, for instance, to \cite{BBDHL05}, 
\cite{En99}, \cite[Sect.\ VI.6]{EN00}, \cite{HKM00}, 
\cite{KMN02}, \cite{MS96}, \cite{Na89}, \cite{Na90}, \cite{St10}, 
\cite[Ch.\ 2]{Tr08}, \cite{Tr09}, \cite{Ve04}, and \cite{Wr86}. 

Still assuming Hypothesis \ref{h4.1}, we next introduce the operator $G_{A,R}$ in $\cH_A \oplus \cH_1$ by
\begin{align}
\begin{split} 
& G_{A,R} = \begin{pmatrix} 0 & J_A \\ - J_A^* & - R \end{pmatrix} 
= \begin{pmatrix} 0 & J_A \\ - A^* \wti A & - R \end{pmatrix},    \\
& \dom(G_{A,R}) = \dom \big(A^* \wti A \big) \oplus \dom(A) 
\subseteq \cH_A \oplus \cH_1.     \lb{4.9}
\end{split}
\end{align}
In particular, one notes that
\begin{equation}
G_{A,R} = G_{|A|,R}. 
\end{equation}

%%%%%%%%%%%%%%
\begin{theorem} \lb{t4.4}
Assume Hypothesis \ref{h4.1}. Then
\begin{equation}
Q_{|A|,R} =  U_{\wti{|A|}} \, i \, G_{A,R} U_{\wti{|A|}}^{-1}.    
\lb{4.10}
\end{equation}
In particular, the operator $G_{A,R}$ is densely defined and closed in the 
energy space 
$\cH_A \oplus \cH_1$. 

If in addition $R$ is accretive, then also $- G_{A,R}$ is accretive and 
$G_{A,R}$ generates a contraction semigroup 
in $\cH_A \oplus \cH_1$, denoted by $\exp (G_{A,R} \, t)$, $t\geq 0$.
\end{theorem}
%%%%%%%%%%%%%%
\begin{proof}
To prove \eqref{4.10}, one can closely follow the proof of Theorem \ref{t2.8} in the special case $R=0$. In particular, since we will use $Q_{|A|,R}$ (instead of 
$Q_{A,R}$) this permits us to replace the pair $(A, A^*)$ by $(|A|, |A|)$ and hence replace the projection $P_{\ker(Q_{A,0})}$ by $0$ (cf.\ \eqref{2.70}, \eqref{2.71}). Alternatively, one can also argue as follows (cf.\ \eqref{2.65}, \eqref{2.66}, and 
\eqref{2.68}--\eqref{2.70}).
\begin{align}
U_{\wti{|A|}}^{-1} \dom (Q_{|A|,R}) &= \begin{pmatrix} 0 & i \, \big(\wti{|A|}\big)^{-1} \\ 
I_{\cH_1} & 0 \end{pmatrix} \dom (|A|) \oplus \dom (|A|)   \no \\
& = \begin{pmatrix} i \, \big(\wti{|A|}\big)^{-1} \dom(|A|) \\ \dom (|A|) \end{pmatrix}     
= \begin{pmatrix} i \, J_A |A|^{-1} \dom(|A|) \\ \dom (|A|) \end{pmatrix}    \no \\
& = \begin{pmatrix} i \, J_A (A^* A)^{-1/2} \dom \big((A^* A)^{1/2}\big) \\ \dom (A) \end{pmatrix}  \no \\
& = \begin{pmatrix} i \, J_A \dom (A^* A) \\ \dom (A) \end{pmatrix}    
= \begin{pmatrix} i \, \dom (A^* A \iota_A) \\ \dom (A) \end{pmatrix}    
= \begin{pmatrix} i \, \dom \big(A^* \wti A \big) \\ \dom (A) \end{pmatrix}    \no \\
& = \dom (G_{A,R}),    \lb{4.13}
\end{align}
where we used \eqref{2.26B} in the next to last step. Analogously to \eqref{2.35} one then obtains 
\begin{align}
U_{\wti{|A|}} i \, G_{A,R} U_{\wti{|A|}}^{-1} 
 &= i \begin{pmatrix} 0 & I_{\cH_1} \\ -i \, \wti{|A|} & 0 \end{pmatrix} 
\begin{pmatrix} 0 & J_A \\ - A^* \wti A & - R \end{pmatrix} 
\begin{pmatrix} 0 & i \big(\wti{|A|}\big)^{-1} \\ I_{\cH_1} & 0 \end{pmatrix}  \no \\
& = i \begin{pmatrix} 0 & I_{\cH_1} \\ -i \, \wti{|A|} & 0 \end{pmatrix} 
\begin{pmatrix} J_A & 0 \\ - R & -i \, |A|^2\iota_A \big(\wti{|A|}\big)^{-1} \end{pmatrix}   \no \\
& = i \begin{pmatrix} 0 & I_{\cH_1} \\ - i \, \wti{|A|} & 0 \end{pmatrix}  
 \begin{pmatrix} J_A & 0 \\ - R & -i \, |A| \end{pmatrix}   \no \\
& = \begin{pmatrix} - i \, R & |A| \\  \wti{|A|} J_A & 0 \end{pmatrix}     
= \begin{pmatrix} - i \, R & |A| \\ |A| & 0 \end{pmatrix}    \no \\ 
& = Q_{|A|,R},    \lb{4.14}
\end{align}
using \eqref{2.67a} in the next to last step. 

Closedness of $G_{A,R}$ in $\cH_A \oplus \cH_1$ then follows from \eqref{4.14} and that 
of $Q_{|A|,R}$ (cf.\ Theorem \ref{t4.3}). Similarly, if $R$ is accretive, then the contraction semigroup 
property of $G_{A,R}$ follows from the one of $- i \, Q_{|A|,R}$ in Theorem \ref{t4.3}, using \eqref{4.14} again.
\end{proof}
%%%%%%%%%%%%%%

We note again that the unitary equivalence in \eqref{4.10} has been observed by Huang \cite[Proposition\ 3.1]{Hu97}. While the contraction semigroup result for 
$G_{A,R}$ in Theorem \ref{t4.4} is well-known (see, e.g., \cite[Sect.\ VI.3]{EN00}), 
we presented it in some detail to illustrate the usefulness of the unitary equivalence relation with the Dirac-type operator $Q_{|A|,R}$ which leads to a rather simple 
proof.   

We also mention the analog of the result \eqref{4.10} when using $A$ and $A^*$ in place of $|A|$. Introducing the operator $Q_{A,R}$ in $\cH_1 \oplus \cH_2$ by  
\begin{equation}
Q_{A,R} = \begin{pmatrix} - i \, R & A^* \\ A & 0 \end{pmatrix}, \quad 
\dom (Q_{A,R}) = \dom(A) \oplus \dom(A^*) \subseteq \cH_1 \oplus \cH_2   
\lb{4.14A} 
\end{equation} 
one obtains the following result.

%%%%%%%%%%%%%%
\begin{theorem} \lb{t4.5}
Assume Hypothesis \ref{h4.1}. Then
\begin{equation}
Q_{A,R} \big(I_{\cH_1} \oplus [I_{\cH_2} - P_{\ker(A^*)}]\big) 
=  U_{\wti A} \, i \, G_{A,R} U_{\wti A}^{-1}.    \lb{4.14a}
\end{equation}
\end{theorem}
%%%%%%%%%%%%%%
\begin{proof}
It suffices to combine Theorem \ref{t2.8}, \eqref{2.29}, \eqref{A.21}, and 
\begin{equation}
U_{\wti A} \begin{pmatrix} 0 & 0 \\ 0 & -i R \end{pmatrix} U_{\wti A}^{-1} 
=  \begin{pmatrix} -i R & 0 \\ 0 & 0 \end{pmatrix}.     \lb{4.14b} 
\end{equation}
\end{proof}
%%%%%%%%%%%%%%

Theorem \ref{t4.5} appears to be a new result. 

Next, we briefly recall the notion of classical solutions of first-order and second-order Cauchy problems.  First, let $G$ be a densely defined, closed, linear operator in a Hilbert space $\cH$. Then the abstract Cauchy Problem (${\rm ACP}$) associated with $G$, by definition, is the initial value problem 
(cf., e.g., \cite[Sect.\ II.6]{EN00}, \cite[Ch.\ I]{Fa85}, \cite[Sect.\ II.1]{Go85})
\begin{align}
\begin{split}
&\dot x (t) = G x (t), \; t\geq 0,  \\
& x (0) = x_0 \in \cH.    \lb{4.15}
\end{split}
\end{align} 
Here we denote $\dot x (t) = \big(\tfrac{d}{dt} x\big)(t)$. 

By definition, a {\it classical solution} of the ${\rm ACP}$ \eqref{4.15} is then a map $x: [0,\infty) \to \cH$ which satisfies.
\begin{align}
& x \in C^1 ([0,\infty); \cH),   \no \\
& x (t) \in \dom (G), \; t \geq 0,   \lb{4.16}\\
& \text{$x$ satisfies \eqref{4.15}.}   \no
\end{align}
In particular, if $G$ is the generator of a strongly continuous semigroup 
$T(t) = e^{Gt}$, $t \geq 0$, 
in $\cH$, then for any $x_0 \in \dom (G)$, the unique classical solution 
$x = x(t)$ of \eqref{4.15} 
is given by 
\begin{equation}
x (t) = e^{Gt} x_0, \; t\geq 0, \quad x_0 \in \dom(G).   \lb{4.17}
\end{equation}
Moreover, the classical solution of \eqref{4.15} exists if and only if 
$x_0 \in \dom (G)$. 

Similarly, let $R$ and $S$ be densely defined, closed, linear operators in $\cH$. Then the abstract second-order Cauchy Problem (${\rm ACP}_2$) associated with $R$ and $S$ is by definition the initial value problem (cf., e.g., \cite{CP89}, 
\cite{En92}, \cite{En94}, \cite[Sect.\ VI.3]{EN00}, \cite[Chs.\ II, III, VIII]{Fa85}, 
\cite[Sect.\ II.7]{Go85}, \cite{XL98}), 
\begin{align}
\begin{split}
&\ddot{u} (t) + R \dot u (t) + S u (t) =0, \; t\geq 0,  \\
& u(0) = u_0 \in \cH, \; \dot u (0) = u_1 \in \cH.    \lb{4.18}
\end{split}
\end{align} 

By definition, a {\it classical solution} of the ${\rm ACP}_2$ \eqref{4.18} is then a map $u: [0,\infty) \to \cH$ which satisfies.
\begin{align}
\begin{split}
& u \in C^2 ([0,\infty); \cH),   \\
& u(t) \in \dom (S), \; t \geq 0, \, \text{ and } \, Su \in C([0,\infty); \cH),   \lb{4.19} \\
& \dot u (t) \in \dom (R), \; t\geq 0, \, \text{ and } \, R \dot u \in C([0,\infty); \cH),  \\
& \text{$u$ satisfies \eqref{4.18}.}  
\end{split}
\end{align}

%%%%%%%%%%%%%%
\begin{corollary} \lb{c4.6} 
Assume Hypothesis \ref{h4.1} and suppose that $R$ is accretive. In addition, let 
\begin{equation}
\begin{pmatrix} y_0 \\  z_0 \end{pmatrix} \in \dom (G_{A,R}) = 
\dom (A^* A \iota_A) \oplus \dom(A) \subseteq \cH_A \oplus \cH_1.   \lb{4.20}
\end{equation}
Then 
\begin{equation}
Y (t) = \begin{pmatrix} y (t) \\ z (t) \end{pmatrix} = e^{G_{A,R}t} 
\begin{pmatrix} y_0 \\  z_0 \end{pmatrix}, \; t\geq 0,    \lb{4.21}
\end{equation}
is the unique classical solution of the ${\rm ACP}$,  
\begin{align}
\begin{split}
&\dot Y (t) = G_{A,R} Y (t), \; t\geq 0,  \\
& Y (0) = \begin{pmatrix} y_0 \\  z_0 \end{pmatrix},    \lb{4.22}
\end{split}
\end{align} 
associated with $G_{A,R}$, and 
\begin{equation}
u (t) = \iota_A y (t), \; t\geq 0,    \lb{4.23}
\end{equation}
is the unique classical solution of the ${\rm ACP}_2$,  
\begin{align}
\begin{split}
&\ddot{u} (t) + R \dot u (t) + A^* A u (t) =0, \; t\geq 0,  \\
& u(0) = \iota_A y_0, \; \dot u (0) = z_0,     \lb{4.24}
\end{split}
\end{align} 
associated with $R$ and $A^* A$.
\end{corollary}
%%%%%%%%%%%%%%
\begin{proof}
One only needs to verify \eqref{4.23}, \eqref{4.24}, and uniqueness of $u$. 
From \eqref{4.21} one infers 
\begin{align}
\begin{split}
& y(\cdot)\in C^1([0,\infty); \cH_A), \quad y(t) \in \dom(J_A^*) = 
\dom(A^* A \iota_A),  \\
& z(t) \in \dom(J_A) = \dom(A), \; t \geq 0, \, \text{ with } \, \dot y = J_A z, 
\end{split} 
\end{align}
and 
\begin{equation}
z(\cdot)\in C^1([0,\infty); \cH_1) \, \text{ with } \, \dot z = - A^* A \iota_A y - R z. 
\end{equation}
Hence by \eqref{4.23} one obtains 
\begin{equation}
u(\cdot) = \iota_A y(\cdot) \in C^1([0,\infty); \cH_1) 
\, \text{ with } \, \dot u = \iota_A \dot y = \iota_A J _A z =z,
\end{equation}
since $\iota_A \in \cB(\cH_A, \cH_1)$. Moreover, this shows that 
$\dot u(\cdot) \in C^1([0,\infty); \cH_1)$, 
implying 
\begin{equation} 
u(\cdot) \in C^2([0,\infty); \cH_1) \, \text{ and } \, 
\ddot u = \dot z = - A^* A \iota_A y - R z = -A^* A u -R \dot u. 
\end{equation}
Consequently, $u (\cdot) = \iota_A y (\cdot)$ is a classical solution of \eqref{4.24}.

Finally, uniqueness of $u(\cdot)$ is shown as in \cite[Prop.\ VI.3.2]{EN00}.
First of all one notes that $\wti{A} y(\cdot) = A u(\cdot) \in C^1([0,\infty); \cH_2)$, 
implying
$\wti{A} \dot y(\cdot) = A \dot u(\cdot) \in C([0,\infty); \cH_2)$. Hence one has
\begin{equation}
R \dot u (\cdot) = \big[R (|A| + I_{\cH_1})^{-1}\big] (|A| + I_{\cH_1}) \dot u(\cdot) 
\in  C([0,\infty); \cH_1), 
\lb{4.31} 
\end{equation}
which in turn implies
\begin{equation}
A^* A u (\cdot) = - \ddot u (\cdot) - R \dot u (\cdot) \in C([0,\infty); \cH_1).      \lb{4.32}
\end{equation}
Now suppose that $u(\cdot)$ satisfies \eqref{4.24} with $(y_0,z_0)=(0,0)$. Then
\begin{equation}
\int_0^t ds \, \begin{pmatrix} J_A u (s) \\ \dot u (s)\end{pmatrix} 
= \begin{pmatrix} J_A \int_0^t ds \, u(s) \\ u(t) \end{pmatrix} \in \dom (A^* A \iota_A) \oplus \dom (A), 
\; t \geq 0,     \lb{4.34}
\end{equation}
and
\begin{align}
G_{A,R} \begin{pmatrix} J_A \int_0^t ds \, u(s) \\ u(t) \end{pmatrix}    
& = \begin{pmatrix} J_A u (t) \\ - A^* A \int_0^t ds \, u(s) - R u(t) \end{pmatrix}    \no \\\
& =  \begin{pmatrix} J_A u (t) \\ - \int_0^t ds \, A^* A u(s) - \int_0^t ds \, R \dot u(s) \end{pmatrix}  \no \\
& = \begin{pmatrix} J_A u (t) \\ \dot u (t)\end{pmatrix}, \; t \geq 0.     \lb{4.35}
\end{align}
Thus, $\left(\begin{smallmatrix} J_A u (\cdot) \\ \dot u (\cdot)\end{smallmatrix}\right)$ is a mild solution 
(cf., e.g., \cite[Ch.\ II.6]{EN00}) 
of $\dot Y (\cdot) = G_{A,R} Y(\cdot)$ satisfying $Y(0) =0$. But then, $Y(t) = 0$, $t \geq 0$ 
(cf.\ \cite[Proposition\ VI.3.2]{EN00}), and hence $u(t) =0$, $t \geq 0$ (cf.\ \cite[Prop.\ II.6.4]{EN00}). 
\end{proof}
%%%%%%%%%%%%%%

Again, Corollary \ref{c4.6} is well-known (see, e.g., \cite[Sect.\ VI.3]{EN00}); for completeness, and due to its importance, we presented its proof in some detail.

%%%%%%%%%%%%%%
\begin{remark} \lb{r4.7}
In the special case $R=0$, and assuming Hypothesis \ref{h4.1}, the self-adjointness 
of $i \, G_{A,0}$ 
in $\cH_A \oplus \cH_1$ then yields the unitary group $e^{G_{A,0} t}$, $t\in\bbR$, in 
$\cH_A \oplus \cH_1$. Explicitly, using 
\begin{equation}
J_A^* J_A = A^* A, \quad J_A J_A^* = J_A A^* A \, \iota_A = \iota_A^{-1} A^* A \, \iota_A, 
\end{equation}
$ e^{G_{A,0} t}$ in the energy space $\cH_A \oplus \cH_1$, is of the form
\begin{align}
\begin{split}
e^{G_{A,0} t} &= \begin{pmatrix} J_A & 0 \\ 0 & I_{\cH_1} \end{pmatrix} 
\begin{pmatrix} \cos(|A| t) & |A|^{-1} | \sin(|A|t) \\ 
-|A| \sin(|A| t) & \cos(|A| t) \end{pmatrix} 
\begin{pmatrix} \iota_A & 0 \\ 0  & I_{\cH_1} \end{pmatrix}   \\
& = \begin{pmatrix} J_A \cos(|A| t) \, \iota_A & J_A |A|^{-1} | \sin(|A|t) \\ 
-|A| \sin(|A| t) \, \iota_A & \cos(|A| t) \end{pmatrix}, \quad t\in \bbR,     \lb{4.37}
\end{split} 
\end{align}
using the spectral theorem to define appropriate functions of the self-adjoint 
nonnegative operator $|A| = (A^* A)^{1/2}$ in $\cH_1$. By \eqref{4.10}, 
$e^{G_{A,0} t}$ in $\cH_A \oplus \cH_1$ is unitarily equivalent to 
$e^{-i Q_{|A|,0}}$ in $\cH_1 \oplus \cH_1$ and explicitly given by
\begin{align}
& e^{-i Q_{|A|,0} t} = \begin{pmatrix} 
\cos (|A| t) & - i \sin (|A| t) \\ 
- i \sin (|A| t) & \cos (|A| t) \end{pmatrix}, \quad t \in \bbR. 
\end{align}
\end{remark}
%%%%%%%%%%%%%%%%

While we primarily focused on (dissipative) damping operators satisfying 
Hypothesis \ref{h4.1}\,$(ii)$, we emphasize that a variety of different conditions on $R$ have also been studied in the literature. We refer, for instance, to \cite{BI88}, 
\cite{BIW95}, \cite{BE04}, \cite{CLL98}, \cite{CT89}, \cite{En99}, 
\cite[Sect.\ VI.3]{EN00}, \cite[Ch.\ VIII]{Fa85}, \cite{GS03}, \cite{HS99}, \cite{HS04}, 
\cite{Hu88}, \cite{Hu96}, \cite{Hu97}, \cite{JT07}, \cite{JT09}, \cite{JTW08}, 
\cite{TW03}, \cite{Ve04}, \cite{Ve08}, \cite{WT03}. 

%%%%%%%%%%%%%%%%

We continue with an illustrative example in which $R$ is assumed to commute with  $|A|$. 

\begin{example} \lb{e4.8}
Assume Hypothesis \ref{h4.1} and consider the special case where $R$ is an 
appropriate function of $|A|$, that is, $R = 2 F(|A|) \geq 0$. Abbreviating 
\begin{equation}
\Gamma (|A|) = \big[|A|^2 - F(|A|)^2 \big]^{1/2}    \lb{4.41}
\end{equation}
one then obtains 
\begin{align}
& e^{G_{A,R} t}     
= \begin{pmatrix} J_A & 0 \\ 0 & I_{\cH_1} \end{pmatrix} e^{- F(|A|) t}   \no \\
& \qquad \times  
 \begin{pmatrix} \big[\cos(\Gamma (|A|) t) & \Gamma (|A|)^{-1} \sin(\Gamma (|A|) t) \\
 + F(|A|) \Gamma (|A|)^{-1} \sin(\Gamma (|A|) t)\big]   \\[2mm] 
-|A|^2 \Gamma (|A|)^{-1} \sin(\Gamma (|A|) t) 
& \big[\cos(\Gamma (|A|) t)  \\
&- F(|A|) \Gamma (|A|)^{-1} \sin(\Gamma (|A|) t)\big] \end{pmatrix} 
\no \\
& \qquad \times 
\begin{pmatrix} \iota_A & 0 \\ 0  & I_{\cH_1} \end{pmatrix}  \no \\
& \quad = e^{- F(|A|) t}    \no \\ 
& \qquad \times \begin{pmatrix} J_A
\big [\cos(\Gamma (|A|) t)   & J_A \Gamma (|A|)^{-1} \sin(\Gamma (|A|) t) \\
+ F(|A|) \Gamma (|A|)^{-1} \sin(\Gamma (|A|) t)\big] \iota_A   \\[2mm]
-|A|^2 \Gamma (|A|)^{-1} \sin(\Gamma (|A|) t) \iota_A & 
\big[\cos(\Gamma (|A|) t)  \\
& - F(|A|) \Gamma (|A|)^{-1} \sin(\Gamma (|A|) t)\big] 
\end{pmatrix},     \no \\ 
& \hspace*{10cm} t \geq 0.     \lb{4.42}
\end{align}
We note that $\cos(\Gamma(|A|) t)$ and 
$\Gamma(|A|)^{-1} \sin(\Gamma(|A|) t)$ are in fact functions of $\Gamma(|A|)^2$ 
and hence the precise specification of the square root branch in \eqref{4.41} does 
not enter in \eqref{4.42}. In addition, using the spectral theorem for $\Gamma(|A|)$, 
one obtains that 
\begin{equation} 
\Gamma(|A|)^{-1} \sin(\Gamma(|A|) t) = f(\Gamma(|A|)) 
= \int_{[0,\infty)} f\Big(\big[\lambda^2 - F(\lambda)^2\big]^{1/2}\Big)  
\, dE_{|A|}(\lambda) \in\cB(\cH_1) 
\end{equation} 
is well-defined without assuming that $\Gamma(|A|)$ is 
boundedly invertible in $\cH_1$ by choosing 
\begin{equation} 
f(\mu) = \begin{cases} \sin(\mu t)/\mu, & \mu > 0,   \\
t, & \mu =0. \end{cases}  
\end{equation}

Next, we intend to exploit the unitary equivalence between $G_{A,R}$ in 
$\cH_A \oplus \cH_1$ 
and the Dirac-type operator $Q_{|A|,R}$ in $\cH_1 \oplus \cH_1$ in \eqref{4.10} 
and thus we now turn to $e^{-i t Q_{|A|,R}}$. Noticing that
\begin{align}
e^{- i \begin{pmatrix} - 2 i a & 1 \\ 1 & 0 \end{pmatrix} t} 
& = e^{- a t} e^{- i \begin{pmatrix} - i a & 1 \\ 1 &  i a \end{pmatrix} t}  \no \\
& = e^{- a t} \begin{pmatrix} \cos(b t) - a b^{-1} \sin(b t) & -i b^{-1} \sin(b t) \\ 
- i b^{-1} \sin(b t) & \cos(b t) + a b^{-1} \sin(b t) \end{pmatrix},    \\[1mm]
& \hspace*{4.1cm} a \in\bbR, \; b = (1 - a^2)^{1/2}, \; t \geq 0,    \no 
\end{align}  
one obtains
\begin{align}
& e^{-i Q_{|A|,R} t} = e^{-F(|A|) t}   \no \\
& \qquad \times \begin{pmatrix} 
\cos (\Gamma (|A|) t) & - i |A| \Gamma(|A|)^{-1} \sin (\Gamma(|A|) t) \\ 
- F(|A|) \Gamma(|A|)^{-1} \sin (\Gamma(|A|) t) \\[2mm] 
- i |A| \Gamma(|A|)^{-1} \sin (\Gamma(|A|) t) 
& \cos (\Gamma(|A|) t) \\
& + F(|A|) \Gamma(|A|)^{-1} \sin (\Gamma(|A|) t) \end{pmatrix}  \no \\
& \quad = \int_{[0,\infty)} e^{-F(\lambda) t}   \no \\
& \qquad \times \begin{pmatrix} 
\cos (\Gamma (\lambda) t) & - i \lambda \Gamma(\lambda)^{-1} 
\sin (\Gamma(\lambda) t) \\ 
- F(\lambda) \Gamma(\lambda)^{-1} \sin (\Gamma(\lambda) t) \\[2mm] 
- i \lambda \Gamma(\lambda)^{-1} \sin (\Gamma(\lambda) t) 
& \cos (\Gamma(\lambda) t) \\
& + F(\lambda) \Gamma(\lambda)^{-1} \sin (\Gamma(\lambda) t) \end{pmatrix} 
d E_{|A|}(\lambda),  \no \\
& \hspace*{9.2cm} \quad t \geq 0. 
\end{align}

For the norm of the semigroup of $e^{G_{A,R}t}$, $t\geq 0$, one thus obtains 
$($cf.\ \cite[Sect.\ 5.3]{BS87}$)$ 
\begin{align} 
& \big\|e^{G_{A, 2 F(|A|)} t} \big\|_{\cB(\cH_A \oplus \cH_1)} 
= \big\|e^{-i Q_{|A|,R} t}\big\|_{\cH_1 \oplus \cH_1}    \no \\
& \quad = E_{|A|}\mhyphen\esssup_{\lambda \in \sigma(|A|)} 
\, e^{-F(\lambda) t}   \no \\  
& \qquad \times \left\|
\begin{pmatrix} 
\cos (\Gamma (\lambda) t) & - i \lambda \Gamma(\lambda)^{-1} 
\sin (\Gamma(\lambda) t) \\ 
- F(\lambda) \Gamma(\lambda)^{-1} \sin (\Gamma(\lambda) t) \\[2mm] 
- i \lambda \Gamma(\lambda)^{-1} \sin (\Gamma(\lambda) t) 
& \cos (\Gamma(\lambda) t) \\
& + F(\lambda) \Gamma(\lambda)^{-1} \sin (\Gamma(\lambda) t) \end{pmatrix} 
\right\|_{\bbC^2}.    \label{4.50} 
\end{align}

Denoting temporarily the  $2 \times 2$ matrix under the norm in \eqref{4.50} 
by $M(\lambda)$, $\lambda \in \sigma(|A|)$, to compute the norm of $M(\lambda)$ 
one computes the square root of the larger of the two eigenvalues of 
$M(\lambda)^* M(\lambda)$, that is, 
\begin{equation}
\sigma (M(\lambda)^* M(\lambda)) = \big\{s_1(M(\lambda))^2, s_2(M(\lambda))^2\big\}, \quad 
0 \leq s_1(M(\lambda)) \leq s_2(M(\lambda)),
\end{equation}
with $s_j(M(\lambda))$, $j=1,2$, the singular values of $M(\lambda)$, and 
\begin{equation}
\|M(\lambda)\|_{\bbC^2} = s_2(M(\lambda)).     \lb{4.52}
\end{equation}
An explicit computation yields 
\begin{align} 
s_j (M(\lambda)) & = 
\Bigg\{\Bigg[1 + \f{F(\lambda)^2 \sin^2\big([\lambda^2 
- F(\lambda)^2]^{1/2} t\big)}{[\lambda^2 - F(\lambda)^2]}\Bigg]^{1/2}   \lb{4.53} \\
& \qquad + (-1)^j  \f{F(\lambda) \big|\sin^2\big([\lambda^2 
- F(\lambda)^2]^{1/2} t\big)\big|}{[\lambda^2 - F(\lambda)^2]^{1/2}}\Bigg\}, \quad 
j=1,2, \quad \lambda \in \sigma(|A|).     \no 
\end{align}
Combining \eqref{4.50}, \eqref{4.52}, and \eqref{4.53}, one finally obtains  
\begin{align} 
& \big\|e^{G_{A, 2 F(A)} t} \big\|_{\cB(\cH_A \oplus \cH_1)} 
= \big\|e^{-i Q_{|A|,R} t}\big\|_{\cH_1 \oplus \cH_1}    \no \\
& \quad = E_{|A|}\mhyphen\esssup_{\lambda \in \sigma(|A|)} \, e^{-F(\lambda) t}   
\Bigg\{\Bigg[1 + \f{F(\lambda)^2 \sin^2\big([\lambda^2 
- F(\lambda)^2]^{1/2} t\big)}{[\lambda^2 - F(\lambda)^2]}\Bigg]^{1/2}  \no \\[1mm]
& \hspace*{4.2cm} + \f{F(\lambda) \big|\sin^2\big([\lambda^2 
- F(\lambda)^2]^{1/2} t\big)\big|}{[\lambda^2 - F(\lambda)^2]^{1/2}}\Bigg\}    
\lb{4.54} \\
& \quad \leq C \begin{cases} e^{- \omega (G_{A, 2 F(A)}) t}, & 
E_{|A|}\big(\big\{\lambda \in \sigma(|A|) \, \big| \, 
F(\lambda)^2 = \lambda^2\big\}\big) = 0, \\
t e^{-\omega (G_{A, 2 F(A)}) t}, & 
E_{|A|}\big(\big\{\lambda \in \sigma(|A|) \, \big| \, 
F(\lambda)^2 = \lambda^2\big\}\big) > 0,
\end{cases} \quad t \geq 0,      \label{4.55} 
\end{align}
where
\begin{equation}
\omega (G_{A, 2 F(A)}) = \text{$E_{|A|}$-$\essinf$}_{\lambda \in \sigma(|A|)} 
\Big[ F(\lambda) - \big[(F(\lambda)^2 - \lambda^2)_+ \big]^{1/2}\Big],      \lb{4.56}
\end{equation}
and $C\geq 0$ is an appropriate constant. Here 
\begin{equation}
x_+ =\begin{cases} x, & x \geq 0, \\ 0, & x \leq 0. \end{cases}  
\end{equation}
The projection operator-valued measure $E_{|A|}(\cdot)$ in \eqref{4.50}, 
\eqref{4.54}--\eqref{4.56} can be replaced by an equivalent scalar control measure 
$\rho_{|A|}(\cdot)$. For instance, one can choose 
\begin{equation}
d \rho_{|A|} (\cdot) = d \|E_{|A|} (\cdot) g_1\|^2_{\cH_1},   \quad 
g_1 = \sum_{j\in\bbN} 2^{-j} \psi_{1,j},
\end{equation}
with $\{\psi_{1,j}\}_{j\in\bbN}$ a complete orthonormal system in $\cH_1$. 

In particular, $- \omega (G_{A, 2 F(|A|)})$ represents the semigroup growth 
bound $($or type$)$ of $e^{G_{A, 2 F(|A|)} t}$, $t\geq 0$ $($cf., e.g., 
\cite[Definition\ I.5.6]{EN00}$)$. Moreover, alluding to the spectral theorem 
for $Q_{|A|,2 F(|A|)}$, $- \omega (G_{A, 2 F(|A|)})$ coincides with the spectral 
bound $($cf., e.g., \cite[Definition\ II.1.12]{EN00}$)$ of $G_{A, 2 F(|A|)}$ 
$($and hence that of $- i Q_{|A|,2 F(|A|)}$$)$. 
\end{example}
%%%%%%%%%%%%%%%%

In this commutative context we also refer to \cite{Ho91} where matrix multiplication operators generating one-parameter semigroups are studied. 

%%%%%%%%%%%%%%%%
\begin{remark} \lb{r4.9} 
We note that the special example where $R = 2 F(A) = |A|^{\alpha}$, 
$\alpha \in [0,1]$, has been discussed in \cite{BT10} and \cite{Hu97}, and in the 
case $\alpha < 0$ in \cite{PBR04}. The case $\alpha \in \bbR$ is studied in 
\cite{Fa90}. In particular, Huang 
\cite[Corollary\ 3.6]{Hu97} estimated the 
semigroup growth bound for $e^{G_{A,R}t}$, $t\geq 0$, from above using a combination of Gearhart's theorem, the unitary equivalence \eqref{4.10}, and 
certain norm estimates. Since he does not rely on the spectral theorem, his bound 
differs from the exact result in \eqref{4.56} in the case $R = |A|^{\alpha}$, 
$\alpha \in [0,1]$. On the other hand, his technique also yields an upper bound in cases where $R$ and $|A|$ do not commute. 
\end{remark}
%%%%%%%%%%%%%%%%

%%%%%%%%%%%%%%%%%%%%%%%%%%%%%%%%%%%%%%
%%%%%%%%%%%%%%%%%%%%%%%%%%%%%%%%%%%%%%
\section{Abstract Linear Damped Wave Equations. 
The Case $\inf(\sigma(A^* A)) = 0$}
\label{s5}
%%%%%%%%%%%%%%%%%%%%%%%%%%%%%%%%%%%%%%
%%%%%%%%%%%%%%%%%%%%%%%%%%%%%%%%%%%%%%

The principal aim of this section is to relax Hypothesis \ref{h4.1} and remove the hypothesis that $A^*A$ is strictly positive definite. 

Our basic hypothesis for this section reads as follows. 

%%%%%%%%%%%
\begin{hypothesis} \lb{h5.1}
Let $\cH_j$, $j=1,2$, be complex separable Hilbert spaces. \\
$(i)$ Assume that $A: \dom(A) \subseteq \cH_1 \to \cH_2$ is a densely 
defined, closed, linear operator satisfying  
\begin{equation}
\ker(A) = \{0\}   \lb{5.1}
\end{equation} 
and 
\begin{equation}
\inf(\sigma(A^* A)) = 0.   \lb{5.2}
\end{equation}
$(ii)$ Let $R$ be a densely defined, closable operator in $\cH_1$ satisfying
\begin{equation}
\dom (R) \cap \dom (R^*) \supseteq \dom (A).    \lb{5.3}
\end{equation}
\end{hypothesis}
%%%%%%%%%%% 

As in the previous section we start by introducing $Q_{|A|,R}$ in $\cH_1 \oplus \cH_1$ assuming Hypothesis \ref{h5.1}.
\begin{equation}
Q_{|A|,R} = \begin{pmatrix} - i \, R & |A| \\ |A| & 0 \end{pmatrix}, \quad 
\dom (Q_{|A|,R}) = \dom(|A|) \oplus \dom(|A|) \subseteq \cH_1 \oplus \cH_1,   \lb{5.4a} 
\end{equation} 
and the quadratic operator pencil $M(\cdot)$ in $\cH_1$, 
\begin{equation}
M(z) = |A|^2 - iz R - z^2 I_{\cH_1}, \quad \dom(M(z)) = \dom\big(|A|^2\big), \; z \in\bbC. 
\lb{5.4}
\end{equation}
We note that even though the pencil $M(\cdot)$ has unbounded coefficients,
replacing $M(\cdot)$ by $M(\cdot)\big(|A|^2 + I_{\cH}\big)^{-1}$ reduces matters to a 
pencil with bounded coefficients, in particular,  
\cite[Lemma\ 20.1]{Ma88} applies to the spectrum of $M(\cdot)$ in this context.

%%%%%%%%%%%%
\begin{lemma} \lb{l5.2}
Assume Hypothesis \ref{h5.1}. \\
$(i)$ Then $M(z)$, $z\in\bbC$, is a densely defined, closed operator in $\cH_1$. \\
$(ii)$ In addition,
\begin{equation}
M(z)^* = |A|^2 + i \ol{z} R^* - \ol{z}^2 I_{\cH_1}, \quad 
\dom(M(z)^*) = \dom\big(|A|^2\big), \; z \in\bbC.     \lb{5.5} 
\end{equation}
$(iii)$ Let $z\in \rho(M(\cdot))$, then  
\begin{align}
& |A| M(z)^{-1} = \ol{|A| M(z)^{-1}} \in \cB(\cH_1), \quad  
\ol{M(z)^{-1} |A|} = \big[|A| (M(z)^*)^{-1}\big]^* \in \cB(\cH_1),    \no \\
& \ol{|A| M(z)^{-1} |A|} \in \cB(\cH_1).   \lb{5.5a} 
\end{align}
\end{lemma}
%%%%%%%%%%%
\begin{proof}
The first two items can be shown as in Lemma \ref{l4.2} and so we focus 
on the proof of item $(iii)$. The fact that $|A| M(z)^{-1}$ and 
$|A| (M(z)^*)^{-1}$ are  bounded 
operators on $\cH_1$ is immediate by \eqref{5.3} and \eqref{5.5}. Similarly, 
\begin{equation} 
\ol{M(z)^{-1} |A|} = [|A| (M(z)^*)^{-1}]^* \in \cB(\cH_1). 
\end{equation} 
For the third operator in 
\eqref{5.5a} one first observes that $|A|^2 M(z)^{-1} \in \cB(\cH_1)$ by the 
closed graph theorem, implying $|A|^2 |M(z)^{-1}| \in \cB(\cH_1)$. By 
\cite[Theorem\ X.18\,(a)]{RS75} (alternatively, by Heinz's inequality, 
\cite[Theorem\ 9.4(b)]{We80}) also $|A| |M(z)|^{-1/2} \in \cB(\cH_1)$. 
Consequently, one also obtains 
$\ol{|M(z)|^{-1/2} |A|} =(|A| |M(z)|^{-1/2})^* \in\cB(\cH_1)$. Replacing $M(z)$ by 
$M(z)^*$ (cf.\ \eqref{5.5}), one also concludes that 
$|A| |M(z)^*|^{-1/2} \in \cB(\cH_1)$. Next, using the generalized polar 
decomposition for $M(z)^{-1}$ (cf.\ \cite{GMMN09}), 
\begin{equation}
M(z)^{-1} = |M(z)^*|^{-1/2} V_{|M(z)|^{-1}} |M(z)|^{-1/2}, 
\end{equation}
(with $V_{|M(z)|^{-1}}$ the partial isometry in $\cH_1$ in the standard polar decomposition of $M(z)^{-1}$, $M(z)^{-1} = V_{|M(z)|^{-1}} |M(z)|^{-1}$), one 
infers that 
\begin{equation}
\ol{|A| M(z)^{-1} |A|} =  |A| |(M(z)^{-1})^*|^{1/2} V_{|M(z)|^{-1}} 
\ol{|M(z)^{-1}|^{1/2} |A|} \in \cB(\cH_1).
 \end{equation} 
\end{proof}
%%%%%%%%%%%

Next, Theorem \ref{t4.3} requires some modifications. 

%%%%%%%%%%%%
\begin{theorem} \lb{t5.3}
Assume Hypothesis \ref{h5.1}. Then $Q_{|A|,R}$ is injective and closed, and  
\begin{align}
&(Q_{|A|,R} - z I_{\cH_1\oplus \cH_1})^{-1} =  
 \begin{pmatrix} z M(z)^{-1} & \ol{M(z)^{-1} |A|} \\
|A| M(z)^{-1} & z^{-1} \ol{|A| M(z)^{-1} |A|} - z^{-1} I_{\cH_1} \end{pmatrix},   
 \lb{5.6} \\[1mm]
& \hspace*{1.65cm} z \in \rho(Q_{|A|,R})\backslash \{0\}  
= \big\{\zeta \in\bbC \,\big|\, M(\zeta)^{-1} \in \cB(\cH_1)\big\} 
= \rho(M(\cdot)),   \no \\
& \sigma(Q_{|A|,R}) \cup \{0\} = \sigma (M(\cdot)).    \lb{5.7}
\end{align}
Suppose in addition that $R$ is accretive. Then $i\lambda\in \rho(M(\cdot))$ for 
all $\lambda> 0$ and also $i \, Q_{|A|,R}$ is accretive. Moreover,  
$- i  \, Q_{|A|,R}$ then generates a contraction semigroup in $\cH_1 \oplus \cH_1$,
denoted by $\exp (-i \, Q_{|A|,R} \, t)$, $t\geq 0$. 
\end{theorem}
%%%%%%%%%%%%
\begin{proof}
Injectivity and closedness of $Q_{|A|,R}$ follow as in the proof of Theorem \ref{t4.3}.

Next, one recalls that $0 \in \sigma (M(\cdot))$ since $M(0) = |A|^2$ and $0\in \sigma (|A|)$ 
(cf.\ \eqref{5.2}). Temporarily denoting the right-hand side of \eqref{5.6} by $T_{|A|,R}(z)$ 
for $z\in\rho(M(\cdot))$, one concludes that 
$T_{|A|,R}(z) \in  \cB(\cH_1 \oplus \cH_1)$ by Lemma \ref{l5.2}\,$(iii)$. That $T_{|A|,R}(z)=
(Q_{|A|,R} - z I_{\cH_1\oplus \cH_1})^{-1}$ can now be checked directly. This proves 
$\rho(M(\cdot)) \subseteq \rho(Q_{|A|,R})$. 
Conversely, if $z\in \rho(Q_{|A|,R})\backslash\{0\}$, the approach used in 
\eqref{4.8}--\eqref{4.8a} also works in the present more general context and hence yields 
$\rho(Q_{|A|,R})\backslash\{0\} \subseteq \rho(M(\cdot))$. This proves \eqref{5.6} and 
\eqref{5.7}.

If $R$ is accretive, then so is $i \, Q_{|A|,R}$. Next, we rewrite 
$M(i \, \lambda)$, $\lambda >0$ in \eqref{5.4} as 
\begin{align}
M(i \, \lambda) &= |A|^2 + \lambda R + \lambda^2 I_{\cH_1}    \no \\
&= \lambda \big(|A|^2 + I_{\cH_1}\big)^{1/2} 
\big[\big(|A|^2 + I_{\cH_1}\big)^{-1/2} R 
\big(|A|^2 + I_{\cH_1}\big)^{-1/2} + \lambda^{-1} I_{\cH_1}\big]     \no \\
& \quad \times \big(|A|^2 + I_{\cH_1}\big)^{1/2}, \quad \lambda > 0. 
\end{align}
Since $\big(|A|^2 + I_{\cH_1}\big)^{-1/2} R \big(|A|^2 + I_{\cH_1}\big)^{-1/2} \in \cB(\cH_1)$, 
and accretivity of $R$ implies accretivity of 
$\big(|A|^2 + I_{\cH_1}\big)^{-1/2} R \big(|A|^2 + I_{\cH_1}\big)^{-1/2}$, one concludes that  
in fact, 
\begin{equation}
\big(|A|^2 + I_{\cH_1}\big)^{-1/2} R \big(|A|^2 + I_{\cH_1}\big)^{-1/2} \, \text{ is $m$-accretive}
\end{equation} 
(cf.\ \cite[p.\ 279]{Ka80}, \cite{Ph59}) and hence 
\begin{equation}
\big[\big(|A|^2 + I_{\cH_1}\big)^{-1/2} R 
\big(|A|^2 + I_{\cH_1}\big)^{-1/2} + \lambda^{-1} I_{\cH_1}\big]^{-1} \in \cB(\cH_1), 
\quad \lambda > 0.    \lb{5.8}
\end{equation}
Equation \eqref{5.8} implies 
\begin{equation}
M(i \, \lambda)^{-1} \in \cB(\cH_1), \quad \lambda > 0.   \lb{5.9}
\end{equation}
Equation \eqref{5.6} then yields $(Q_{|A|,R} - i \, \lambda)^{-1} \in \cB(\cH_1\oplus \cH_2)$, 
$\lambda > 0$, and this fact combined with the Lumer--Phillips 
Theorem \cite{LP61} (cf.\ \cite[Theorem\ II.3.15]{EN00}) then again 
proves that $- i  \, Q_{|A|,R}$ generates a contraction semigroup. 
\end{proof}
%%%%%%%%%%%%

%%%%%%%%%%%%
\begin{remark} \lb{r5.4}
$(i)$ Without additional restrictions on $R$ it is not possible to decide whether or not 
$0 \in \sigma(Q_{|A|,R})$ (although, one always has $0 \in \sigma (M(\cdot))$). \\
$(ii)$ An alternative argument for \eqref{5.9} can be formulated as follows. An application 
of Cauchy's inequality yields
\begin{align}
& \|M(i \, \lambda) f\|_{\cH_1} \|f\|_{\cH_1}  \geq |(f, M(i \, \lambda) f)_{\cH_1}| 
\geq \Re((f, M(i \, \lambda) f)_{\cH_1})    \no \\
& \quad \geq \lambda \Re((f, R f)_{\cH_1}) + \lambda^2 \|f\|_{\cH_1}^2 
\geq \lambda^2 \|f\|_{\cH_1}^2, \quad \lambda \geq 0, \; f \in \dom\big(|A|^2\big). 
\end{align}
In particular, $\|M(i \, \lambda) f\|_{\cH_1} \geq \lambda^2 \|f\|_{\cH_1}$ yields that 
$M(i \, \lambda)$ is injective for all $\lambda > 0$. The analogous argument proves that 
also $M(i \, \lambda)^*$ is injective for all $\lambda > 0$. Thus,
\begin{equation}
\ker (M(i \, \lambda)^*) = \ran(M(i \, \lambda))^\bot = \{0\}, \quad \lambda > 0,
\end{equation}
implies that $M(i \, \lambda)$, $\lambda > 0$, is a bijection which in turn yields \eqref{5.9} since $M(i \lambda)$ and hence $M(i \lambda)^{-1}$ are closed 
in $\cH_1$. \\
$(iii)$ If in addition, $R$ is bounded with respect to $A$ with relative bound less than one, that is, 
there exist constants $0 \leq a < 1$ and $b \geq 0$ such that
\begin{equation}
\|R f\|_{\cH_1} \leq a \||A| f\|_{\cH_1} + b \|f\|_{\cH_1}, \quad f \in \dom(|A|)
\end{equation}
(here we used that $\| A f\|_{\cH_1} = \||A| f\|_{\cH_1}$, $f \in \dom(A) = \dom(|A|)$, using 
the polar decomposition for $A$ and $|A|$), one concludes that there exists $c>0$ 
such that 
\begin{align}
\big\|R (|A|^2 + \mu^2)^{-1} g\big\|_{\cH_1} & \leq 
a \big\| |A| (|A|^2 + \mu^2)^{-1} g\big\|_{\cH_1}  + b \big\|(|A|^2 + \mu^2)^{-1} g\big\|_{\cH_1} 
\no \\
& \leq a \mu^{-1} + b \mu^{-2} \leq c \mu^{-1}, \quad \mu > 0, \; g \in \cH_1. 
\end{align}
In particular, one can choose $0 < c < 1$ for $0 < \mu$ sufficiently large. This 
then yields 
\begin{align}
M(i \, \lambda)^{-1} & = (|A|^2 + \lambda R + \lambda^2 I_{\cH_1})^{-1}  \no \\
& = (|A|^2 + \lambda^2 I_{\cH_1})^{-1} 
\big[I_{\cH_1} + \lambda R (|A|^2 + \lambda^2 I_{\cH_1})^{-1} \big]^{-1} 
\in \cB(\cH_1)    \lb{5.20} \\
& \hspace*{5.15cm}  \text{for $0 < \lambda$ sufficiently large,}    \no 
\end{align}
since then $\lambda \big\|R (|A|^2 + \lambda^2)^{-1}\big\|_{\cB(\cH_1)} \leq c <1$. 
(One observes that accretivity of $R$ was not used in arriving at \eqref{5.20}.) Together with 
accretivity of $R$ this again permits 
the application of the Lumer--Phillips Theorem to the effect that $- i  \, Q_{|A|,R}$ is 
generating a contraction semigroup.
\end{remark}
%%%%%%%%%%%%

The following remark is not explicitly used in this paper, but its perturbation 
theoretic context is relevant when considering 
$Q_{|A|,R} = \left(\begin{smallmatrix} - i \, R & |A| \\ |A| & 0 
\end{smallmatrix}\right)$ as an operator sum of 
$Q_{|A|,0} = \left(\begin{smallmatrix} 0 & |A| \\ |A| & 0 \end{smallmatrix}\right)$ 
and $\left(\begin{smallmatrix} - i \, R & 0 \\ 0 & 0 \end{smallmatrix}\right)$ and 
invoking the notion of off-diagonal dominance as discussed, for instance, in 
\cite[Sect.2.2]{Tr08}. 

%%%%%%%%%%%%
\begin{remark} \lb{r5.5}
Suppose $S$ is $m$-accretive (resp., self-adjoint) in $\cH$ and $T$ is 
accretive (resp., symmetric) in $\cH$ with $\dom(T) \supseteq \dom(S)$ and 
assume that there exists constants $0 \leq a < 1$, $b \geq 0$ such that 
\begin{equation}
\|T f\|_{\cH} \leq a \|S f\|_{\cH} + b \|f\|_{\cH}, \quad f \in \dom(S).    \lb{5.21}
\end{equation}
Then $S+T$ defined on $\dom(S+T) = \dom(S)$ is $m$-accretive (resp., 
self-adjoint) in $\cH$, in particular, $S+T$ is closed on $\dom(S)$ in $\cH$ 
(cf., e.g., \cite[Sect.\ III.2]{EN00}, \cite[Sect.\ 1.6]{Go85}, 
\cite[Sects.\ IV.1, V.4]{Ka80}).

The choice $T = - S$ shows that one cannot permit $a=1$ in \eqref{5.21} as 
the zero operator is not closed on $\dom(S)$ if the latter is only dense in $\cH$. Moreover, if $a>1$, then $S+T$ 
need not be $m$-accretive on $\dom(S)$ as the following simple example shows. 
Consider  
\begin{align}
& S_0 f= f', \quad f \in \dom (S_0) = H^1_0((0,\infty)) = \big\{g \in L^2([0,\infty); dx) 
\,\big|\, g \in AC([0,R])   \no \\ 
& \hspace*{4.2cm} \text{for all $R>0$}; \, g(0) = 0; \, g' \in L^2([0,\infty); dx)\big\}.
\end{align}
Then $S_0$ is $m$-accretive and generates the semigroup 
\begin{equation}
(e^{- S_0 t} f)(x) = f_{*}(x-t), \quad f \in L^2([0,\infty); dx), \; t \geq 0,
\end{equation} 
in $L^2([0,\infty); dx)$, where $f_{*}$ denotes the extension of $f$ to $\bbR$ 
such that $f_{*}(x)$ vanishes for a.e. $x<0$. The adjoint semigroup is given by
\begin{equation}
(e^{- S_0^* t} f)(x) = f(x+t), \quad f \in L^2([0,\infty); dx), \; t \geq 0,
\end{equation} 
with generator
\begin{align}
& S_0^* f= -  f', \quad f \in \dom (S_0^*) = H^1((0,\infty)) 
= \big\{g \in L^2([0,\infty); dx) 
\,\big|\, g \in AC([0,R])   \no \\ 
& \hspace*{5.8cm} \text{for all $R>0$}; \, g' \in L^2([0,\infty); dx)\big\}.
\end{align}
In particular, $S_0^* \supseteq - S_0$. Next, consider $T_a = a S_0^*$, $a>0$.  
Then for $0<a<1$, $S_0 + T_a = (1-a) S_0$ is $m$-accretive in 
$L^2([0,\infty); dx)$, but for $a>1$, $S_0 + T_a = - (a-1) S_0$ is closed but not 
$m$-accretive in $L^2([0,\infty); dx)$.  
\end{remark}
%%%%%%%%%%%%

Still assuming Hypothesis \ref{h5.1}, we next introduce the operator $G_{A,R}$ in $\cH_A \oplus \cH_1$ by
\begin{align}
\begin{split} 
& G_{A,R} = \begin{pmatrix} 0 & J_A \\ - J_A^* & - R \end{pmatrix} 
= \begin{pmatrix} 0 & J_A \\ - A^* \wti A & - R \end{pmatrix},    \\
& \dom(G_{A,R}) = \dom \big(A^* \wti A \big) \oplus \dom(A) \subseteq \cH_A \oplus \cH_1.
\end{split}
\end{align}

The same proof as for Theorem \ref{t4.4} also yields the following result.  

%%%%%%%%%%%%%%
\begin{theorem} \lb{t5.4}
Assume Hypothesis \ref{h5.1}. Then
\begin{equation}
Q_{|A|,R} =  U_{\wti{|A|}} \, i \, G_{A,R} U_{\wti{|A|}}^{-1}.
\lb{5.22}
\end{equation}
In particular, the operator $G_{A,R}$ is densely defined and closed in the energy space 
$\cH_A \oplus \cH_1$. 

If in addition $R$ is accretive, then $i\lambda\in \rho(M(\cdot))$ for some 
$\lambda> 0$, and also $- G_{|A|,R}$ is accretive and $G_{|A|,R}$ generates a contraction semigroup 
in $\cH_A \oplus \cH_1$, denoted by $\exp (G_{|A|,R} \, t)$, $t\geq 0$.
\end{theorem}
%%%%%%%%%%%%%%

We also mention the analog of the result \eqref{5.22} when replacing 
$|A|$ by $A$.

%%%%%%%%%%%%%%
\begin{theorem} \lb{t5.5}
Assume Hypothesis \ref{h5.1}. Then
\begin{equation}
Q_{A,R} \big(I_{\cH_1} \oplus [I_{\cH_2} - P_{\ker(A^*)}]\big) 
=  U_{\wti A} \, i \, G_{A,R} U_{\wti A}^{-1}.    \lb{5.23}
\end{equation}
\end{theorem}
%%%%%%%%%%%%%%

We continue with an illustrative example.

%%%%%%%%%%%%%%
\begin{example} \lb{e5.7}
Considering the prototypical example of a nonnegative operator $A\geq 0$ 
in a Hilbert space $\cH$, with $\ker(A)=\{0\}$, one can, without loss of 
generality, restrict one's attention 
to the case of $\cH= L^2([0,\infty);d\rho)$ $($with $\rho$ a Borel measure 
on $[0,\infty)$ satisfying $\rho(\{0\}) = 0$$)$ 
and $A$ being the operator of multiplication with the independent variable 
so that 
\begin{equation}
\sigma (A) = \supp(d\rho)  
\end{equation}
$($with $\supp(\cdot)$ denoting the topological, i.e., smallest closed, 
support\,$)$. 
For example, one could simply choose Lebesgue measure 
$d\rho(\la)= d\la$ on $[0,\infty)$. Introducing the weighted $L^2$-spaces
\begin{equation}
\cH^{(n)} = L^2\big([0,\infty); \la^{2n}d\rho(\la)\big), \quad n \in\bbN_0, 
\end{equation}
one verifies that 
\begin{align}
\begin{split}
& \, \cH=\cH^{(0)}, \quad \cH_A = \cH^{(1)}, \quad 
\dom\big(A^n\big) =\cH^{(0)} \cap \cH^{(n)},  \; n \in\bbN,    \\ 
& \dom\big(A \wti{A}\big) = \cH^{(1)} \cap \cH^{(2)}.
\end{split}
\end{align}
Moreover, one notes that there will be elements in $\cH_A \backslash \cH$ if and only if $\inf (\si(A)) = 0$, since otherwise
the natural imbedding $\iota_A$ would be continuous by the closed graph theorem. 
Similarly there will be elements in $\cH \backslash \cH_A$ if and only if $A$ is unbounded $($i.e., if and only if $\sup (\si(A))=\infty$$)$.

Next, one checks that the unique classical solution of the  ${\rm ACP}$
\begin{align}
\begin{split}
&\dot Y (t) = G_{A,0} Y (t), \; t\geq 0,  \\
& Y (0) = \begin{pmatrix} y_0 \\  z_0 \end{pmatrix} \in\dom(G_{A,0}),
\end{split}
\end{align} 
is given by
\begin{equation}
Y(t) = e^{G_{A,0} t} Y(0), \quad t\geq 0, 
\end{equation}
with
\begin{equation}
e^{G_{A,0} t} =  
\begin{pmatrix} \cos(t \la) & \la^{-1} \sin(t \la) \\ 
-\la \sin(t \la) & \cos(t \la) \end{pmatrix}, \quad t\geq 0.
\end{equation}
Similarly, the classical solution of ${\rm ACP}_2$
\begin{align}
\begin{split}
&\ddot{u} (t) + A^2 u (t) =0, \; t\geq 0,  \\
& u(0) = u_0 \in \dom\big(A^2\big), \; \dot u (0) = u_1 \in \dom(A),
\end{split}
\end{align} 
is given by
\begin{equation}
u(t,\la) = \cos(t \la) u_0(\la) + \la^{-1} \sin(t \la) u_1(\la), \quad t\geq 0, \; 
\la\geq 0.
\end{equation}
Hence, if one chooses $Y(0) \in \cH^{(0)} \oplus \cH^{(1)}$ with support in $(0,1)$ 
$($such that automatically, $Y(0) \in \cH^{(1)} \oplus\cH^{(2)}$$)$, 
then the first component of $Y(t)$ will solve ${\rm ACP}_2$ if and only if
$y_0\in\cH^{(0)}$. In particular, there are classical solutions of  ${\rm ACP}$ 
which to not correspond to classical solutions of ${\rm ACP}_2$ if 
$\inf (\si(A))=0$.
\end{example}
%%%%%%%%%

Concerning conserved quantities in connection with the abstract wave equations we now mention the following result. 

%%%%%%%%%%%%%%%%
\begin{lemma} \lb{l5.7} 
Let $\cH_j$, $j=1,2$, be complex separable Hilbert spaces and assume that 
$A: \dom(A) \subseteq \cH_1 \to \cH_2$ is a densely defined, closed, linear operator. 
Let $B: \dom(B) \subseteq \cH_1 \to \cH_1$ be some closed operator which commutes with $|A|$ in the sense that $B |A| \subseteq |A| B$. In addition, let $R$ be a densely defined, closable operator in $\cH_1$ satisfying $\dom (R) \supseteq \dom (A)$.
Suppose $u$ is a classical solution of
\begin{equation}
\ddot u + R \dot u + A^* A u = 0,  
\end{equation}
such that $B |A| u(\cdot) \in C^1([0,\infty);\cH_1)$, $B \dot u(\cdot) \in C^1([0,\infty);\cH_1)$
and $R \dot u(\cdot) \in \dom(B)$. Then
\begin{equation}
\f{d}{dt} \big[\| B \dot u\|^2_{\cH_1} + \| A B u\|^2_{\cH_1}\big]    
= - 2\Re (B \dot u, B R \dot u)_{\cH_1}.     \lb{5.26}
\end{equation}
In particular, the right-hand side of \eqref{5.26} vanishes if $R=0$, that is, 
in the absence of damping, and hence 
$\big[\| B \dot u\|^2_{\cH_1} + \| A B u\|^2_{\cH_1}\big]$ represents a 
family of conserved quantities for $\ddot u + A^* A u = 0$. 
\end{lemma}
%%%%%%%%%%%%%%%%
\begin{proof}
One computes
\begin{align}
& \f{d}{dt} \big[\| B \dot u\|^2_{\cH_1} + \| A B u\|^2_{\cH_1}\big] 
= \f{d}{dt} \big[\| B \dot u\|^2_{\cH_1} + \| |A| B u\|^2_{\cH_1}\big]   \no \\
& \quad = 2 \Re \big( (B \dot u, B \ddot u)_{\cH_1} 
+ (|A| B u, |A| B \dot u)_{\cH_1}  \big) \no \\
& \quad = 2 \Re \big( (B \dot u, B [- R \dot u - |A|^2 u])_{\cH_1}
+ (B u, |A|^2 B \dot u)_{\cH_1} \big)  \no \\
& \quad = 2\Re \big(- (B \dot u, B R \dot u)_{\cH_1}
- (B \dot u, [B |A|^2 -|A|^2 B] u)_{\cH_1}\big)  \no \\
& \quad = - 2\Re\big( (B \dot u, B R \dot u)_{\cH_1}\big).     \lb{5.27}
\end{align}
\end{proof}
%%%%%%%%%%%%%%%%

One observes that the special case $B = I_{\cH_1}$ in \eqref{5.26} is usually associated
with the energy of the abstract wave equation 
$\ddot u + A^* A u = 0$ (resp., $\ddot u + R \dot u + A^* A u = 0$). Typical examples 
for $B$ would be $B = |A|^\alpha$, $\alpha \in [0,1]$. 

%%%%%%%%%%%%%%%%
\begin{remark} \lb{r5.8}
More generally, let $\alpha \in \bbR$ if $B^* B \geq \varepsilon I_{\cH_1}$ 
for some $\varepsilon > 0$ 
and $\alpha \geq 0$ if $\inf(\sigma(B^* B)) = 0$. Assuming that $u$ satisfies 
\begin{equation}
\ddot u + R \dot u + \big[(B^* B)^2 + C^* C (B^* B)^{\alpha}\big] u = 0,  
\end{equation}
and assuming additional appropriate conditions on $u, \dot u, R, B, C$, one obtains 
\begin{align}
& \f{d}{dt} \big[\| |B|^{\alpha} \dot u\|^2_{\cH_1} + \| |B|^{\alpha + 2} u\|^2_{\cH_1} 
+ \| |C| |B|^{2\alpha}u\|^2_{\cH_1}\big]    \no \\
& \quad = - 2 \Re (|B|^{\alpha} R \dot u, |B|^{\alpha} \dot u)_{\cH_1}.     \lb{5.29}
\end{align}
Again, the right-hand side of \eqref{5.29} vanishes if $R=0$, that is, 
in the absence of damping. 

A situation equivalent to the special case $\alpha = 1$ has recently 
been studied in \cite{LW09} in the concrete context of plate equations. 
In this connection we recall that 
$\| |C| |B|^{2\alpha}u\|^2_{\cH_1} = \| C |B|^{2\alpha}u\|^2_{\cH_1}$, etc. \\[1mm] 
%%%%%%%%%%%%%%%%
{\it Sketch of proof of \eqref{5.29}.}
\begin{align}
& \f{d}{dt} \big[\| |B|^{\alpha} \dot u\|^2_{\cH_1} + \| |B|^{\alpha + 2} u\|^2_{\cH_1} 
+ \| |C| |B|^{2\alpha}u\|^2_{\cH_1}\big]    \no \\
& \quad = 2\Re \Big(  (|B|^{\alpha} \ddot u, |B|^{\alpha} \dot u)_{\cH_1}
 + (|B|^{\alpha + 2} \dot u, |B|^{\alpha + 2} u)_{\cH_1}  \no \\
& \qquad + (|C| |B|^{2\alpha} \dot u, |C| |B|^{2\alpha} u)_{\cH_1} \Big)  \no \\
& \quad = 2\Re \Big( (|B|^{\alpha }[- R \dot u - |B|^4 u - |C|^2 |B|^{2\alpha} u], |B|^{\alpha} \dot u)_{\cH_1}  \no \\
& \qquad + (|B|^{\alpha + 2} \dot u, |B|^{\alpha + 2} u )_{\cH_1}
+ (|C| |B|^{2\alpha} \dot u, |C| |B|^{2\alpha} u)_{\cH_1} \Big)  \no \\
& \quad = 2\Re\Big( - (|B|^{\alpha} \dot u, |B|^{\alpha} R \dot u)_{\cH_1}
- (|B|^{\alpha} |C|^2 |B|^{2\alpha} u, |B|^{\alpha} \dot u)_{\cH_1} \no \\
& \qquad + (|C| |B|^{2\alpha} \dot u, |C| |B|^{2\alpha} u)_{\cH_1} \Big)  \no \\
& \quad = - 2 \Re (|B|^{\alpha} R \dot u, |B|^{\alpha} \dot u)_{\cH_1}.     \lb{5.30}
\end{align}
\end{remark}
%%%%%%%%%%%%%%%%

%%%%%%%%%%%%%%%%%%%%%%%%%%%%%%%%%%%%%%%
%%%%%%%%%%%%%%%%%%%%%%%%%%%%%%%%%%%%%%%
\section{Equipartition of Energy for Supersymmetric Dirac-Type Operators 
and Abstract Wave Equations}   \lb{s6}
%%%%%%%%%%%%%%%%%%%%%%%%%%%%%%%%%%%%%%%
%%%%%%%%%%%%%%%%%%%%%%%%%%%%%%%%%%%%%%%

In this section we briefly revisit the notion of asymptotic equipartition for abstract wave equations (in the absence of damping) and show that it implies the same phenomenon for a class of supersymmetric Dirac-type operators. 

We start with our basic hypothesis.

%%%%%%%%%%%
\begin{hypothesis} \lb{h6.1}
Let $\cH_j$, $j=1,2$, be complex separable Hilbert spaces and  
assume that $A: \dom(A) \subseteq \cH_1 \to \cH_2$ is a densely 
defined, closed, linear operator. 
\end{hypothesis}
%%%%%%%%%%% 

Assuming Hypothesis \ref{h6.1}, we introduce the supersymmetric Dirac operator 
(also known as ``supercharge'') by 
\begin{equation}
Q = \begin{pmatrix} 0 & A^*\\ A & 0 \end{pmatrix}, \quad 
\dom(Q)= \dom(A) \oplus \dom(A^*) \subseteq \cH_1 \oplus \cH_2. 
\end{equation}
(For simplicity we now use the simplifying notation $Q$ rather than the symbol 
$Q_{A,0}$ in previous sections.) 
As discussed in Appendix \ref{sA}, $Q$ is self-adjoint in $\cH_1 \oplus \cH_1$.
 
A number of Dirac operators, including the free one (i.e., one without electromagnetic potentials), one with a Lorentz scalar potential, one describing electrons in a magnetic field, one describing neutrons in an electric field, and the one modeling particles with anomalous electric moment in a magnetic field can all be put in this form (cf.\ 
\cite[Section\ 5.5]{Th92} for details). 

The solution of the corresponding time dependent Dirac equation
\begin{equation}
i \frac{d}{dt} \Psi(t) = Q \Psi(t), \quad \Psi(t) = (\psi_1(t), \psi_2(t))^\top \in \dom(Q), 
\; t\in\bbR, 
\end{equation}
is given by
\begin{equation}
\Psi(t) = e^{-i Q t} \Psi(0), \quad t \in \bbR,
\end{equation}
with $e^{-i Q t}$, $t\in\bbR$, a unitary group in $\cH_1 \oplus \cH_2$. 

One of the principal aims in this section is to prove the following result.

%%%%%%%%%%
\begin{theorem} \lb{t6.2}
Assume Hypothesis \ref{h6.1}.
Suppose $\Psi(t) = e^{-i Q t} \Psi(0)$ with $\Psi(t)=(\psi_1(t),\psi_2(t))^\top$, $t\in\bbR$,  and $\Psi(0) \in \cH_1 \oplus \cH_2$ arbitrary. Then the following assertions 
$(i)$--$(iv)$ are equivalent. \\
$(i)$ $\lim_{t\to\pm\infty} \|\psi_j(t)\|^2_{\cH_j} 
= \|\Psi(0)\|^2_{\cH_1 \oplus \cH_2}/2$, $j=1,2$. \\ 
$(ii)$ $\wlim_{t\to\infty}e^{-i Q t} =0$.  \\
$(iii)$ $\wlim_{t\to\infty}e^{-i |A| t} = 0$. \\
$(iv)$ $\wlim_{t\to\infty} \cos(|A| t) = 0$.  \\
Similarly,
\begin{align}
\begin{split} \label{6.5}
& \lim_{t\to \pm \infty} \frac{1}{t} \int_0^t ds \, \|\psi_j(s)\|^2_{\cH_j} 
= \frac{1}{2} \|\Psi(0)\|^2_{\cH_1 \oplus \cH_2}, \; j=1,2,    \\
& \quad \text{if and only if $0$ is not an eigenvalue of $Q$.}
\end{split}
\end{align}
\end{theorem}
%%%%%%%%%%%
\begin{proof}
First of all one notes that neither $(i)$--$(iv)$ nor \eqref{6.5} can hold if $0$ is an 
eigenvalue of $Q$. Hence we assume without loss of generality that 
\begin{equation}
\ker(Q)=\ker(A) \oplus \ker(A^*) = \{ 0\}.  
\end{equation}
Moreover, since $\big(e^{-i Q t}\big)^* = e^{i Q t}$, it 
suffices to study the
limit $t\to \infty$ in Theorem \ref{t6.2}\,$(i)$ and \eqref{6.5}. Next, we recall \eqref{A.2}--\eqref{A.6}, 
\eqref{A.22a}, \eqref{A.22b}, the polar decomposition $A= V_A |A|$, 
$A^* = (V_A)^* |A^*|$ where, due to our assumption 
$\ker(A) =\ker(A^*) = \{0\}$ and hence $V_A \in \cB(\cH_1, \cH_2)$ is unitary. In addition, 
we use the notation $H_1=A^* A$, $H_2=A A^*$ (cf.\ Appendix \ref{sA} for details). 
Then, by the spectral theorem applied to $Q$,
\begin{equation}
e^{-i Q t} = \cos(|Q| t) - i \sin(|Q| t) V_Q, \quad t \in \bbR, 
\end{equation}
and by \eqref{A.22b}, one obtains
\begin{equation}
e^{-i Q t} =  \begin{pmatrix} \cos(|A| t) & -i \sin(|A| t) (V_A)^* \\[1mm] 
-i \sin(|A^*| t) V_A & \cos(|A^*| t) \end{pmatrix}, \quad t \in\bbR.
\end{equation}
Taking scalar products of $e^{-i Q t}$ with vectors of the type $(f, 0)^\top$ 
and $(0, g)^\top$ then shows that 
\begin{align} 
\begin{split}
&\wlim_{t\to\infty} e^{-i Q t} = 0  \\ 
& \quad \text{if and only if} \,  
\begin{cases} \wlim_{t\to\infty} \cos(|A| t) = 0, 
\; \wlim_{t\to\infty} \cos(|A^*| t) = 0, \\
\wlim_{t\to\infty} \sin(|A| t) = 0, 
\; \wlim_{t\to\infty} \sin(|A^*| t) = 0.  \end{cases} 
\end{split} 
\end{align}
However, since $\big(e^{- i |T| t}\big)^* = e^{i |T| t}$ for any densely defined closed 
operator $T$, one actually infers that
\begin{align} 
\begin{split}
& \wlim_{t\to\infty} e^{-i Q t} = 0   \, \text{ if and only if } 
\big\{ \wlim_{t\to\infty} \cos(|A| t) = 0, 
\; \wlim_{t\to\infty} \cos(|A^*| t) = 0\big\}  \\
& \quad \text{ if and only if }  \big\{ \wlim_{t\to\infty} e^{ -i |A| t} = 0, 
\; \wlim_{t\to\infty} e^{ -i |A^*| t} = 0\big\}.     
\end{split} 
\end{align} 
Finally, since $H_1= A^* A$ and $H_2= A^* A$ are unitarily equivalent (recalling 
\eqref{A.10a} and the fact that $P_{\ker(A)} = P_{\ker(A^*)} = 0$), this actually yields 
that 
\begin{equation} 
\wlim_{t\to\infty} e^{-i Q t} = 0 \, \text{ if and only if }  
\wlim_{t\to\infty} e^{-i |A| t} = 0 \, \text{ if and only if }   
\wlim_{t\to\infty} \cos(|A| t) = 0.   \lb{6.12} 
\end{equation}

Given $\Psi(0)=(\psi_1,\psi_2)^\top \in \cH_1 \oplus \cH_2$, one then computes 
\begin{align} \no
& \| \psi_1(t) \|^2_{\cH_1} = \| \cos(|A| t) \psi_1 
- i \sin(|A| t) (V_A)^* \psi_2 \|^2_{\cH_1}  \\ \no
& \quad = \frac{1}{4} \big\| \big[e^{i |A| t} + e^{-i |A| t}\big] \psi_1
+ \big[e^{i |A| t} - e^{-i |A| t}\big] (V_A)^* \psi_2 \big\|^2_{\cH_1}   \\ \no
& \quad = \frac{1}{4} \big\| \big[e^{2i |A| t} + I_{\cH_1} \big] \psi_1 
+ \big[e^{2i |A| t} - I_{\cH_1}\big] (V_A)^* \psi_2 \big\|^2_{\cH_1}   \\ \no
& \quad = \frac{1}{4} \big\| [\psi_1 - (V_A)^* \psi_2] 
+ e^{2i |A| t} [\psi_1 + (V_A)^* \psi_2] \big\|^2_{\cH_1}    \\ \no
& \quad = \frac{1}{4} \| [\psi_1 - (V_A)^* \psi_2]\|_{\cH_1}  
+ \frac{1}{4} \big\| e^{2i |A| t} [\psi_1 + (V_A)^* \psi_2]\big\|^2_{\cH_1}    \\ \no
& \qquad {} + \frac{1}{2} \Re \big(\big([\psi_1 - (V_A)^* \psi_2], 
e^{2i |A| t} [\psi_1 + (V_A)^* \psi_2]\big)_{\cH_1}\big)      \\ \no
& \quad = \frac{1}{2} \big(\|\psi_1\|^2_{\cH_1} + \| (V_A)^* \psi_2\|^2_{\cH_1}\big)  \\ \no
& \qquad + \frac{1}{2} \Re \big(\big([\psi_1 - (V_A)^* \psi_2], 
 e^{2i |A| t} [\psi_1 + (V_A)^* \psi_2]\big)_{\cH_1}\big)   \\ \no
& \quad = \frac{1}{2} \big(\|\psi_1\|^2_{\cH_1} + \| \psi_2\|^2_{\cH_2}\big)
 + \frac{1}{2} \Re \big(\big([\psi_1 - (V_A)^* \psi_2], 
 e^{2i |A| t} [\psi_1 + (V_A)^* \psi_2]\big)_{\cH_1}\big)    \\  
& \quad = \frac{1}{2} \|\Psi(0)\|^2_{\cH_1 \oplus \cH_2}
 + \frac{1}{2} \Re \big(\big([\psi_1 - (V_A)^* \psi_2], 
 e^{2i |A| t} [\psi_1 + (V_A)^* \psi_2]\big)_{\cH_1}\big).    \lb{6.13}
\end{align}
Thus, $\wlim_{t\to\infty} e^{-i |A| t} = 0$ yields 
$\lim_{t\to\infty} \| \psi_1(t) \|^2_{\cH_1} = \frac{1}{2} \|\Psi(0)\|^2_{\cH_1 \oplus \cH_2}$, 
and hence also 
$\lim_{t\to\infty} \| \psi_2(t) \|^2_{\cH_2} = \frac{1}{2} \|\Psi(0)\|^2_{\cH_1 \oplus \cH_2}$, 
since $e^{- i Q t}$, $t\in\bbR$, is unitary on $\cH_1 \oplus \cH_2$.  

Conversely, choose $\varphi, \psi \in \cH_1$ and set  
$\psi_1 = (\psi + \varphi)/2$ and 
$\psi_2 = V_A (\varphi - \psi)/2$. Then \eqref{6.13} shows that 
$\lim_{t\to\infty} \| \psi_1 (t) \|^2_{\cH_1} = \frac{1}{2} \|\Psi(0)\|^2_{\cH_1 \oplus \cH_2}$
implies 
\begin{equation}
\lim_{t\to\infty}  \Re \big(\big(\psi, e^{2i |A| t} \varphi\big)_{\cH_1}\big) = 0, 
\quad \psi,\varphi\in\cH_1.
\end{equation}
In particular,
\begin{equation}
\lim_{t\to\infty}  \Re \big(\big(\psi, e^{2i |A| t} \psi\big)_{\cH_1}\big) 
= \lim_{t\to\infty} (\psi, \cos(2 |A| t) \psi)_{\cH_1} = 0, 
\quad \psi\in\cH_1.
\end{equation} 
By polarization for sesquilinear forms, this is equivalent to
\begin{equation}
\lim_{t\to\infty} (\psi, \cos(2 |A| t) \varphi)_{\cH_1} = 0, \quad \psi,\varphi\in\cH_1,
\end{equation}
and thus by \eqref{6.12} also to
\begin{equation}
\lim_{t\to\infty} \big(\Psi, e^{i Q t} \Psi\big)_{\cH_1 \oplus \cH_2} = 0, \quad 
\Psi\in \cH_1 \oplus \cH_2,
\end{equation}
proving the equivalence of $(i)$--$(iv)$.

Applying von Neumann's mean ergodic theorem in the weak sense (cf., e.g., 
\cite[Corollary\ 5.2]{Da80}, \cite[Theorem\ 1.8.20]{Go85}) to \eqref{6.13} yields 
the Ces{\` a}ro limit
\begin{align}
& \lim_{t\to \pm \infty} \frac{1}{t} \int_0^t ds \, \|\psi_1 (s)\|^2_{\cH_1} 
= \frac{1}{2} \|\Psi(0)\|^2_{\cH_1 \oplus \cH_2}   \no \\
& \qquad + \frac{1}{2} \Re \bigg(\lim_{t\to \pm \infty} \frac{1}{t} \int_0^t ds \, 
\big([\psi_1 - (V_A)^* \psi_2], 
 e^{2i |A| s} [\psi_1 + (V_A)^* \psi_2]\big)_{\cH_1}\bigg)   \no \\
& \quad = \frac{1}{2} \|\Psi(0)\|^2_{\cH_1 \oplus \cH_2}  
+ \frac{1}{2} \Re \big( \big([\psi_1 - (V_A)^* \psi_2], E_{|A|} (\{0\}) 
 [\psi_1 + (V_A)^* \psi_2]\big)_{\cH_1}\big)    \no \\
& \quad = \frac{1}{2} \|\Psi(0)\|^2_{\cH_1 \oplus \cH_2}  
\end{align}
if $E_{|A|} (\{0\}) = 0$. (Here we used the notation $E_S(\cdot)$ for the 
strongly right continuous spectral 
family associated with the self-adjoint operator $S$.) 
Conversely, choose again $\varphi, \psi \in \cH_1$ and set  
$\psi_1 = (\psi + \varphi)/2$ and 
$\psi_2 = V_A (\varphi - \psi)/2$. Then \eqref{6.13} shows that 
\begin{equation}
\lim_{t\to\infty} \f{1}{t} \int_0^t ds \, \| \psi_1 (s) \|^2_{\cH_1} 
= \frac{1}{2} \|\Psi(0)\|^2_{\cH_1 \oplus \cH_2}
\end{equation} 
implies 
\begin{equation}
 \Re \bigg(\lim_{t\to\infty} \f{1}{t} \int_0^t ds \, 
 \big(\psi, e^{2i |A| s} \varphi\big)_{\cH_1}\bigg) = 0, 
\quad \psi,\varphi\in\cH_1.
\end{equation}
In particular,
\begin{align}
& \Re \bigg(\lim_{t\to\infty} 
\f{1}{t}  \int_0^s ds \, \big(\psi, e^{2i |A| s} \psi\big)_{\cH_1}\bigg) 
= \lim_{t\to\infty} \f{1}{t} \int_0^t ds \, (\psi, \cos(2 |A| s) \psi)_{\cH_1} = 0,  \no\\ 
& \hspace*{9.5cm} \psi\in\cH_1.  
\end{align} 
By polarization for sesquilinear forms, this is equivalent to
\begin{equation}
\lim_{t\to\infty} \f{1}{t} \int_0^t ds \, (\psi, \cos(2 |A| s) \varphi)_{\cH_1} 
= 0, \quad \psi,\varphi\in\cH_1.     \lb{6.21}
\end{equation}
Since generally, as a corollary of von Neumann's weak ergodic theorem 
(cf.\ \cite{GRS78})
\begin{equation}
\lim_{t\to\infty} \f{1}{t} \int_0^t ds \, (\psi, \cos(|A| s) \varphi)_{\cH_1} 
= (\psi, E_{|A|} (\{0\}) \varphi)_{\cH_1}, \quad \psi,\varphi\in\cH_1, 
\end{equation}
\eqref{6.21} yields $E_{|A|} (\{0\}) = 0$. The same computation with $\psi_1(\cdot)$ 
replaced by $\psi_2(\cdot)$ then proves that 
\begin{equation}
\lim_{t\to\infty} \f{1}{t} \int_0^t ds \, \| \psi_2 (s) \|^2_{\cH_2} 
= \frac{1}{2} \|\Psi(0)\|^2_{\cH_1 \oplus \cH_2} \, \text{ if and only if } \, 
E_{|A^*|} (\{0\}) = 0,
\end{equation} 
proving \eqref{6.5} (cf.\ \eqref{A.17}). 
\end{proof}
%%%%%%%%%

We note that 
\begin{equation}
e^{-i |Q| t} = \cos(|Q| t) - i \sin(|Q| t)
=  \begin{pmatrix} e^{- i |A| t} & 0 \\  
0 & e^{- i |A^*| t} \end{pmatrix}, \quad t \in \bbR,   
\end{equation}
and hence \eqref{6.12} also yields that  
\begin{equation} 
\text{if $\ker(Q)=\{0\}$, then } \, \wlim_{t\to\infty} e^{-i Q t} = 0  
\, \text{ if and only if } \,  
\wlim_{t\to\infty} e^{-i |Q| t} = 0. 
\end{equation}

%%%%%%%%%
\begin{remark} \lb{r6.3}
The proof of Theorem \ref{t6.2} is similar in spirit to the proof for equipartition of energy for abstract wave equations \cite{Br67} (see also \cite{Go69}, \cite{Go70}, 
\cite[Theorems\  7.12 and 7.13]{Go85}, \cite{Go86}, \cite{GLS93}, \cite{GR80}, 
\cite{GS76}, \cite{GS79}, \cite{GS82a}, \cite{GS82}, \cite{GS87}, \cite{Sh68},   
and the references therein). In fact, since the two problems are unitarly equivalent, 
one follows from the other. For the benefit of the reader we decided to provide the 
proof in the context of supersymmetric Dirac-type operators. 
\end{remark}
%%%%%%%%%%

For completeness we finally recall the corresponding result concerning the asymptotic equipartition for abstract wave equations in the absence of damping, which motivated the derivation of Theorem \ref{t6.2}.

Consider the initial value problem 
\begin{align}
\begin{split}
& \ddot u (t) + A^*A u (t) = 0, \quad  t \in \bbR,  \\
& u(0) = f_0 \in \dom (A^* A), \;  \dot u (0) = f_1 \in \dom(A).     \lb{6.26} 
\end{split}  
\end{align}
Introducing kinetic and potential energies, $K_u(t)$ and $P_u(t)$, associated 
with a (strong) solution $u(\cdot)$ of \eqref{6.26} at time $t \in \bbR$, 
\begin{equation}
K_u(t) = \|\dot u (t)\|^2_{\cH_1}, \quad P_u(t)  = \|A u\|^2_{\cH_1} 
= \||A| u\|^2_{\cH_1}, \quad t\in\bbR,
\end{equation} 
one recalls 
conservation of the total energy (cf.\ Lemma \ref{l5.7})
\begin{equation}
K_u(t) + P_u(t) = K_u (0) + P_u(0), \quad t \in\bbR. 
\end{equation}
Moreover, the initial value problem \eqref{6.26} is said to admit {\it asymptotic equipartition of energy} if
\begin{equation}
\lim_{t\to\pm\infty} K_u(t) = \lim_{t\to\pm\infty} P_u(t) 
= \f{1}{2} [K_u (0) + P_u(0)]. 
\end{equation}
Asymptotic equipartition of energy has extensively been discussed in the literature, 
we refer, for instance, to \cite{BG84}, \cite{Go69}, \cite{Go70}, 
\cite{Go85a}, \cite{Go86}, \cite{Go93}, \cite{GLS93}, \cite{GRS78}, \cite{GR80}, 
\cite{GS76}, \cite{GS79}, \cite{GS81}, 
\cite{GS82a}, \cite{GS82}, \cite{GS87}, and \cite{Sh68}. In particular, the following 
theorem appeared in Goldstein \cite[Theorems\ 7.12 and 7.13]{Go85}.  

%%%%%%%%%%
\begin{theorem} \lb{t6.4} 
Assume Hypothesis \ref{h6.1} and let $u(\cdot): \bbR \to \cH_1$ be a solution 
of \eqref{6.26}.\ Then the following assertions $(i)$ and $(ii)$ are equivalent.  \\
$(i)$ $\lim_{t\to\pm\infty} K_u(t) = \lim_{t\to\pm\infty} P_u(t) 
= [K_u (0) + P_u(0)]/2$.  \\
$(ii)$ $\wlim_{t\to\infty}e^{-i |A| t} = 0$. \\
Similarly,
\begin{align}
\begin{split} \label{6.30}
& \lim_{t\to \pm \infty} \frac{1}{t} \int_0^t ds \,K_u(s) 
=  \lim_{t\to \pm \infty} \frac{1}{t} \int_0^t ds \, P_u(s)  
= \f{1}{2}[K_u (0) + P_u(0)]    \\
& \quad \text{if and only if $0$ is not an eigenvalue of $A$.}
\end{split}
\end{align}
\end{theorem}
%%%%%%%%%%

%%%%%%%%%%%%%%%%%%%%%%%%%%%%%%%%%%%%%%
%%%%%%%%%%%%%%% appendices %%%%%%%%%%%%%%%%
\appendix
%%%%%%%%%%%% Appendix A %%%%%%%%%%%%%%
\section{Supersymmetric Dirac-Type Operators in a Nutshell} \lb{sA}
\renewcommand{\theequation}{A.\arabic{equation}}
\renewcommand{\thetheorem}{A.\arabic{theorem}}
\setcounter{theorem}{0} \setcounter{equation}{0}
%%%%%%%%%%%%%%%%%%%%%%%%%%%%%%%%%%%%%%
%%%%%%%%%%%%%%%%%%%%%%%%%%%%%%%%%%%%%%

In this appendix we briefly summarize some results on supersymmetric 
Dirac-type operators and commutation methods due to \cite{De78}, 
\cite{GSS91}, \cite{Th88}, and \cite[Ch.\ 5]{Th92} (see also \cite{HKM00}). 

The standing assumption in this appendix will be the following.

%%%%%%%%%%%%%
\begin{hypothesis} \lb{hA.1}
Let $\cH_j$, $j=1,2$, be separable complex Hilbert spaces and 
\begin{equation}
T: \cH_1 \supseteq \dom(T) \to \cH_2    \lb{A.1}
\end{equation} 
be a densely defined closed linear operator. 
\end{hypothesis}
%%%%%%%%%%%%%

We define the self-adjoint Dirac-type operator in $\cH_1 \oplus \cH_2$ by 
\begin{equation}
Q = \begin{pmatrix} 0 & T^* \\ T & 0 \end{pmatrix}, \quad \dom(Q) = \dom(T) \oplus \dom(T^*).     
\lb{A.2}
\end{equation}
Operators of the type $Q$ play a role in supersymmetric quantum mechanics (see, e.g., the extensive list of references in \cite{BGGSS87}). Then,
\begin{equation}
Q^2 = \begin{pmatrix} T^* T & 0 \\ 0 & T T^* \end{pmatrix}     \lb{A.3}
\end{equation}
and for notational purposes we also introduce
\begin{equation}
H_1 = T^* T \, \text{ in } \, \cH_1, \quad H_2 = T T^* \, \text{ in } \, \cH_2.    \lb{A.4}
\end{equation}
In the following, we also need the polar decomposition of $T$ and $T^*$, that is, the representations
\begin{align}
T& = V_T |T| = |T^*| V_T = V_T T^* V_T  \, \text{ on } \, \dom(T) = \dom(|T|),  
\label{A.4a} \\
T^*& = V_{T^*} |T^*| = |T|V_{T^*} = V_{T^*} T V_{T^*}   \, \text{ on } \, 
\dom(T^*) = \dom(|T^*|),   \label{A.4b} \\
|T|& = V_{T^*} T=T^* V_T = V_{T^*} |T^*| V_T \, \text{ on } \, \dom(|T|),    \label{A.4c} \\
 |T^*|& = V_T T^* = T V_{T^*} = V_T |T| V_{T^*}  \, \text{ on } \, \dom(|T^*|),   
 \label{A.4d}
\end{align}
where
\begin{align}
& |T| = (T^* T)^{1/2}, \quad |T^*| = (T T^*)^{1/2},     \lb{A.6} \\
& V_{T^*} = (V_T)^*,     \lb{A.6a}
\end{align}
and 
\begin{equation}
V_{T^*} V_T = P_{\ol{{\ran}(|T|)}}=P_{\ol{{\ran}(T^*)}} \, ,  \quad
V_T V_{T^*} = P_{\ol{{\ran}(|T^*|)}}=P_{\ol{{\ran}(T)}} \, .      \label{A.4e}
\end{equation}
In particular, $V_T$ is a partial isometry with initial set $\ol{{\ran}(|T|)}$
and final set $\ol{{\ran}(T)}$ and hence $V_{T^*}$ is a partial isometry with initial 
set $\ol{\ran(|T^*|)}$ and final set $\ol{\ran(T^*)}$. In addition, 
\begin{equation}
V_T = \begin{cases} \ol{T (T^* T)^{-1/2}} = \ol{(T T^*)^{-1/2}T} & \text{on }  (\ker (T))^{\bot},  \\
0 & \text{on }  \ker (T).  \end{cases}     \lb{A.7}
\end{equation} 

Next, we collect some properties relating $H_1$ and $H_2$.

%%%%%%%%%%%
\begin{theorem} [\cite{De78}] \lb{tA.2}  
Assume Hypothesis \ref{hA.1} and let $\phi$ be a bounded Borel measurable 
function on $\bbR$. \\ 
$(i)$ One has
\begin{align}
& \ker(T) = \ker(H_1) = (\ran(T^*))^{\bot}, \quad \ker(T^*) = \ker(H_2) = (\ran(T))^{\bot},  \lb{A.8} \\
& V_T H_1^{n/2} = H_2^{n/2} V_T, \; n\in\bbN, \quad 
V_T \phi(H_1) = \phi(H_2) V_T.  \lb{A.9} 
\end{align}
$(ii)$ $H_1$ and $H_2$ are essentially isospectral, that is, 
\begin{equation}
\sigma(H_1)\backslash\{0\} = \sigma(H_2)\backslash\{0\},    \lb{A.10}
\end{equation}
in fact, 
\begin{equation}
T^* T [I_{\cH_1} - P_{\ker(T)}] \, \text{ is unitarily equivalent to } \, T T^* [I_{\cH_2} - P_{\ker(T^*)}]. 
\lb{A.10a}
\end{equation} 
In addition,
\begin{align}
& f\in \dom(H_1) \, \text{ and } \, H_1 f = \lambda^2 f, \; \lambda \neq 0,   \no \\
& \quad \text{implies }  \,  T f \in \dom(H_2) \, \text{ and } \, H_2(T f) = \lambda^2 (T f),    \lb{A.11} \\
& g\in \dom(H_2)\, \text{ and } \, H_2 \, g = \mu^2 g, \; \mu \neq 0,     \no \\
& \quad \text{implies }  \, T^* g \in \dom(H_1)\, \text{ and } \, H_1(T^* g) = \mu^2 (T^* g),    \lb{A.12} 
\end{align}
with multiplicities of eigenvalues preserved. \\
$(iii)$ One has for $z \in \rho(H_1) \cap \rho(H_2)$,
\begin{align}
& I_{\cH_2} + z (H_2 - z I_{\cH_2})^{-1} \supseteq T (H_1 - z I_{\cH_1})^{-1} T^*,    \lb{A.13} \\
& I_{\cH_1} + z (H_1 - z I_{\cH_1})^{-1} \supseteq T^* (H_2 - z I_{\cH_2})^{-1} T,    \lb{A.14}
\end{align}
and 
\begin{align}
& T^* \phi(H_2) \supseteq \phi(H_1) T^*, \quad 
T \phi(H_1) \supseteq \phi(H_2) T,   \lb {A.15} \\
& V_{T^*} \phi(H_2) \supseteq \phi(H_1) V_{T^*}, \quad 
V_T \phi(H_1) \supseteq \phi(H_2) V_T.   \lb {A.15a}
\end{align}
\end{theorem}
%%%%%%%%%%% 

As noted by E.\ Nelson (unpublished), Theorem \ref{tA.2} follows from the spectral theorem and the 
elementary identities, 
\begin{align}
& Q = V_Q |Q| = |Q| V_Q,    \lb{A.16} \\
& \ker(Q) = \ker(|Q|) = \ker (Q^2) = (\ran(Q))^{\bot} 
= \ker (T) \oplus \ker (T^*),    \lb{A.17}  \\
& I_{\cH_1 \oplus \cH_2} + z (Q^2 - z I_{\cH_1 \oplus \cH_2})^{-1} 
= Q^2 (Q^2 -z I_{\cH_1 \oplus \cH_2})^{-1} \supseteq Q (Q^2 -z I_{\cH_1 \oplus \cH_2})^{-1} Q,  \no \\
& \hspace*{9.3cm}   z \in \rho(Q^2),    \lb{A.18} \\
& Q \phi(Q^2) \supseteq \phi(Q^2) Q,   \lb{A.19}
\end{align}
where
\begin{equation}
V_Q = \begin{pmatrix} 0 & (V_T)^* \\ V_T & 0 \end{pmatrix} 
= \begin{pmatrix} 0 & V_{T^*} \\ V_T & 0 \end{pmatrix}.   \lb{A.20}
\end{equation}

In particular,
\begin{equation}
\ker(Q) = \ker(T) \oplus \ker(T^*), \quad  
P_{\ker(Q)} = \begin{pmatrix} P_{\ker(T)} & 0 \\ 0 & P_{\ker(T^*)} \end{pmatrix},    \lb{A.21}
\end{equation}
and we also recall that
\begin{equation}
\sigma_3 Q \sigma_3 = - Q, \quad \sigma_3 = \begin{pmatrix} I_{\cH_1} & 0 \\ 0 & - I_{\cH_2} 
\end{pmatrix},     \lb{A.22}
\end{equation}
that is, $Q$ and $-Q$ are unitarily equivalent. (For more details on Nelson's 
trick see also \cite[Sect.\ 8.4]{Te09}, \cite[Subsect.\ 5.2.3]{Th92}.) 
We also note that
\begin{equation}
\psi(|Q|) = \begin{pmatrix} \psi(|T|) & 0 \\ 0 & \psi(|T^*|) \end{pmatrix}    \lb{A.22a} 
\end{equation}
for Borel measurable functions $\psi$ on $\bbR$, and 
\begin{equation}
\ol{[Q |Q|^{-1}]} = \begin{pmatrix} 0 & (V_T)^*\\ V_T & 0 \end{pmatrix} = V_Q 
\, \text{ if } \, \ker(Q) = \{0\}.   \lb{A.22b} 
\end{equation}

Finally, we recall the following relationships between $Q$ and $H_j$, $j=1,2$.

%%%%%%%%%%%%
\begin{theorem} [\cite{BGGSS87}, \cite{Th88}] \lb{tA.3}
Assume Hypothesis \ref{hA.1}. \\
$(i)$ Introducing the unitary operator $U$ on $(\ker(Q))^{\bot}$ by
\begin{equation}
U = 2^{-1/2} \begin{pmatrix} I_{\cH_1} & (V_T)^* \\ -V_T & I_{\cH_2} \end{pmatrix} 
\, \text{ on } \,  (\ker(Q))^{\bot},     \lb{A.23}
\end{equation}
one infers that
\begin{equation}
U Q U^{-1} = \begin{pmatrix}  |A| & 0 \\ 0 & - |A^*| \end{pmatrix} 
\, \text{ on } \,  (\ker(Q))^{\bot}.     \lb{A.24}
\end{equation}
$(ii)$ One has
\begin{align}
\begin{split}
(Q - \zeta I_{\cH_1 \oplus \cH_2})^{-1} = \begin{pmatrix} \zeta (H_1 - \zeta^2 I_{\cH_1})^{-1} 
& T^* (H_2 - \zeta^2 I_{\cH_2})^{-1}  \\  T (H_1 - \zeta^2 I_{\cH_1})^{-1}  & 
\zeta (H_2 - \zeta^2 I_{\cH_2})^{-1}  \end{pmatrix},&    \\
\zeta^2 \in \rho(H_1) \cap \rho(H_2).&   \lb{A.25}
\end{split}
\end{align}
$(iii)$ In addition, 
\begin{align}
\begin{split} 
& \begin{pmatrix} f_1 \\ f_2 \end{pmatrix} \in \dom(Q) \, \text{ and } \, 
Q \begin{pmatrix} f_1 \\ f_2 \end{pmatrix} = \eta \begin{pmatrix} f_1 \\ f_2 \end{pmatrix}, \; \eta \neq 0,  
 \\
& \quad \text{ implies } \, f_j \in \dom (H_j) \, \text{ and } \, H_j f_j = \eta^2 f_j, \; j=1,2.    \lb{A.26}
\end{split} 
\end{align}
Conversely,
\begin{align}
\begin{split} 
& f \in \dom(H_1) \, \text{ and } H_1 f = \lambda^2 f, \; \lambda \neq 0, \\
& \quad \text{implies } \, \begin{pmatrix} f \\ \lambda^{-1} T f \end{pmatrix} \in \dom(Q) \, \text{ and } \, 
Q \begin{pmatrix} f \\ \lambda^{-1} T f \end{pmatrix} 
= \lambda \begin{pmatrix} f \\ \lambda^{-1} T f \end{pmatrix}.  
\end{split} 
\end{align}
Similarly,
\begin{align}
\begin{split} 
& g \in \dom(H_2) \, \text{ and } H_2 \, g = \mu^2 g, \; \mu \neq 0, \\
& \quad \text{implies } \, \begin{pmatrix} \mu^{-1} T^* g \\ g \end{pmatrix} \in \dom(Q) \, \text{ and } \, 
Q \begin{pmatrix} \mu^{-1} T^* g \\ g \end{pmatrix} 
= \mu \begin{pmatrix} \mu^{-1} T^* g \\ g \end{pmatrix}.  
\end{split} 
\end{align}
\end{theorem}
%%%%%%%%%%%%

%%%%%%%%%%%% Appendix B %%%%%%%%%%%%%%
\section{Adjoints and Closures of Operator Products} \lb{sB}
\renewcommand{\theequation}{B.\arabic{equation}}
\renewcommand{\thetheorem}{B.\arabic{theorem}}
\setcounter{theorem}{0} \setcounter{equation}{0}
%%%%%%%%%%%%%%%%%%%%%%%%%%%%%%%%%%%%%%
%%%%%%%%%%%%%%%%%%%%%%%%%%%%%%%%%%%%%%

The purpose of this appendix is to describe some situations in which equality 
holds between $(TS)^*$ and $S^* T^*$ and similarly, describe relations 
between $\ol{(TS)}$ and $\ol{T} \, \ol{S}$. 

We recall that if $C: \cH \supseteq \dom(C) \to \cH'$ is a closed operator (with 
$\cH$ and $\cH'$ complex, separable Hilbert spaces), then a linear subspace 
$\cD$ of $\dom(C)$ is called a {\it core} for $C$ if 
$\ol{C|_{\cD}} = C$.

%%%%%%%%%%%%%%%
\begin{lemma} \lb{lB.1} 
Let $\cH, \cH', \cH''$ be complex, separable Hilbert spaces, and introduce 
the linear operators 
\begin{equation}
S: \cH \supseteq \dom(S) \to \cH',  \quad 
T: \cH' \supseteq \dom(T) \to \cH''.    \lb{B.1} 
\end{equation} 
$(i)$ Assume that $T$ and $TS$ are densely defined. Then $S$ is densely 
defined and 
\begin{equation}
(TS)^* \supseteq S^* T^*.    \lb{B.2}
\end{equation}
$(ii)$ Suppose that $S$ is densely defined and $T \in \cB(\cH', \cH'')$. Then
\begin{equation}
(TS)^* = S^* T^*.    \lb{B.3}
\end{equation}
$(iii)$ Assume that $T$ and $TS$ are densely defined. In addition, 
suppose that $S$ is injective $($i.e., $\ker(S) = \{0\}$$)$ and 
$S^{-1} \in \cB(\cH',\cH)$. Then $S$ is densely defined and 
\begin{equation}
(TS)^* = S^* T^*.    \lb{B.4} 
\end{equation}
$(iv)$ Suppose that $TS$ is densely defined, assume that 
$\ol{\dom(T) \cap \ran(S)} = \cH'$, and introduce 
\begin{equation}
\hatt T = T|_{\dom(T) \cap \ran(S)}.   \lb{B.5}
\end{equation} 
Moreover, assume that $S$ is injective and that $S^{-1}$ is a bounded operator. 
Then $\ol{\dom(S)} = \cH$, $\ol{\ran(S)} = \cH'$, $\ol{S^{-1}} \in \cB(\cH',\cH)$, 
and 
\begin{equation}
(TS)^* = \big(\hatt T S\big)^* = S^* \big({\hatt T}\big)^* \supseteq S^* T^*.    
\lb{B.6}
\end{equation}
Suppose, in addition, that $T$ is closable. Then 
\begin{equation}
(TS)^* = S^* T^* \, \text{ if and only if $\dom(T) \cap \ran(S)$ is a 
core for $\ol T$.}     \lb{B.7}
\end{equation} 
$(v)$ Assume that $S$ and $T$ are densely defined, suppose $S$ 
is closed, and assume in addition that $\ran(S)$ has finite codimension 
$($i.e., $\dim\big(\ran(S)^\bot\big) < \infty$$)$. Then $TS$ is densely defined 
and 
 \begin{equation}
(TS)^* = S^* T^*.    \lb{B.7a}
\end{equation}
\end{lemma}
%%%%%%%%%%%%%%%
\begin{proof} We refer to \cite[Theorem 4.19\,$(a)$]{We80} for a proof of 
item $(i)$. 

Item $(ii)$ is a classical result, see, for instance, 
\cite[Lemma X.II.1.6]{DS88} and \cite[Theorem 4.19\,$(b)$]{We80}. 

Item $(iii)$ is mentioned in \cite[Exercise 4.18]{We80} and is a special case of 
item $(iv)$ to be proven next. 

To prove item $(iv)$ one can argue as follows. Since $S$ is injective, $\dom(S^{-1}) = \ran(S)$ is dense in $\cH'$, and 
$S^{-1}$ is a bounded operator, $S^{-1}$ is closable and hence 
$\dom\big(\ol{S^{-1}}\big) = \ol{\dom(S^{-1})} = \cH'$ (cf.\ 
\cite[Theorem\ 5.2]{We80}). Thus,  
$\ol{S^{-1}} \in \cB(\cH',\cH)$ by the closed graph theorem. 

The fact that $TS = \hatt T S$, $\hatt T \subseteq T$ (implying 
$T^* \subseteq \big(\hatt T\big)^*$), and generally, $S,\hatt T, \hatt T S$ all being densely defined implies that   
$\big(\hatt T S\big)^* \supseteq S^* \big(\hatt T\big)^*$ (cf.\ 
item $(i)$), \eqref{B.6} will follow once one proves that 
$\big(\hatt T S\big)^* \subseteq S^* \big(\hatt T\big)^*$. For this purpose let 
$f \in \dom\big(\big(\hatt T S\big)^*\big)$ and $g \in  \dom\big(\hatt T S\big)$, 
then
\begin{align}
\big(f,\hatt T S g\big)_{\cH''} &= \big(\big(\hatt T S\big)^* f, g\big)_{\cH} 
= \big(\big(\hatt T S\big)^* f, S^{-1} S g\big)_{\cH} 
= \big(\big(\hatt T S\big)^* f, \ol{S^{-1}} S g\big)_{\cH}   \no \\
&= \big(\big(\ol{S^{-1}}\big)^* \big(\hatt T S\big)^* f, S g\big)_{\cH'}.  
\lb{B.8}
\end{align}
Thus,
\begin{equation}
\big|\big(f, \hatt T S g\big)_{\cH''}\big| \leq 
\big\|\big(\ol{S^{-1}}\big)^* \big(\hatt T S\big)^* f\big\|_{\cH'} \, 
\|S g\|_{\cH'}.   \lb{B.9}
\end{equation}
Since $\dom\big(\hatt T S\big) = S^{-1} \dom\big(\hatt T\big)$, one obtains that 
as $g$ runs through all of $\dom\big(\hatt T S\big)$, $Sg$ runs through all of 
$ \dom\big(\hatt T\big)$. Hence, \eqref{B.9} implies that 
$f \in \dom\big(\big(\hatt T\big)^*\big)$ and thus \eqref{B.8} yields 
\begin{equation}
\big(\big(\hatt T\big)^* f, S g\big)_{\cH'} 
= \big(\big(\ol{S^{-1}}\big)^* \big(\hatt T S\big)^* f, S g\big)_{\cH'}.  
\lb{B.10}
\end{equation}
Since (as a consequence of the hypothesis $\ol{\dom(T) \cap \ran(S)} = \cH'$), $\ol{\ran(S)} = \cH'$, \eqref{B.10} implies 
\begin{align}
\big(\hatt T\big)^* f &= \big(\ol{S^{-1}}\big)^* \big(\hatt T S\big)^* f 
= \big(S^{-1}\big)^* \big(\hatt T S\big)^* f    
= (S^*)^{-1} \big(\hatt T S\big)^* f.    \lb{B.11}
\end{align}
Here we used that $\big(\ol{A}\big)^* = A^*$ if $A$ is densely defined and closable 
(cf.\ \cite[Theorem\ 5.3\,(c)]{We80}) and that $\big(B^{-1}\big)^* = (B^*)^{-1}$ 
if $B$ is injective and densely defined with dense range (implying injectivity 
of $B^*$, cf.\ \cite[Theorem\ 4.17\,(b)]{We80}). Equation \eqref{B.11} yields 
$S^* \big(\hatt T\big)^* f = \big(\hatt T S\big)^* f$ and hence 
$\big(\hatt T S\big)^* \subseteq S^* \big(\hatt T\big)^*$.

Next, assume in addition that $T$ is closable (and hence $\big(\ol T\big)^* = T^*$). Then if $\dom(T) \cap \ran(S)$ is a core for $\ol T$, $\ol{\hatt T} = \ol T$ yields
\begin{equation}
\big(\hatt T\big)^* = \big(\ol{\hatt T}\big)^* = \big(\ol T\big)^* = T^*,   \lb{B.12}
\end{equation}
and hence \eqref{B.6} implies $(TS)^* = S^* T^*$. Conversely, suppose that 
$\big(\hatt T\big)^* = T^*$. Then
\begin{equation}
\ol{\hatt T} = \big(\big(\hatt T\big)^*\big)^* = (T^*)^* = \ol T     \lb{B.13}
\end{equation}
proves that $\dom(T) \cap \ran(S)$ is a core for $\ol T$. \\
For a proof of item $(v)$ we refer to \cite{Gu69}, \cite{Ho68}, \cite{Ho69}, 
\cite{Sc70}, and \cite{vCG70}. In this context we note that $\ran(S)$ is closed 
in $\cH'$ (since $S$ is assumed to be closed and 
$\dim\big(\ran(S)^\bot\big) < \infty$, cf.\ \cite[Corollary IV.1.13]{Go85b}) and 
hence does not have to be assumed to be closed, and similarly, it is not necessary to assume that $T$ is closed as is done in some references.  
\end{proof}
%%%%%%%%%%%%%%% 

We note again that Lemma \ref{B.1}\,$(iv)$ is a refinement of 
\cite[Exercise\ 4.18]{We80}, listed as item $(iii)$ in Lemma \ref{B.1}; it may 
be of independent interest. 

For additional results guaranteeing $(TS)^* = S^* T^*$ (including the Banach 
space setting), we refer, for instance, to \cite{Gu69}, \cite{Ho68}, \cite{Ho69}, 
\cite{KS63}, \cite{Pf81}, and \cite{vCG70} (in particular, the case of nondensely defined operators is discussed in detail in \cite{Pf81}). 

Next, we briefly consider situations which relate $\ol{ST}$ with $\ol{S} \, \ol{T}$ 
(much less appears to have been studied in this context). 

%%%%%%%%%%%%%%%
\begin{lemma} \lb{lB.2} 
Let $\cH, \cH', \cH''$ be complex, separable Hilbert spaces, and introduce 
the linear operators 
\begin{equation}
S: \cH \supseteq \dom(S) \to \cH',  \quad 
T: \cH' \supseteq \dom(T) \to \cH''.    \lb{B.15} 
\end{equation} 
$(i)$ Assume that $S$ is bounded, $\ol S \in \cB(\cH,\cH')$, and that $T$ 
is closed. Then $T{\ol S}$ is closed, implying that $TS$ is closable and that 
\begin{equation}
\ol{TS} \subseteq T{\ol S}.     \lb{B.16}
\end{equation}
$(ii)$ Assume that $S$ is injective with $S^{-1}$ bounded and 
$\ol{S^{-1}} \in \cB(\cH',\cH)$.
Furthermore, suppose $T|_{\dom(T) \cap \ran(S)}$ is closable and
\begin{equation}
T|_{\dom(T) \cap \ran(\ol S)} \subseteq \ol{T|_{\dom(T) \cap \ran(S)}}. \lb{B.17}
\end{equation}
In addition, assume that $TS$ is closable. Then
\begin{equation}
T{\ol S} \subseteq \ol{TS}.       \lb{B.18}
\end{equation} 
\end{lemma}
%%%%%%%%%%%%%%%
\begin{proof}
For the purpose of proving item $(i)$ we suppose that 
$\{f_n\}_{n\in\bbN} \subset \dom(T {\ol S})$ 
such that $\slim_{n\to\infty} f_n = f \in \cH$ and $T {\ol S} f_n = h \in \cH''$. By 
the definition of $\dom(T{\ol S})$, this implies that 
$\{f_n\}_{n\in\bbN} \subset \dom({\ol S})$, and since $\ol S \in \cB(\cH,\cH')$, 
one concludes that $\slim_{n\in\bbN} {\ol S} f_n = {\ol S} f \in \cH'$. Since 
$\slim_{n\to\infty}T ({\ol S} f_n) = \slim_{n\to\infty}T {\ol S} f_n = h$, closedness 
of $T$ implies that ${\ol S} f \in \dom(T)$ and 
$\slim_{n\to\infty} T({\ol S} f_n) = T ({\ol S} f)$, that is, 
$f \in \dom(T {\ol S})$ and $\slim_{n\to\infty} T {\ol S} f_n = T {\ol S} f$. Thus, 
$T{\ol S}$ is closed. 

Since $T S \subseteq T {\ol S}$ and the latter is closed, $TS$ is closable and 
\begin{equation}
\ol{T S} \subseteq \ol{T {\ol S}} = T {\ol S}.   \lb{B.19}
\end{equation}
To prove item $(ii)$ let $f\in \dom(T \ol S)$ and $g=T \ol S f$.
Then $h= \ol S f \in \dom(T) \cap \ran(\ol S)$, and by assumption \eqref{B.17}
we can find $\{h_n\}_{n\in\bbN} \in \dom(T) \cap \ran(S)$ such that 
$\slim_{n\to\infty} h_n = h$ in $\cH'$ and $\slim_{n\to\infty} T h_n = T h = g$ in 
$\cH''$. Since $S^{-1}$ is bounded, the sequence $f_n = S^{-1} h_n$ converges 
strongly to $\ol{S^{-1}} h = \big(\ol S\big)^{-1} h = f$ in $\cH$, and by construction, 
$T S f_n = T h_n$, 
$n\in\bbN$, satisfies $\slim_{n\to\infty} T S f_n = g$. Thus, $f \in \dom(\ol{T S})$
and $\ol{T S} f = g = T \ol S f$.
\end{proof}
%%%%%%%%%%%%%%%

We note that closedness of $T{\ol S}$ in Lemma \ref{lB.2}\,$(i)$ has been noted 
in \cite[Proposition\ B.2]{EN00}.

\medskip

%%%%%%%%%%%%%%%%%%%%%%%%%%%%%%%%%%%%%
\noindent {\bf Acknowledgments.}
We are indebted to Klaus Engel and Delio Mugnolo for providing us with a number of pertinent references.  
%and to Viktor Losert for a stimulating discussion concerning 
%Rajchman measures.

This paper was initiated when taking part in Êthe international research program on Nonlinear Partial Differential Equations at the Centre for Advanced Study (CAS) at the Norwegian Academy of Science and Letters in Oslo during the academic year 2008Ð-09. F.G.\ and G.T.\ gratefully acknowledge  the great hospitality at CAS during a five, respectively, four-week stay in May--June, 2009.
%%%%%%%%%%%%%%%%%%%%%%%%%%%%%%%%%%%%%

%%%%%%%%%%%%%%%%%%%%%%%%%%%%%%%%%%%%%%
%%%%%%%%%%%%%%%%%%%%%%%%%%%%%%%%%%%%%%

\end{document}